\documentclass[11pt,UTF-8,reqno]{amsart}
\usepackage{enumerate}
\usepackage{mhequ}
\usepackage[margin=1.0in]{geometry}
\usepackage{amssymb,url,color, booktabs,nccmath}
\usepackage{mathrsfs}
\usepackage{enumitem}
\usepackage{graphicx}
\usepackage{mathtools}
\usepackage{tikz}
\usepackage{graphicx}
\usetikzlibrary{decorations.pathreplacing,decorations.markings}
\usepackage{array}      
\usepackage{colortbl}   
\usepackage{xcolor}
\usepackage{makecell}

\usepackage{color}
\usepackage[colorlinks=true]{hyperref}
\hypersetup{
    linkcolor=blue,          
    citecolor=red,        
    filecolor=blue,      
    urlcolor=cyan
}

\definecolor{lightgray}{gray}{0.9}
\definecolor{faintgray}{gray}{0.7}

\definecolor{darkergreen}{rgb}{0.0, 0.5, 0.0}

\numberwithin{equation}{section}

\newcommand{\be}{\begin{eqnarray}}
\newcommand{\ee}{\end{eqnarray}}
\newcommand{\ce}{\begin{eqnarray*}}
\newcommand{\de}{\end{eqnarray*}}
\newtheorem{theorem}{Theorem}[section]
\newtheorem{meta-theorem}{Meta-Theorem}[section]
\newtheorem{lemma}[theorem]{Lemma}
\newtheorem{remark}[theorem]{Remark}
\newtheorem{definition}[theorem]{Definition}
\newtheorem{proposition}[theorem]{Proposition}
\newtheorem{Examples}[theorem]{Example}
\newtheorem{corollary}[theorem]{Corollary}

\def\proj{\mathbf{p}}

\def\Re{{\mathrm{Re}}}

\def\var{{\mathrm{var}}}

\def\eps{\varepsilon}

\def\p{\partial}

\def\[{{\Big[}}
\def\]{{\Big]}}
\def\<{{\langle}}
\def\>{{\rangle}}
\def\({{\Big(}}
\def\){{\Big)}}

\def\bx{{\mathbf{x}}}
\def\tr{\mathrm {tr}}

\def\Ric{{\rm Ricci}}

\def\dif{{\mathord{{\rm d}}}}

\def\Hess{{\mathord{{\rm Hess}}}}
\def\min{{\mathord{{\rm min}}}}

\def\no{\nonumber}
\def\={&\!\!=\!\!&}
\def\YMH{\textnormal{\small \textsc{ymh}}}
\def\YM{\textnormal{\small \textsc{ym}}}

 \newcommand{\eqdef}{\stackrel{\mbox{\tiny def}}{=}}

\def\cE{{\mathcal E}}

\def\cL{{\mathcal L}}

\def\cP{{\mathcal P}}
\def\cQ{{\mathcal Q}}

\def\cS{{\mathcal S}}

\def\mN{{\mathbb N}}

\def\mR{{\mathbb R}}
\def\mS{{\mathbb S}}
\def\mT{{\mathbb T}}

\def\mZ{{\mathbb Z}}

\def\1{{\mathbf{1}}}

\def\sG{{\mathscr G}}

\def\E{\mathbf E}

\def\geq{\geqslant}
\def\leq{\leqslant}
\def\ge{\geqslant}
\def\le{\leqslant}
\def\dg{d(\mfg)}

\def\E{\mathbf{E}}
\def\CP{\mathcal{P}}
\def\Tr{\mathrm{Tr}}
\def\q{\mathfrak{q}}

\def\div{\mathord{{\rm div}}}

\def\bx{{\mathbf{x}}}
\def\tr{\mathrm {tr}}

\def\Ric{{\rm Ricci}}

\def\dif{{\mathord{{\rm d}}}}

\def\Hess{{\mathord{{\rm Hess}}}}
\def\min{{\mathord{{\rm min}}}}

\def\tr{{\rm Tr}}
\def\no{\nonumber}
\def\={&\!\!=\!\!&}
\def\bt{\begin{theorem}}
\def\et{\end{theorem}}
\def\bl{\begin{lemma}}
\def\el{\end{lemma}}
\def\br{\begin{remark}}
\def\er{\end{remark}}
\def\bx{\begin{Examples}}
\def\ex{\end{Examples}}
\def\bd{\begin{definition}}
\def\ed{\end{definition}}
\def\bp{\begin{proposition}}
\def\ep{\end{proposition}}
\def\bc{\begin{corollary}}
\def\ec{\end{corollary}}

\def\so{\mathfrak{so}}

\def\mfg{\mathfrak{g}}

\def\geq{\geqslant}
\def\leq{\leqslant}
\def\ge{\geqslant}
\def\le{\leqslant}

\def\div{\mathord{{\rm div}}}

\def\R{\mathbb R}

\def\<{\langle} \def\>{\rangle}
\def\x{{\bf x}}
\def\y{{\bf y}}

\def\cov{\textrm{cov}}
\def\${|\!|\!|}

\allowdisplaybreaks

\makeatletter
\def\section{\@startsection{section}{1}%
  \z@{1.7\linespacing\@plus\linespacing}{.5\linespacing}%
  {\normalfont\scshape\centering}}
\def\subsection{\@startsection{subsection}{2}%
  \z@{1\linespacing\@plus.7\linespacing}{-.5em}%
  {\normalfont\bfseries}}
\makeatother

\tikzset{
dot/.style={circle,fill,inner sep=1pt},
  on each segment/.style={
    decorate,
    decoration={
      show path construction,
      moveto code={},
      lineto code={
        \path [#1]
        (\tikzinputsegmentfirst) -- (\tikzinputsegmentlast);
      },
      curveto code={
        \path [#1] (\tikzinputsegmentfirst)
        .. controls
        (\tikzinputsegmentsupporta) and (\tikzinputsegmentsupportb)
        ..
        (\tikzinputsegmentlast);
      },
      closepath code={
        \path [#1]
        (\tikzinputsegmentfirst) -- (\tikzinputsegmentlast);
      },
    },
  },
  mid arrow/.style={postaction={decorate,decoration={
        markings,
        mark=at position .5 with {\arrow[#1]{stealth}}
      }}},
midarrow/.style={postaction={decorate,decoration={markings,mark=at position 0.5 with {\arrow{>}}}}},
}

\begin{document}

%
%
%
%


\date{\today}

\title{Langevin dynamics of lattice Yang--Mills--Higgs and applications}

\author{Hao Shen}
\address[H. Shen]{Department of Mathematics, University of Wisconsin - Madison, USA}
\email{pkushenhao@gmail.com}

\author{Rongchan Zhu}
\address[R. Zhu]{Department of Mathematics, Beijing Institute of Technology, Beijing 100081, China 
	}
\email{zhurongchan@126.com}

\author{Xiangchan Zhu}
\address[X. Zhu]{State Key Laboratory of Mathematical Sciences, Academy of Mathematics and Systems Science,
	Chinese Academy of Sciences, Beijing 100190, China
	}
\email{zhuxiangchan@126.com}

\maketitle

\begin{abstract}
\normalsize
In this paper,
we investigate the Langevin dynamics of various lattice formulations of the Yang--Mills--Higgs model,  with an inverse Yang--Mills coupling $\beta$ and a Higgs parameter $\kappa$.
The Higgs component is either a bounded field taking values in a compact target space, or an unbounded field taking values in a vector space in which case the model also has a Higgs mass parameter $m$.
We study the regime where  $(\beta,\kappa)$ are small  in the first case
or  $(\beta,\kappa/m)$ are small  in the second case.
We prove the exponential ergodicity of the dynamics on the whole lattice via functional inequalities. 
We establish exponential decay of
 correlations for a broad class of observables, namely, the infinite volume measure exhibits a strictly positive mass gap.
 Moreover, when the target space of the Higgs field is compact,
 appropriately rescaled observables exhibit factorized correlations in the large $N$ limit.
These extend the earlier results \cite{SZZ22} on pure lattice Yang--Mills
to the case with a coupled Higgs field.

Unlike pure lattice Yang--Mills where the field is always bounded,
in the case where the coupled Higgs component is unbounded,
the control of its behavior is much harder and requires new techniques.
 Our approach involves a disintegration argument and a delicate analysis of correlations to effectively control the unbounded Higgs component.
\end{abstract}

\setcounter{tocdepth}{2}
\tableofcontents

\section{Introduction}
\label{sec:intro}

Given  a compact Lie group $G$ with Lie algebra $\mfg$, the Euclidean quantum Yang--Mills--Higgs (YMH) model on $\mR^d$ with structure group $G$ is formally defined by a ``probability measure''
\begin{equ}\label{eq:YMH_measure}
	\dif\mu_\YMH(A,\Phi)  \propto
	\exp\big(
	S_{\YMH}(A,\Phi)
	\big)\, \dif A\, \dif \Phi.
\end{equ}
Here  the ``connection'' field (or Yang--Mills field) $A$ is a $\mfg$-valued 1-form,
 the Higgs field $\Phi$ takes values in a vector space endowed with an inner product, on which $G$ has an orthogonal representation,
and $\dif A \dif\Phi$ is a formal Lebesgue measure.
In a more general geometric setting,
$A$ is a connection on a given principal bundle on $\mR^d$,
the representation  of $G$ determines an associated vector bundle,
and $\Phi$ is a section of the associated  bundle;
in this paper, we assume that the principal and the associated bundles are trivial.
The continuum YMH action $S_{\YMH}$ is given by
\footnote{In physicists's standard convention one put the minus sign in \eqref{eq:YMH_measure} instead of \eqref{eq:YM_energy}.
Here we prefer to write the measure in the form $\exp(S_{\YMH})$ for convenience.}
\begin{equ}\label{eq:YM_energy}
	S_{\YMH}(A,\Phi) \eqdef  - \int_{\mR^{d}}
	\Big(|F_A(x)|^2 +  |\dif_A \Phi (x)|^2 
	+V(\Phi(x))
	\Big)\, \dif x\;.
\end{equ}
Here,
 $F_A$
is the curvature 2-form of $A$.
Also, $\dif_A$ denotes the covariant derivative on the Higgs fields defined by the connection $A$,
 namely, $\dif_A\Phi = \dif\Phi+A\Phi$ where $\dif$ is the exterior derivative and $A\Phi$ is understood
 via the induced representation of $\mfg$ on the target space of $\Phi$.
The potential
$V(\Phi(x))$ is polynomial in $\Phi$ which is invariant under rotation of the space in which $\Phi$ takes values,
for instance $V(\Phi(x))= m|\Phi (x)|^2 + \frac12 |\Phi(x)|^4$.
When $G$ is abelian, this is also called the Coleman--Weinberg model.
The model without the term $ |\dif_A \Phi (x)|^2$ and with $V=0$ will be called pure Yang--Mills (YM) model.

An important feature of the model is gauge invariance.
Assuming that all the objects are smooth,
a \textit{gauge transformation} is a $G$-valued function on $\R^d$  which acts on
$(A,\Phi)$ by
\begin{equ}[e:gauge-transformation]
	g \circ A  \eqdef gAg^{-1}  - (\dif g)g^{-1} \;,
	\qquad
	\mbox{and}
	\qquad
	g\circ \Phi  \eqdef  g\Phi\;,
\end{equ}
and one can verify $S_{\YMH}(g \circ  A, g \circ  \Phi) = S_{\YMH}(A,\Phi)$.
Due to gauge invariance of the action,  observables which are gauge invariant  are of particular interest. Two well-known classes of  gauge invariant observables are called Wilson loops and Wilson lines: again assuming that everything is smooth,
the Wilson loop observable
given a loop in $ \mR^d$ is given by
a class function
(which means  a function on a group $G$ that is constant on the conjugacy classes of $G$, for example the trace) of the holonomy of $A$ along the loop,
and the Wilson line observable
given a path in $ \mR^d$ is given by the holonomy of $A$ along the path
contracted with the values of $\Phi$ at the beginning and the end points of the path.

There are multiple reasons which make the YMH model \eqref{eq:YMH_measure} prominent and interesting.
From the physics point of view,
the fundamental Standard Model for elementary particles
is described by an action functional with three fields: Yang--Mills, fermions, and Higgs,
and our model \eqref{eq:YM_energy} can be viewed as a simplified version without the fermion field.
%
Also, from the geometric perspective, since a given connection $A$
specifies a way of differentiating sections of a bundle via the covariant derivative $\dif_A$, the term $|\dif_A \Phi|^2$ is the natural minimal coupling between the two fields, i.e. the connection field $A$ and the bundle section $\Phi$.
Turning to the probabilistic perspective, the pure YM model in 2D exhibits integrability,
a crucial property that led to a fruitful list of rigorous results (c.f. ~\cite{Driver89,GKS89,Sengupta97,Levy03}, and more recently \cite{MR3861073,Chevyrev19YM,PPSY2023}),
but such integrability is broken once the model is perturbed by a Higgs field,
even in 2D, and thus more tools are required to study the Yang--Mills-Higgs model.
We will discuss more related literature in Section~\ref{sec:Discussion}.

In this paper we
study several lattice versions of the YMH model.
Let  $\Lambda \subset \mZ^d$ be a finite lattice.
Postponing the precise definitions to Section~\ref{sec:Lattice YMH},
a lattice YMH model consists of two fields (see Fig.~\ref{fig:YMH}):
\begin{enumerate}
\item
The lattice Yang--Mills field $Q$, as the discrete analogue of $A$ in \eqref{eq:YM_energy}.
The field $Q$ is a $G$-valued field on  edges.
In other words, it assigns each lattice edge $e=(x,y)$ an element $Q_e\in G$.
Moreover, we impose $Q_{e^{-1}} = Q_e^{-1}$ where $e^{-1}$ is the reverse edge of $e$.
\item
The lattice Higgs field $\Phi$, as the  discrete analogue of the field $\Phi$ in \eqref{eq:YM_energy}.
The field $\Phi$ is an $M$-valued field on  vertices,
where $M$ is a fixed finite dimensional Riemannian manifold.
In other words,  $\Phi$ assigns each lattice vertex $x$ an element $\Phi_x\in M$.
\end{enumerate}
We also assume that we have a $G$-action on $M$, which preserves the metric on $M$.
In this paper,
we will restrict to  
 the case where $G=SO(N)$ for arbitrary $N>2$,
 but
we will consider several choices for the target space of the Higgs field. 
In particular, we will consider
\begin{equ}[e:M]
M\in \{\mR^N, \mS^{N-1}, G\},
\end{equ}
with $\mS^{N-1}$ being $(N-1)$-dimensional sphere with radius $1$.
For simplicity, we will also focus on
a quadratic potential $V(\Phi)$.

\begin{figure}[h]
\begin{tikzpicture}[scale=1.5]
\foreach \y in {0,1,2,3} {
    \draw (0,\y) -- (5,\y);
}
\foreach \x in {0,1,2,3,4,5} {
    \draw (\x,0) -- (\x,3);
}
\node at (1.5,1.5) {$Q_p$};
\draw [>=stealth,midarrow] (1,1) to (2,1); \node at (1.5,0.8) {$Q_{e_1}$};
\draw [>=stealth,midarrow] (2,1) to (2,2); \node at (2.25,1.5) {$Q_{e_2}$};
\draw [>=stealth,midarrow] (2,2) to (1,2); \node at (1.5,2.2) {$Q_{e_3}$};
\draw [>=stealth,midarrow] (1,2) to (1,1); \node at (0.8,1.5) {$Q_{e_4}$};
\draw [>=stealth,midarrow] (4,1) to (4,2); \node at (4.2,1.5) {$Q_{e}$};
   \fill (4,1) circle (1pt); \node at (3.8,0.8) {$\Phi_{x}$};
   \fill (4,2) circle (1pt);  \node at (3.8,2.2) {$\Phi_{y}$};
\end{tikzpicture}
\caption{The edge variables $Q$, the vertex variables $\Phi$,
the plaquette variable
$Q_p=\prod_{i=1}^4Q_{e_i}$, and the discrete covariant derivative $Q_e \Phi_y - \Phi_x$ for $e=(x,y)$.}
\label{fig:YMH}
\end{figure}

Given  coupling constants
 $(\beta, \kappa, m)$,
the model is
 given by well-defined probability measures  of the form
\begin{equ}[e:LYMH]
	\dif\mu_{\Lambda}(Q,\Phi)
	:= Z_{\Lambda}^{-1}
	\exp\Big(\mathcal S_{\YMH} (Q,\Phi)\Big)
	\prod_{e}  
	\dif\sigma(Q_e) 
	\prod_{z} \dif\Phi_z\,
\end{equ}
where $\dif\sigma$ is the Haar measure on $G$,
$\dif\Phi_z$ is a suitable measure on $M$  invariant under the $G$ action,
and $Z_{\Lambda}$ is a normalization constant.
The action
$\cS_{\YMH}$ is a suitable discrete version of \eqref{eq:YM_energy} of the following form
\begin{equ}[e:YMH0]
	\mathcal S_{\YMH}(Q,\Phi) =
	N\beta
	 \sum_{p} 
	\Tr(Q_p)
	- N\kappa
	\sum_{e} 
	|Q_e \Phi_y - \Phi_x|^2
	-N m
	 \sum_{z}|\Phi_z|^2 \;. 
\end{equ}
Here, $z,e,p$ are vertices, edges, and plaquettes, and
their  precise meanings (especially the positive orientations of edges and plaquettes)
will be given in Section~\ref{sec:Lattice YMH}.
In the first term,
for a plaquette $p$ consisting of four consecutive edges $e_1,\cdots,e_4$,
we write $Q_p=\prod_{i=1}^4Q_{e_i}$.
In the second term, $e=(x,y)$, and $Q_e \Phi_y - \Phi_x$ is a discrete covariant derivative,
in which $Q_e \Phi_y$ is well-defined thanks to the action of $G$ on $M$.
If the last term is a constant (for example $|\Phi_z|^2=1$
when  $M=\mS^{N-1}$), then we can drop it by a change of normalization.

Remark that  we restrict to $G=SO(N)$ for simplicity of exposition.
We focus on the list \eqref{e:M} of target spaces $M$
since they are prototypical choices that exemplify the different challenges and strategies.
Our techniques
can be  extended to other choices of  Lie groups such as $SU(N)$, $Sp(N)$,
and other target spaces $M$ (see Remark~\ref{rem:1} for some discussions).

One of the main objects of study is
the Langevin dynamics, or stochastic quantization  of the measure \eqref{e:LYMH}.
The goals of the paper include the study of

(1) ergodicity properties of the dynamics;

(2) functional inequalities, i.e. the log-Sobolev and the Poincar\'e inequalities of the measures;

(3) mass gap, namely exponential decay of correlations.

These are, as we will show, intertwined questions.
For pure Yang--Mills on the lattice (i.e. $\kappa=m=0$), such a dynamical approach was developed in \cite{SZZ22}.

In general, functional inequalities are powerful tools to study spatial or temporal mixing properties. Recall that  a measure satisfying the log-Sobolev inequality means that
$$
\E (F^2 \log F^2) \le C\, \E |\nabla F|_{L^2}^2
$$
for $C>0$ and smooth cylinder functions $F$ with $\E (F^2)=1$,
and it satisfying the Poincar\'e inequality means that
\begin{equ}[e:basic-Poin]
\E (F^2) \le C\, \E |\nabla F|_{L^2}^2
\end{equ}
for $C>0$ and smooth cylinder functions $F$ with $\E (F)=0$.
(Left-hand side is the variance if $\E (F)\neq 0$.)
The meaning of the Dirichlet form $\E |\nabla F|_{L^2}^2$ (especially the gradient $\nabla$) depends on the particular model,
see Section~\ref{sec:YMHSDEs} for the precise discussion.
The advantage of having the Poincar\'e inequality is that it immediately implies the $L^2$
exponential ergodicity, i.e.  $\E (P_t F)^2\leq e^{-2 t/C}\E(F^2)$ where $P_t$ is the semigroup associated to the Markov process generated by the Dirichlet form,
which easily follows from differentiating in $t$ and integration by parts in  the Dirichlet form (c.f. \cite[Theorem 1.1.1]{Wang});
this then connects the inequality with temporal mixing.
For general readers, we recall that in the simplest setting where one has
a measure $\frac{1}Z e^{S(x)}\dif x$ on an Euclidean space,
\begin{equ}[e:strategy]
\mbox{$S$ is strictly concave (Bakry--\'Emery condition)} \Rightarrow \mbox{log-Sobolev}
\Rightarrow \mbox{Poincar\'e}\;.
\end{equ}
On a compact Riemannian manifold, there is no continuous strictly concave function
\cite{Yau74}, but the geometric version of the Bakry--\'Emery condition
states that if the Ricci curvature of the underlying manifold  minus the Hessian of $S$
is strictly positive, one still has the log-Sobolev inequality:
\begin{equ}[e:strategyGeo]
\mbox{Strictly positive lower bound of  $\Ric-\Hess_{S}$} \; \Rightarrow \mbox{log-Sobolev}
\Rightarrow \mbox{Poincar\'e}\;.
\end{equ}
This is the strategy in \cite{SZZ22}, where the manifold is a product Lie group with positive Ricci curvature, which dominates the Hessian of the pure Yang--Mills action in the strong coupling regime, namely, under the small $\beta$ condition \cite[Assumption~1.1]{SZZ22}
\begin{equ}[e:K_YM]
K_{\YM}\eqdef \frac14(N-2)-8(d-1)|\beta|N  >0. \tag{B\'E-YM}
\end{equ}
($K_{\YM}$ was denoted as $K_{\cS}$ therein.)
As some geometric background of \eqref{e:strategyGeo}:  a
Riemannian manifold equipped with the measure $\frac{1}Ze^{S}\dif \sigma$
is called a weighted manifold, where $\dif \sigma$ is the Riemannian volume form;
  the generator $\cL$ for the related Langevin dynamics is given by $\Delta+\langle \nabla {S},\nabla \rangle$. In this geometric setting
 the Bakry--\'Emery condition is just the strictly positive lower bound of $\Ric(\cL):=\Ric-\Hess_{S}$, which is also called the ``curvature for the weighted manifold'' (see \cite[Section 1.16]{MR3155209}  for further discussion).

In the setting of this paper, we will encounter situations for which the simple strategy \eqref{e:strategy} or \eqref{e:strategyGeo} does not go through, but we still manage to show the desired mixing results. Below we state our results and discuss such challenges.
We will sometimes write B\'E  for Bakry--\'Emery,
and LSI for the log-Sobolev inequality.

Regarding the mass gap problem, the general message is that temporal mixing (ergodicity) and spatial mixing (mass gap) are often related. This is how the dynamic, or stochastic quantization, comes into play and serves as a powerful tool to demonstrate mass gap.
Continuing with the general discussion above, assume the Poincar\'e inequality \eqref{e:basic-Poin} and let $P_t$ be the Markov semigroup of a dynamic that leaves the measure invariant,
and consider observables  $f,g$ depending on the fields around two separated space points $x,y$.
 Then 
\begin{equ}[e:gen-mass]
\mbox{Cov}(f,g)
 = \E (fg)-\E(f)\E(g)= \E (P_t (fg))-\E (P_t f)\E( P_t g)
  =  \mbox{Cov}(P_t f,P_t g)+\mbox{Comm}
\end{equ}
where ``Comm'' arises since $P_t (fg) \neq P_t f P_t g$.
By Cauchy--Schwarz and the $L^2$ exponential ergodicity mentioned above, one has
$ \mbox{Cov}(P_t f,P_t g) 
\lesssim e^{-|x-y|/C}$,
if we choose $t\sim |x-y|$. 
This is the strategy for pure Yang--Mills in \cite{SZZ22},
where the main work is to show that the commutator term ``Comm'' also satisfies the same bound (as long as $t \lesssim |x-y|$),
and this requires detailed analysis for the generator of the Yang--Mills semigroup $P_t$.
Such argument goes back to \cite[Section 8.3]{GZ} for bounded continuous spin models,
where the analysis for the generator of course depends on the particular model.

We will explain more on the challenges in implementing the idea outlined in \eqref{e:gen-mass} to deduce mass gap from temporal mixing in our particular setting,
but, let us mention here that the converse is also often true:  one can often use spatial mixing to prove functional inequalities (c.f. \cite{Zegarlinski1996, Yosida} and  references therein).
It turns out that we will employ this strategy below as well: 
in the case where the Higgs field $\Phi$ takes values in an unbounded space $\R^N$
(for which we will see that a direct application of B\'E does not work),
we will first use the general strategies \eqref{e:strategy} and \eqref{e:gen-mass} to prove  
 functional inequalities and spatial mixing for the Higgs marginal (Lemma~\ref{lem:co}),
and then invoke this spatial mixing for the Higgs marginal to prove 
the Poincar\'e  inequality for the whole Yang--Mills--Higgs measure,
and finally use this Poincar\'e  inequality to establish a mass gap for the  YMH model.


\subsection{Main results}

The  Langevin dynamics for the Yang--Mills--Higgs model
are stochastic differential equations (SDEs)
with the following general form:
\begin{equ}[e:LSDE]
	\dif (Q,\Phi) = \nabla \mathcal S_\YMH (Q,\Phi) \dif t + \sqrt 2\dif \mathfrak B\;,
\end{equ}
where $\nabla$ is the gradient and $\mathfrak B$ is  the  Brownian motion
on the configuration space of $(Q,\Phi)$ (see Section~\ref{sec:mani}).
For various choices of the Higgs target space $M$,
we derive more explicit forms for equation \eqref{e:LSDE}
in Section~\ref{sec:YMHSDEs}: see \eqref{SDE:RN}, \eqref{SDES}, \eqref{SDEG}.
We will first show that
\eqref{e:LSDE} is globally well-posed, and has \eqref{e:LYMH} as invariant measure.
Namely: 

\begin{theorem}[Well-posedness and invariant measure]
\label{main:1}

On a finite lattice $\Lambda\subset \mZ^d$,
let $M$ be one of the target spaces in \eqref{e:M}, namely
 $M\in \{\mR^N, \mS^{N-1}, G\}$.
Then, \eqref{e:LSDE} is globally well-posed.

Again for these choices of $M$ and on a finite lattice $\Lambda$,
assuming further $m,\kappa>0$ in the case  $M=\mR^N$, the measure $\mu_\Lambda$ in \eqref{e:LYMH} is invariant under the SDE system \eqref{e:LSDE}.

Moreover, global well-posedness of  \eqref{e:LSDE}  can be extended to the entire lattice $\mZ^d$.

Finally, the invariant measures $(\mu_\Lambda)_{\Lambda\subset \mZ^d}$  \eqref{e:LYMH} form a tight family as the finite lattice $\Lambda$ grows to $\mZ^d$,
and any tightness limit is an invariant measure of the dynamics \eqref{e:LSDE} on the entire lattice $\mZ^d$.
\end{theorem}

We again refer to Section~\ref{sec:Lattice YMH} for more precise definitions of the lattices and the models;
and see Lemma~\ref{lem:exist},
Lemma~\ref{lem:inv}, Proposition~\ref{lem:4.7}, Lemma~\ref{lem:tightR}, Theorem~\ref{th:in1} and Section~\ref{sec:SDESG} for 
more complete formulations of
the above results.

Now we turn to our main results  regarding
 the long time behavior of the Langevin dynamics and functional inequalities.

First, we consider the case $M=\R^N$.
In contrast to pure lattice Yang--Mills, where the field $Q$ is consistently bounded since it takes values in a compact Lie group $G$, as explored in \cite{Seiler,SZZ22}, the scenario changes when the Higgs field assumes values in an unbounded space, such as $\mathbb{R}^N$. This introduces a significant challenge. In \cite{SZZ22}, the well-established Bakry--\'Emery condition is sufficient to derive the log-Sobolev inequalities. However, due to the unbounded nature of the Higgs field, directly applying the Bakry--\'Emery condition
(i.e. following \eqref{e:strategy} or \eqref{e:strategyGeo})
 becomes unfeasible.
 In fact, by employing the dynamics  \eqref{e:LSDE}, we can establish certain {\it moment bounds} on the Higgs field, as demonstrated in Lemma~\ref{lem:b:phi}.
A key strategy for proving functional inequalities with $\mathbb{R}^N$-valued Higgs fields then involves a disintegration
\begin{equ}[e:disint]
\dif\mu_\Lambda(Q,\Phi)=\dif\nu(Q)\dif \mu_Q(\Phi).
\end{equ}
 In this process, we initially fix the Yang-Mills field $Q$ and focus on the measure $\mu_Q$, which represents the regular conditional probability given $Q$. Subsequently, we prove the following results.

\begin{theorem}\label{main:2}
For $M=\mR^N$ with $m>0, \kappa>0, N>2$, suppose that
\begin{equ}[e:small-bkm]
8(d-1)|\beta|+\frac\kappa{m}+\frac{2\kappa^2}{m^2}
<\frac14-\frac1{2N}.
\end{equ}
The measure
$\mu_\Lambda$ \eqref{e:LYMH} satisfies the Poincar\'e inequality
with proportionality constant independent of the lattice size,
so that
any infinite volume limit $\mu$ satisfies the  Poincar\'e inequality. The  infinite volume dynamics is exponentially ergodic in $L^2(\mu)$.
\end{theorem}

See Theorem~\ref{th:poin}, Corollary~\ref{co:ergodi}
for the precise statements. In fact, the log-Sobolev and the Poincar\'e inequalities hold for both
$\mu_Q$ and $\nu$ in \eqref{e:disint},
where the proportionality constants in these inequalities for $\mu_Q$ only depend on $m$ and not on $Q$. See Lemma~\ref{log:1}.
In fact, the proof of these inequalities for $\nu$
does not directly follow from \eqref{e:strategy}
(the main reason is that the measure $\nu$ involves averaging over $\Phi$),
and we actually
 use the moment bounds of the Higgs part derived via the dynamics.
For convenience of the readers, we have the following diagram:
\begin{equ}[e:PfLem43]
\left.
    \begin{aligned}
        \text{B\' E} \stackrel{m>0}{\Longrightarrow}
        \text{LSI and Poincar\'e for } \mu_Q \, \eqref{log:muQ}& \\
        \mbox{Estimates on dynamics}
        \stackrel{Lemma~\ref{lem:b:phi}}{\Longrightarrow} \E_{\mu_Q} (\Phi_x^2) \leq 1/m
        &\\
        \text{B\' E} &
    \end{aligned}
    \right\} \Longrightarrow  \text{LSI and Poincar\'e for }\nu  \, \eqref{log:nu}
\end{equ}

However, despite that the log-Sobolev inequality holds
for the marginal law $\mu_Q$
for each fixed $Q$,
it does not automatically extend to $\mu_\Lambda$.
To prove Theorem~\ref{main:2}, further meticulous analysis of correlation functions is required.
As a consequence of the functional inequalities for $\mu_Q$,  one has a mass gap
for the measure $\mu_Q$ (the conditional measure with $Q$ fixed),
namely, the correlation of two observables under $\mu_Q$ decays exponentially in the separation distance,
see Lemma~\ref{lem:co}
(which only requires $m>0$ without any condition on $\beta,\kappa$).
This is crucial for the delicate analysis of the correlation functions in the proof of  Poincar\'e inequalities for $\mu_\Lambda$ and its infinite volume limit $\mu$, namely
Theorem~\ref{th:poin}.
The proof of this theorem is technical, but it can be outlined as follows:
\begin{equ}[e:PfThm47]
	\left.
	\begin{aligned}
	\mbox{Poincar\'e for  $\mu_Q$ and $\nu$ (Lemma~\ref{log:1})} &
	\\
		\text{ Poincar\'e for } \mu_Q \stackrel{m>0}{\Longrightarrow}
		 \mbox{mass gap for $\mu_Q$ (Lemma~\ref{lem:co})}
	 & \\
		\text{Decomposition of domain into smaller boxes} &
		\\
		\text{Delicate analysis of correlation functions} &
	\end{aligned}
	\right\} \Longrightarrow  
	\parbox{4cm}{Poincar\'e for $\mu_\Lambda$ and $\mu$ \\ (Theorem~\ref{main:2} or \ref{th:poin})}
\end{equ}

\medskip

Turning to the situations where the Higgs field takes values in a bounded target space $M$,
generally speaking, the strategy \eqref{e:strategyGeo} still goes through, up to technical details arising from the geometry of the manifolds.
By calculating the curvature of $M$ and estimating the Hessian operators, we verify the Bakry--\'Emery condition, and prove the following result,
under small $(\beta,\kappa)$ conditions.

\begin{theorem}\label{main3}
Let $M$ be the sphere  $\mS^{N-1}$ or the Lie group $G=SO(N)$ with $N>2$.
Assume the condition \eqref{e:K_YM}, as well as 
\begin{equ}[e:small-kappa]
 \left(4d\,\alpha +  (N-2)K_{\YM}^{-1}\right)|\kappa|
<
\frac12 - \frac1N
\end{equ}
where $\alpha=1$ for $M=\mS^{N-1}$ and $\alpha=4$ for $M=G$.
Then:


(1) The log-Sobolev and the Poincar\'e inequalities hold for $\mu_\Lambda$.

(2) The invariant measure $\mu$ of the infinite volume dynamics \eqref{e:LSDE} is unique and satisfy the log-Sobolev and the Poincar\'e inequalities.

(3) The infinite volume dynamics is exponentially ergodic w.r.t. $L^\infty(\mu)$.
\end{theorem}

See Lemma~\ref{lem:ba-S}, Lemma~\ref{lem:log-G}, Corollary~\ref{cor:SG-erg} for the precise statements.
Combining Theorem~\ref{main3}(2) and Theorem~\ref{main:1}, we know that for $M\in \{\mS^{N-1},G\}$, the infinite volume limit of $\mu_\Lambda$ is unique.

Remark that 
\eqref{e:small-bkm} is equivalent with assuming
 \eqref{e:K_YM} together with 
$ (\frac\kappa{m}+\frac{2\kappa^2}{m^2}) K_{\YM}^{-1}< \frac{1}N$.
So in both Theorem~\ref{main:2} and Theorem~\ref{main3}
we are assuming 
the {\it same} 
 small $\beta$ (i.e. strong Yang--Mills coupling) condition \cite[Assumption~1.1]{SZZ22}
as in the pure Yang--Mills  case,
together with additional small $\kappa/m$ or $\kappa$ conditions.
In fact $\kappa/m$ in \eqref{e:small-bkm} or $\kappa$ in \eqref{e:small-kappa}
 has to tend to $0$ if $K_{\YM}\searrow 0$. 

The constants for the log-Sobolev and the Poincar\'e inequalities 
in Theorem~\ref{main3}
will be given in 
\eqref{e:KS-S} \eqref{e:KS-G}.
We will not give an implicit constant for the Poincar\'e inequality
in Theorem~\ref{main:2} since as discussed above it does not simply follow from 
Bakry--\'Emery.

Now we turn to the problem of mass gap. The main result is as follows.

\begin{theorem}[Mass gap]\label{main4}
Assume \eqref{e:small-bkm} if  $M = \mR^N$, or \eqref{e:K_YM}\eqref{e:small-kappa} 
if  $M\in \{\mS^{N-1}, G \}$.

Then there is a positive mass gap
for any infinite volume limit of $\mu_\Lambda$.
Namely
the correlation of two observables which may depend on both $Q$ and $\Phi$ decays at least exponentially in their separation distance. 

If the observables are gauge invariant, exponential decay of correlation holds with
the condition \eqref{e:K_YM}\eqref{e:small-kappa} replaced by
a weaker condition \eqref{e:K_YM} and $|\kappa|<K_{\YM}/(2N)$.
\end{theorem}

See Theorem~\ref{co:mass}, Corollary~\ref{c:ex}, Corollary~\ref{cor:fix-gauge}, Remark~\ref{rem:CompareGG}.
In the case $M=\mS^{N-1}$, 
\cite[Theorem~3.18]{Seiler} has shown such a mass gap result 
for small $\beta,\kappa$ using cluster expansions.
However note that our 
conditions \eqref{e:small-bkm}  or \eqref{e:K_YM}\eqref{e:small-kappa} 
are uniform in $N$.
This allows us to extend our analysis to the large $N$ limit of observables. In Section \ref{sec:large}, we discuss such an application, namely we show the factorization property of observables as $N\to \infty$ for YMH models with $M\in \{\mS^{N-1},G\}$.  See Theorem~\ref{co:1}.
Remark that obtaining uniform in $N$ conditions on the coupling constants is one of the main improvements on the previous results \cite{Seiler}. 
(On the other hand, \cite[Theorem~3.18]{Seiler} also covers another regime where $\kappa,\beta$ are arbitrary but $\kappa/\beta$ is large, for certain special Lie groups such as $U(1)$ and $SU(2)$, see Section~\ref{sec:discuss-mass}.) 

For pure Yang--Mills considered in \cite{SZZ22}, where the field $Q$ is bounded, mass gap was proved using the Poincar\'e inequality.
To briefly review the argument at the heuristic level, recall  \eqref{e:gen-mass} that the key is to bound the commutator due to $P_t (fg) \neq P_t f P_t g$.
Writing $P_t = e^{t\cL}$, note that the Laplacian 
in the generator $\cL$ does not satisfy the product rule,
and therefore  to bound the commutator it is crucial to understand the quantity
$(\nabla_e P_t f)(\nabla_e P_t g)$ which is the discrepancy of the product rule, where $\nabla_e$ is a derivative  w.r.t. $Q_e$. 
Without $P_t$ (i.e. $t=0$), this quantity would vanish, since the supports of $f$ and $g$ are separated so $e$ cannot lie in both supports.
In other words, the argument essentially boils down to estimating commutators 
$[\nabla_e, P_t ]$ or  $[\nabla_e,\cL]$.
The key estimate in \cite[Cor~4.11]{SZZ22} is of the following form:
\begin{equ}[e:idea-comm]
\|\nabla_e P_t f \|_{L^\infty} 
\le \|\nabla_e f \|_{L^\infty} 
+ \int_0^t \sum_{\bar e\sim e} D_{e,\bar e} \|\nabla_{\bar e} P_s f\|_{L^\infty} ds
\end{equ}
where $\bar e$ and $e$ are the same or neighbor edges and $D_{e,\bar e}$ are 
coefficients bounded uniformly in $e,\bar e$.
Here one should think of $e$ as in the support of $g$, which is separated from the support of $f$ by a distance of order $|x-y|$,
so $\nabla_e f$ is zero. 
Iterating the above bound creates more and more neighbor edges, 
until one has iterated it $\approx |x-y|$ times, and then the first term on the RHS is non-zero.
This observation leads to a bound  on $\nabla_e P_t f$ of the form $\sum_{n=|x-y|}^\infty t^n/n!$ which decays as $e^{-c|x-y|}$ if $t\lesssim |x-y|$.

For $M=\mR^N$, as $\Phi_x$ is unbounded, we do not derive mass gap using the Poincar\'e inequality in Theorem \ref{main:2} directly along the line explained below \eqref{e:gen-mass}.
This is because the commutator estimate is more subtle;
one would not obtain a bound as in \eqref{e:idea-comm} with uniformly bounded coefficients
$D_{e,\bar e}$.
 Instead, we use disintegration \eqref{e:disint} to decompose the covariance $\cov_\Lambda(F,H)$ of two test functions $F, H$  as
	\begin{align}\label{eq:cor-i}
		\cov_\Lambda(F,H)=\E_\nu\cov_{\mu_Q}(F,H)+\cov_{\nu}(\E_{\mu_Q}F,\E_{\mu_Q}H),
	\end{align}
	(see \eqref{cov1} below).
As explained in \eqref{e:PfLem43} 
the Poincar\'e inequality holds  for both $\mu_Q$ and $\nu$.
The exponential decay 
of the first term on the RHS of \eqref{eq:cor-i}
will be shown  using the mass gap of $\mu_Q$.
For the second term on the RHS of \eqref{eq:cor-i},
writing $P_t^\nu$ for the Markov semigroup of the Langevin dynamics w.r.t. $\nu$ 
and $\cL^\nu$ the associated generator,  
similarly as in \eqref{e:gen-mass} we write
\begin{equ}\label{eq:cor-i1}
\cov_{\nu}(\E_{\mu_Q}F,\E_{\mu_Q}H)
=\cov_{\nu}(P_t^\nu \E_{\mu_Q}F,P_t^\nu \E_{\mu_Q}H) + \mbox{Comm}
\end{equ}
where the first term will be bounded using 
the Poincar\'e inequality for  $\nu$, 
and then the problem boils down to estimating the commutator
between the derivatives on the Lie group and the generator $\cL^\nu$.


It turns out that this commutator estimate is quite nontrivial, since  
after  integrating out $\Phi$
through the above disintegration, $\cL^\nu$ is nonlocal,
and we must utilize the mass gap of $\mu_Q$ to control this non-locality. 
Additionally, the moment bound on the Higgs field
derived from the dynamics is also useful in this control.
The proof strategy for  Theorem~\ref{main4} with $M=\R^N$
 can be summarized as
\begin{equ}[e:PfThm51]
	\left.
	\begin{aligned}
	 \eqref{eq:cor-i} \mbox{ and }\eqref{eq:cor-i1} 
	 & \\
	\mbox{mass gap for $\mu_Q$,  Poincar\'e for $\nu$}&
	\\
	\mbox{mass gap for $\mu_Q$} \Longrightarrow	\text{commutator estimates} &
	\end{aligned}
	\right\} \Longrightarrow  \text{ mass gap for }\mu
\end{equ}

We also note that
 Corollary~\ref{c:ex} and Corollary~\ref{cor:fix-gauge} both show mass gap for the case where $\Phi$ takes values
 in $G=SO(N)$, with the representation of $G$ on itself by left multiplication. However in the proof of Corollary~\ref{cor:fix-gauge} we exploited a gauge fixing
 (often called  unitary gauge or $U$-gauge), thanks to the fact that for any $N>1$  the stabilizer subgroup of $G$ for the identity element of $G$ is trivial.

 \br\label{rem:1}
 Our result can be readily extended to the case where $G=SU(N)$ as in \cite{SZZ22}.
Our proofs for the Poincar\'e inequality and existence of a positive mass gap
 also apply to alternative choices of the Higgs target spaces, see Remark~\ref{rem:Other-M}.
 \er

\subsection{Related results, literature, and discussions}
\label{sec:Discussion}

We conclude the introduction by discussing more  results, literature, problems that are related with the lattice Yang--Mills--Higgs models.

\subsubsection{Pure Yang--Mills on lattice}
The study of the pure Yang--Mills model on the lattice with various Lie groups $G$ using probabilistic methods
 has drawn much attention recently, see  \cite{Cha,Jafar,ChatterjeeJafar,Chatterjee16,MR4278289},
 \cite{SSZloop,SZZ22}, \cite{CPS2023}.
 Many questions such as the existence of a mass gap, functional inequalities, Dyson--Schwinger equation, large $N$ behavior, $1/N$ expansion, area laws, surface correspondence are investigated (see \cite{OS1978,Seiler} for earlier results on some of these questions).  The papers  \cite{SSZloop,SZZ22} developed a dynamical or stochastic approach which is close to the present paper.
 In this paper we focus on mass gap and functional inequalities (with some applications to large $N$) in the presence of a Higgs field, and the other questions just mentioned that were studied in the above references for pure YM remain to be investigated when coupled with matter field (such as a Higgs field or fermionic field;  in fact this is one of the open problems raised by \cite{CPS2023}) in the future.
Regarding mass gap results for the $U(1)$ lattice YM model,
let us  recall that  in 4D \cite{Guth1980,MR0649811} proved that there is a massless phase when $\beta$ is sufficiently large,
whereas in 3D \cite{MR0641914} proved that a non-zero mass gap exists for all $\beta$.

\subsubsection{Lattice Yang--Mills--Higgs in physics}
We also mention some physics work on lattice Yang--Mills--Higgs model, although the literature is vast.
Most physics literature is concerned with the $U(1),SU(2),SU(3)$ cases due to the relevance with the Standard Model.
The paper \cite{fradkin1979phase} contains some physics-level
discussion on different phases of the YMH model. 
The book \cite[Chapter~6.1]{MontvayMunster}  (see also more recent \cite{maas2014two}) summarizes a phase diagram  found by Monte-Carlo numerical methods
for different values  of the parameters (corresponding to our $\beta,\kappa$ and the shape of $V$ in \eqref{eq:YM_energy}).
We also refer to  \cite{BFS1979nuclear} which discussed some known and conjectural results about YMH.


\subsubsection{Related results on mass gap and Higgs mechanism}
\label{sec:discuss-mass}
In physics, in particular in the Standard Model,  the Higgs field plays the role of giving rise to masses of the particles via spontaneous symmetry breaking,
known as Higgs mechanism. See the standard physics textbook \cite[Chapter 23]{zinnbook}. Heuristically, the mechanism shows that one can exploit the Higgs part of the model, for instance by shifting the Higgs field to a particular vacuum value or performing a unitary gauge or ``U-gauge'' transformation, to create a new mass term in the Yang--Mills part (or essentially equivalently, on the lattice, a new term which renders the compact valued YM field likely to be close to identity).
This  mechanism is rigorously carried out in the abelian case
\cite{balaban1984mass,MR0836009,MR4025784} (see \cite{MR0912507} for a survey),
and  \cite[Section~4]{OS1978} \cite{Seiler} for certain non-abelian cases such as $SU(2)$.
In Seiler's book \cite[Theorem~3.18]{Seiler}, the utilization of U-gauge is effective and applicable to some choices of Lie groups (such as $G=SU(2)$) and representations,
yielding a mass gap
  in the regime where $\kappa/\beta$ is large enough. The restriction to particular Lie groups and representations is due to technical but important assumptions therein, mainly the complete symmetry breakdown condition.
As we are finishing this paper, Chatterjee \cite{Cha2024} proved the Higgs mechanism in a much more interesting regime where 
$\kappa/\beta$ is kept small
as $\kappa,\beta$ grow to infinity for $G=SU(2)$; see also the discussion in Sec.~\ref{sec:discuss-cont} below.
In this paper, 
 we  instead consider $G=SO(N)$ with general values of $N$ and a variety of representations;
 it is not our main purpose  to discuss the Higgs mechanism (only in Section~\ref{sec:gauge fixing} we use a similar gauge fixing idea). 
In fact for general Lie groups and representations, as the complete symmetry breakdown condition would be violated, it is not clear to us if the methods of \cite{OS1978,Seiler}  still apply.

\subsubsection{Estimates by dynamics}\label{sec:es-1}

It is noteworthy to highlight
that in this paper we allow
$\Phi$ to be unbounded i.e. $\mR^N$ valued, which poses one of the main challenges addressed in this paper.
The aforementioned work \cite{Seiler} assumes that
$\Phi$ takes values in a sphere, namely a bounded spin system
for which the analysis is simpler. (We also obtain new results
for sphere valued Higgs, see Theorems~\ref{main3} and \ref{main4}.)
Also, in \cite[Section~4]{OS1978}, a particular type of potential $V$ is chosen, which is a soft constraint on the Higgs field to lie within a sphere of radius $\bar R$ in the probabilistic sense.
The results \cite{balaban1984mass,MR4025784}  in the abelian case
has a term $\lambda|\Phi|^4$ in their potential with $\lambda>0$,
which strongly suppresses large fields.
In our results Theorems~\ref{main:2} and \ref{main4} with $M=\mR^N$,
we do not assume any such bounds on large fields, except for a mass term which only assumes very weak damping;
instead we prove bounds on the Higgs fields using the dynamics. Our proofs for the Poincar\'e inequality and mass gap also apply to  the case where $V(\Phi_x)=m |\Phi_x|^2+\lambda |\Phi_x|^4$, with $\lambda\geq0$ and $m\in\mathbb{R}$, and 
	 other choices of the Higgs target spaces. We refer to Remark \ref{rem:Other-M} for more discussions.

Note that strategies of using dynamics to control the fields and obtain qualitative  (including large $N$ limit) results on QFT models flourish recently, for instance \cite{GH21,AK20,HS21,MR4612651,MR4385125, MR4470243,SZZ23}.



\subsubsection{Discrete groups}
While our paper considers the Lie group $SO(N)$ as the structure group,
discrete structure groups
are also of much interest. See the recent papers \cite{FLV21,forsstromAbelian} for abelian lattice Higgs model with structure group $G=\mathbf Z_n$.
See also \cite{Adhikari2021} which is also for gauge theory with discrete groups coupled with Higgs fields. In these works computations of Wilson loop expectations to leading order at large $\beta$ and $\kappa$ have been obtained.

\subsubsection{Models in the continuum}\label{sec:discuss-cont}
As discussed in \cite{SZZ22}, it would indeed be intriguing to explore whether the Poincar\'e or the log-Sobolev inequalities, as well as the presence of a mass gap, persist as the lattice spacing tends to zero, particularly in scenarios where the continuum limits of these models are either established (or expected). See Remark~\ref{rem:cont} for some discussions; in particular it remains an interesting question to investigate whether the methods developed here can be applied to the large $\beta$ and $\kappa$ regimes.
For the abelian Yang--Mills--Higgs, the construction of continuum  and infinite volume limit for complex-valued $\Phi$ with potential of the form $V(\Phi)=(\Phi^2-1)^2$ was achieved
by \cite{BFS79,BFS80,BFS81} in 2D,
and then by \cite{King86I,King86II} in 3D (see also \cite{MR0677999,MR0679203,MR0701926} which proved lower and upper bounds uniform in lattice spacings on vacuum energy in 3D.) More recently, for  abelian YMH in 2D \cite{CC2022gauge} recently proved that the H\"older--Besov norm of the gauge field is  bounded uniformly in the lattice spacing.
The very recent paper  \cite{Cha2024} mentioned in Section~\ref{sec:discuss-mass}
proved that for $G=SU(2)$ in a  particular scaling regime one obtains {\it massive} Gaussian field in the continuum limit,
which is the first rigorous proof of the validity of the Higgs mechanism  in the continuum  in $d\ge 3$, although it remains open if there is a different scaling regime in which the limit is non-Gaussian.
Regarding the dynamical approach, the  Langevin dynamics for the Yang--Mills--Higgs models on
the three-dimensional continuous torus
have been recently constructed in \cite{CCHS3d};
see also \cite{CCHS2d} and \cite{Chevyrev2023}
for results on Langevin dynamics of the pure Yang--Mills in 2D.
Notably, for the pure Yang--Mills theory on $\mathbb{R}^2$, the
independence of separated gauge invariant observables
(in particular mass gap) is well-known, see \cite{MR3982691,LevySengupta17}.
Regarding other QFT models in the continuum,
recent progress includes the proof of the log-Sobolev inequalities for the $\Phi^4_{2,3}$ and sine-Gordon models (\cite{RolandPhi4, RolandSG}), as well as the 1D nonlinear $\sigma$-model (\cite{AnderssonDriver, Hairer16, StringManifold}) where the log-Sobolev inequalities, ergodicity, and non-ergodicity depending on the curvature of the target manifolds were established in \cite{RWZZ17, CWZZ18}.


\bigskip
{\bf Structure of the paper.}
The paper is organized as follows.
In Section \ref{sec:not}, we provide the definition of the lattice Yang--Mills--Higgs model and list some common choices for the Higgs target spaces,
and we recall some relevant geometric background. Section \ref{sec:YMHSDEs} establishes the global well-posedness of SDEs \eqref{e:LSDE} and proves Theorem \ref{main:1}. 
Section \ref{sec:ergo} is dedicated to proving the exponential ergodicity of SDEs \eqref{e:LSDE} through functional inequalities. In Section \ref{sec:mass}, we present the proof of the mass gap. Finally Section \ref{sec:large} discusses the large $N$ limits of gauge-invariant observables.

\bigskip
{\bf Acknowledgments.}
We would like to thank Scott Smith for very helpful discussions,
including the derivations and estimates for the lattice YMH dynamics and many other aspects. We would also like to thank Roland Bauerschmidt for discussions on functional inequalities and related problems.
H.S. gratefully acknowledges financial support from NSF grants
DMS-1954091 and CAREER DMS-2044415. R.Z. and X.Z. are grateful to
the financial supports   by National Key R\&D Program of China (No. 2022YFA1006300) and the support from the NSFC (No. 12426205).
R.Z. gratefully acknowledges financial support from the NSFC (12271030).  X.Z. is grateful to
the financial supports   by   the NSFC (No. 12595281,
 12288201). This work is funded by the Deutsche Forschungsgemeinschaft (DFG, German Research Foundation) – Project-ID 317210226--SFB 1283.

\section{Notation and definitions}\label{sec:not}

In this section we precisely set up the notation and define the model.

\subsection{Lattice Yang--Mills--Higgs model}
\label{sec:Lattice YMH}

Let $\mT^d=[-1,1]^d$ be the $d$-dimensional torus
and $\Lambda_{L}=\mZ^d\cap L\mT^d$ be a finite $d$ dimensional lattice
with side length $L>0$ and unit lattice spacing, and we will consider various functions on it with periodic boundary conditions.
We will often write $\Lambda=\Lambda_{L}$ for short. Our analysis can be extended to other boundary conditions, such as free or Dirichlet boundaries, as the proofs presented in the following sections remain valid with small modifications. We have chosen periodic boundary conditions for the sake of simplicity.

Each lattice edge $e$ of $\mZ^d$ is  oriented with the starting vertex denoted by $u(e)$ and the end vertex denoted by $v(e)$. We write $e=(u(e),v(e))$.
Let $E^+$ (resp. $E^-$) be the set of positively (resp. negatively) oriented edges,
and  denote by $E_{\Lambda}^+$, $E_{\Lambda}^-$ the corresponding subsets
of edges with both vertices in ${\Lambda}$.  Define $E\eqdef E^+\cup E^-$, $E_\Lambda\eqdef E^+_\Lambda\cup E^-_\Lambda$.

A {\it path}  is defined to be a sequence of edges $e_1e_2\cdots e_n$ with $e_i\in E$ and $v(e_i)=u(e_{i+1})$ for $i=1,2,\cdots, n-1$. The path is called closed if $v(e_n)=u(e_1)$.
A {\it plaquette} is a closed path of length four which traces out the boundary of a square. Let $\cP$ denote the set of plaquettes and let $\cP^+$ be the subset of plaquettes $p=e_1e_2e_3e_4$ such that
$u(e_1)$ is lexicographically the smallest among all the vertices in $p$ and  $v(e_1)$ is the second smallest.
Also, let $\CP_{\Lambda}$ be the set of plaquettes whose vertices are all in $\Lambda$, and
$\CP^+_{\Lambda}=\cP^+\cap \CP_{\Lambda}$.

Let
$G=SO(N)$. The lattice Yang--Mills--Higgs model is defined by the following action:
\begin{equ}[e:CS]
\mathcal S_{\YMH}(Q,\Phi) =
N\beta  \sum_{p\in \CP^+_\Lambda} \Tr(Q_p)
- N\kappa \sum_{e\in E^+_\Lambda }
|Q_e \Phi_y - \Phi_x|^2  -m N\sum_{z\in \Lambda}|\Phi_z|^2
\end{equ}
where $\beta,\kappa,m$ are coupling constants.
We now specify the meaning of all the objects appearing in \eqref{e:CS}.

The discrete Yang--Mills field $Q=(Q_e : e\in E_{\Lambda}^+)$ is a collection of $G$-matrices.
The first term on the RHS of \eqref{e:CS} is the lattice Yang--Mills model proposed by Wilson  \cite{Wilson1974}, where
 $Q_p \eqdef Q_{e_1}Q_{e_2}Q_{e_3}Q_{e_4}$ for a plaquette $p=e_1e_2e_3e_4$, and
$x,y$ are the beginning and ending points of $e$ (we will write  $e=(x,y)$).  Throughout the paper we define $Q_{e}\eqdef Q_{e^{-1}}^{-1}$ for $e \in E^{-}$, where $e^{-1}$ denotes the edge with orientation reversed.
It is well-known that
by the Baker--Campbell--Hausdorff formula
this term (suitably re-centered around identity $I_N$) is the discrete approximation of the continuum Yang--Mills action in \eqref{eq:YM_energy} (c.f. \cite[Section~3]{Chatterjee18}), namely,
with $Q_e=e^{\eps A_j}$ where $e\in E^+$ is parallel with the $j$th axis and $A_j$ takes values in the Lie algebra $\mfg$ of $G$,
\begin{equ}[e:YMeps]
\sum_p \Tr (I_N-Q_p) \approx
\frac{\eps^{4-d}}4  \int |F_A(x)|^2 \dif x,
\end{equ}
where $\eps$ is the lattice spacing. In this paper we fix $\eps=1$ and drop the identity $I_N$ since it only contributes to the action a constant.
Note that there are various other versions of the lattice Yang--Mills actions besides the above Wilson action,
for example Villain's action defined via heat kernels (c.f. \cite[Section~8]{Driver89}, \cite{Levy03,Levy06}),
or Manton's action \cite{Manton80} defined via the Riemannian metric on $G$; see also \cite[Section~2]{Chevyrev2023}.
But we will not consider these choices in this paper.

Turning to the other terms in \eqref{e:CS},
the field $\Phi$ takes values in a space $M$ where the Lie group $G$ acts on.
The space $M$ can be either a vector space endowed with an inner product
in which case the $G$-action is an orthogonal representation on $M$,
or a non-linear space, more precisely a Riemannian manifold on which $G$ has a metric preserving action.
Note that $|Q_e \Phi_y - \Phi_x|^2$
is a generalization of the usual kinetic term $(\nabla\Phi)^2$,
but this term is now coupled with the field $Q$ in order to obtain a gauge invariant model.
We will call
\begin{equ}[e:cov-der]
Q_e \Phi_y - \Phi_x
\end{equ}
 the (discrete) gauge covariant derivative of $\Phi$,
which is named due to its gauge covariance property,
i.e. transforming $(Q,\Phi)$ in the covariant derivative by a gauge transformation $g$ amounts to
the action of $g$ on the entire covariant derivative, see \eqref{e:cov-D}.

We list and discuss some  of the common and natural choices of the Higgs target spaces $M$, endowed with a distance $|\cdot |$, and the action of $G$ on it.
 \footnote{``Defining representation'' is also called the ``fundamental representation'' of $G=SO(N)$ on $\mR^N$ in physics (c.f. \cite[Section~15.4]{MR1402248}). Note that this is different from the  ``fundamental representation'' in representation theory which refers to an irreducible representation of a semisimple Lie group whose highest weight is a fundamental weight.} $\!\!$
 \footnote{We endow the distance  $|\cdot|_{\R^N}$ on
 $\mS^{N-1}$ since it is more convenient than the sphere
(Riemannian) distance on $\mS^{N-1}$; of course they are equivalent.}

\begin{enumerate}[itemsep=5pt]
\item \label{case1}
$M=\R^N$.
\\
$|\cdot |$ is  the Euclidean distance $|\cdot|_{\R^N}$.
\\
The $G$-action is given by the defining representation, namely, $Q\Phi$ simply means the multiplication of a matrix $Q$ and a vector $\Phi$.
\item  \label{case2}
$M= \mS^{N-1} \subset \mR^N$, the unit sphere.
\\
 $|\cdot |$ is the induced Euclidean distance $|\cdot|_{\R^N}$.
\\
The $G$-action is the same as above, in other words, the
rotation of the sphere by $SO(N)$.
\item \label{case3}
$M=G$.
\\
 $|\cdot |$ is the distance for a bi-invariant Riemannian metric on $G$ (see Section~\ref{sec:mani}).
\\
  The $G$-action is given by the left multiplication of $G$ on itself,
  namely, $Q\Phi$ is a matrix product.

\item \label{case4}
$M=G$.
\\
 $|\cdot |$ is bi-invariant Riemannian as above.
\\
  The $G$-action is given by the adjoint action, namely  $Q \circ \Phi
 \eqdef
Q  \Phi Q^{-1}$.
\item  \label{case5}
$M=\mfg$ ($N\times N$ skew-symmetric matrices).
\\
 $|\cdot |$ is the Hilbert-Schmidt inner product (see \eqref{e:HS}).
\\
  The $G$-action is given by the adjoint representation
   $Q \circ \Phi \eqdef
Q  \Phi Q^{-1}$.
\item  \label{case6}
$M=\{N\times N \mbox{ symmetric matrices}\}$.
\\
 $|\cdot |$  Hilbert-Schmidt as above.
\\
  The $G$-action is by adjoint as above.
\end{enumerate}

Note that in case \eqref{case2}, $|\Phi|^2=1$ and in case \eqref{case3},
$|\Phi|^2=\Tr (\Phi\Phi^t)=N$,
so they only contribute to the overall normalization factor and we can simply drop
the last term in \eqref{e:CS}.
The sphere constraint as in case \eqref{case2} often leads to simplification
 (c.f. \cite[Section~3]{Seiler})
  since it means that the Higgs field $\Phi$ is well bounded.
It could be understood as a formal limit as $\lambda\to \infty$
for a potential of the form $-V(\Phi) = -\lambda (|\Phi|^2-1)^2$ in \eqref{e:CS} with $\Phi\in \mR^N$,
so that $V(\Phi)\to+\infty$ if $\Phi\notin\mS^{N-1}$.
Similarly,  the choice $M=G$ is also bounded which is convenient for the analysis.

Also, we note that some of these models degenerate into other interesting models:
without  $Q$ (taking $Q\equiv 1$),
\begin{itemize}
\item Case \eqref{case1}  reduces  to the $N$-component massive Gaussian free field;
\item Case \eqref{case2} reduces to $N$-component Heisenberg model c.f. \cite{MR574175};
\item Case \eqref{case3} becomes a version of  orthogonal  multi-matrix model,
        c.f. \cite[Section~9]{CGM2009} or \cite{GuionnetNovak};
\item Case \eqref{case6}   becomes  a chain of (real) random matrices, c.f. \cite{Eynard2009}.
\end{itemize}
Case \eqref{case5} is of particular interest in physics (e.g. $N=2,3$),
c.f. \cite[(15.83)]{MR1402248} (in continuum) or \cite{Drouffe1984} (on a lattice).

Moreover, the case \eqref{case3} is special since for any $N>1$, a so-called ``complete breakdown of symmetry'' property holds
(namely the stability group of identity $I_N\in M=G$ is trivial, see \cite[Section~3.c]{Seiler}) which makes it easier to perform a U-gauge fixing trick, see Section~\ref{sec:Discussion} and Section~\ref{sec:gauge fixing}.

Although we do not study scaling limits in this paper,
remark that \eqref{e:cov-der} are indeed discrete approximations of
$\dif_A \Phi$ in \eqref{eq:YM_energy} under proper scalings.
Taking $Q_e=e^{\eps A_j}$ where $e=(x,y)\in E^+$ as  in \eqref{e:YMeps}, one has
\begin{equs}
Q_e \Phi_y - \Phi_x &=e^{\eps A_j}\Phi_y - \Phi_x
=
(1+\eps A_j)\Phi_y - \Phi_x + O(\eps^2)
\\
&\qquad\qquad\qquad\;\;\;
=
\eps(\partial_j \Phi_x + A_j\Phi_x)
+ O(\eps^2)
\qquad\qquad\qquad\qquad\quad
 (\mbox{for } \Phi_x\in \R^N),
\\
Q_e \Phi_y Q_e^{-1}- \Phi_x
&=e^{\eps A_j} \Phi_y e^{-\eps A_j}  - \Phi_x
=
(\Phi_y-\Phi_x )+ \eps A_j \Phi_y - \eps\Phi_y A_j  + O(\eps^2)
\\
&\qquad\qquad\qquad\qquad\quad\, =
\eps \partial_j \Phi_x + \eps [A_j,\Phi_x] + O(\eps^2)
\qquad\qquad\qquad
 (\mbox{for } \Phi_x\in \mfg).
\end{equs}
For $M=G$,
letting $\Phi_z = e^{\eps \Psi_z}$
with $\Psi_z\in \mfg$,
one has, by the Baker--Campbell--Hausdorff formula,
\begin{equs}
Q_e \Phi_y Q_e^{-1}- \Phi_x
&=e^{\eps A_j} e^{\eps \Psi_y} e^{-\eps A_j}  - e^{\eps \Psi_x}
=
e^{\eps \Psi_y + \eps^2 [A_j,\Psi_y] +O(\eps^3)}
-
e^{\eps \Psi_x}
\\
&=
\eps(\Psi_y-\Psi_x) + \eps^2 [A_j,\Psi_y] +O(\eps^3)
=
\eps^2 \partial_j \Psi_x + \eps^2 [A_j,\Psi_x] +O(\eps^3).
\end{equs}
Thus the second term in \eqref{e:CS} is approximately
\begin{equ}[e:Higgseps]
\eps^{2-d} \int_{\mR^{d}}
 |\dif_A \Phi (x)|^2 \dif x
 \quad
M\in\{\R^N,\mfg\},
\qquad
\eps^{4-d} \int_{\mR^{d}}
 |\dif_A \Phi (x)|^2 \dif x
  \qquad
M=G.
\end{equ}

As mentioned below \eqref{eq:YMH_measure},
in the general geometric setting,
given a $G$-principal bundle,
the field $Q$ should be thought as the discrete analogue of the connection  on the principal bundle;
with a ``fiber space'' $M$ and a $G$-action on it,  one can define an associated bundle
and $\Phi$ should be thought as the discrete analogue of a section of the associated bundle.
The above examples are special cases of this construction,
which correspond to associated vector bundles of rank $N$, associated sphere bundles,
associated adjoint bundles with fiber $\mfg$,
associated bundle with fiber $G$ with $G$-action  by adjoint,
or  associated bundle which is simply just the principal bundle itself.
Of course in this paper all the bundles are trivial (i.e. product bundles).

 We then define the probability measure
\begin{equation}\label{measure}
\dif\mu_{\Lambda, N, \beta,\kappa,m}(Q,\Phi)
 := Z_{\Lambda, N,\beta,\kappa,m}^{-1}
\exp\Big(\mathcal S_{\YMH} \Big)
\prod_{e\in E^+_\Lambda} \dif\sigma_N(Q_e)
\prod_{z\in \Lambda} \dif\Phi_z\, .
\end{equation}
Here $Z_{\Lambda, N,\beta,\kappa,m}$ is a normalization factor and $\sigma_N$ is the Haar measure on $G$, and $\dif\Phi_z$ means Lebesgue/uniform/Haar measures on  $M$. We also write $\mu_\Lambda=\mu_{\Lambda, N, \beta,\kappa,m}$ for notation's simplicity.

In all the cases \eqref{case1}--\eqref{case6}, we can write
\begin{equs}[e:rev-e]
|Q_e \Phi_y - \Phi_x|^2
&=\Tr ((Q_e \Phi_y - \Phi_x)(Q_e \Phi_y - \Phi_x)^t)
\\
&=\Tr ((Q_e \Phi_y - \Phi_x)^t(Q_e \Phi_y - \Phi_x))
\\
&=| \Phi_y - Q_e^{-1}\Phi_x|^2
\end{equs}
where the $\Tr$ in the second line can be omitted in case \eqref{case1} and case \eqref{case2}, and $Q_e\Phi_y$, $Q_{e^{-1}}\Phi_x$ are replaced by $Q_e\circ \Phi_y$, $Q_{e^{-1}}\circ \Phi_x$, respectively in case \eqref{case4}--\eqref{case6}.
From the second line above one also has
\begin{equ}[e:h]
|Q_e \Phi_y - \Phi_x|^2
=\Tr[\Phi_y^t \Phi_y - \Phi_y^t Q_e^t \Phi_x - \Phi_x^t Q_e \Phi_y
+ \Phi_x^t \Phi_x],
\end{equ}
where the $\Tr$ can be omitted for $\R^N$ or sphere valued Higgs fields,  and $Q_e\Phi_y$, $Q_{e}^t\Phi_x$ are replaced by $Q_e\circ \Phi_y$, $Q_{e}^t\circ \Phi_x$, respectively in case \eqref{case4}--\eqref{case6},
and $\Phi_x^t \Phi_x=\Phi_y^t \Phi_y = 1$ if
$\Phi_x\in \mS^{N-1}$.

For  any $g:\Lambda \to G$, we define gauge transformations by
\begin{equ}[e:gauge]
Q_e \mapsto g_x Q_e g_y^{-1},
\qquad
\Phi_x \mapsto g_x\Phi_x
\end{equ}
in cases \eqref{case1}--\eqref{case3}, and in case \eqref{case4}--\eqref{case6} the second transformation is replaced by $\Phi_x \mapsto g_x\Phi_x g_x^{-1}$.

\bl
\label{lem:gauge-inv}
The measure $d\mu_{\Lambda, N, \beta,\kappa,\lambda}$ is invariant under gauge transformations.
\el

\begin{proof}
This is well-known but we record the proof for completeness.
The gauge invariance can be seen in the following way.  
For the discrete covariant derivative,  in cases \eqref{case1}--\eqref{case3} it transforms as
\begin{equ}[e:cov-D]
Q_{e}\Phi_{y}-\Phi_{x} \quad  \mapsto \quad
g_{x}Q_{e}g_{y}^{-1}g_{y}\Phi_{y}-g_{x}\Phi_{x}
= g_{x}  (Q_{e}\Phi_{y}-\Phi_{x}).
\end{equ}
So $|Q_{e}\Phi_{y}-\Phi_{x}|^{2}$ is invariant, again since $g_{x} \in SO(N)$.  Similarly in case \eqref{case4}--\eqref{case6}, we have
\begin{equ}[e:cov-D4]
Q_{e} \Phi_{y}Q_{e}^{-1} -\Phi_{x} \quad  \mapsto \quad
g_{x}Q_{e}g_{y}^{-1}  g_{y} \Phi_{y}g_{y}^{-1} (g_{x}Q_{e}g_{y}^{-1})^{-1}
-g_{x}\Phi_{x}g_{x}^{-1}
= g_{x}  (Q_{e}\Phi_{y}Q_e^{-1}-\Phi_{x})g_{x} ^{-1}.
\end{equ}
So we again have the desired invariance.

Finally, for the pure YM action, writing $Q_{p}=Q_{e_{1}}\cdots Q_{e_{4}}$ and expressing the edges as $e_{i}=(u(e_{i}),v(e_{i}))$, we find that it transforms as
$$
\text{Tr}Q_{p} \mapsto \text{Tr}[g_{u(e_{1})}Q_{e_{1}}g_{v(e_{1})}^{-1} \cdots Q_{e_{4}}g_{v(e_{4})}^{-1}]=\text{Tr}[Q_{e_{1}}\cdots Q_{e_{4}}g^{-1}_{v(e_{4})}g_{u(e_{1}) }]=\text{Tr}[Q_{p}]
$$
 where the middle terms cancel because $e_{1}\cdots e_{4}$ is a path and the final terms cancel because of cyclic invariance and the fact that the $e_{1}\cdots e_{4}$ is closed.
\end{proof}

Due to Lemma~\ref{lem:gauge-inv}, observables which are gauge invariant in this model are of particular interest.
There are two classes of well-known gauge invariant observables:
Wilson loops and Wilson lines.
Given a loop $ l  = e_1 e_2 \cdots e_n$ (which means a closed path as defined above \footnote{The definition of loops are slightly different from e.g. \cite{Cha} or \cite{SSZloop} since we do not quotient out cyclic relations and backtracks here. In this paper we do not need to enumerate these loops as in the study of master loop equations, so these issues are unimportant.}),
 the Wilson loop variable $W_l$ is defined as
\begin{align}\label{de:wlo}
W_l = \Tr (Q_{e_1}Q_{e_2}\cdots Q_{e_n})\;.
\end{align}
Now given an open line $l  = e_1 e_2 \cdots e_n$ (meaning a path which is not necessarily closed),
we  define the  Wilson line variable as
\begin{align}\label{de:wl}
W_l = \text{Tr} ( \Phi_{u(e_1)}^t Q_{e_1}Q_{e_2}\cdots Q_{e_n}\Phi_{v(e_n)} )\;.
\end{align}
Here the trace does not have any effect in the vector cases, i.e. cases \eqref{case1}\eqref{case2}.
The loop variables and the line variables
are invariant under gauge transformations \eqref{e:gauge} for any $g$.
Remark that for case \eqref{case3} where $M=G$, the Wilson lines are gauge invariant even without the trace, but we define it with a trace so that it is a scalar random variable. 

We remark that the loop and line variables could be viewed
as ``macroscopic'' versions of the terms
in \eqref{measure}: $\tr(Q_p)$ is just a loop variable
where the loop is  a single plaquette, and the terms
on the RHS of \eqref{e:h} are just line variables
where the line has one or zero edge.

As mentioned in Section~\ref{sec:intro}, we will consider the cases \eqref{case1}\eqref{case2}\eqref{case3} in the rest of the paper, since they are prototypical choices that exemplify the different challenges and strategies; see Remark~\ref{rem:Other-M} for some discussion on the other cases.

\br
For the YM term
$N\beta  \Re \sum_{p} \Tr(Q_p)$ in \eqref{e:CS},
we have chosen the well-known t'Hooft scaling $N\beta$ (\cite{tHooft1974}). Concerning Higgs part in \eqref{e:CS}, we choose the scaling with coefficients $N\kappa$ and $N m$, so that one can factor out an overall $N$
from $\mathcal S_{\YMH}$. This appears to us to be a convention in the physics literature on YM coupled with matter fields, see e.g. \cite[(2.1)]{olmsted1998}. This choice is also such that
in the case of sphere valued Higgs, when the field $Q$ is absent the scaling is consistent with that of
the large $N$ classical Heisenberg model e.g. \cite[(1a)]{MR574175}.
Remark that by a change  of variable $\Phi\to \sqrt N \Phi$,
the coefficients $(N\kappa, N m)$  in \eqref{e:CS} become $(\kappa,m)$, but we do not follow this convention.
\er

\br\label{rem:cont}
The above formal derivations \eqref{e:YMeps}, \eqref{e:Higgseps} 
indicate that to pass our results on functional inequalities or mass gap to continuum, one has to essentially extend the results from small $(\beta,\kappa)$ to large $(\beta,\kappa)$.
\er

\subsection{Some geometric background}\label{sec:mani}

We collect some notations and basic  geometric facts   used throughout the paper.

\subsubsection{Riemannian manifolds.}
Let $M$ be a Riemannian manifold.
We denote by $C^\infty(M)$  the space of real-valued smooth functions on $M$.
For $x\in M$ we denote by $T_xM$  the tangent space at $x$ with inner product $\<\cdot,\cdot\>_{T_xM}$.
For $X \in T_xM$, we write $Xf$ or $X(f)$ for the differentiation of $f$ along $X$ at $x$.
%
Let $\nabla$  be the Levi-Civita connection,
and $\nabla_{Y}X$ be the covariant derivative of a vector field $X$ along the vector field $Y$.
Recall that $(\nabla_{Y}X)(x)$  depends on $Y$ only via $Y(x)$ for $x\in M$ (e.g. \cite[Remark~2.3]{MR1138207}). Set $[X,Y]:=XY-YX$ for vector field $X,Y$.

%
%

For $f\in C^\infty(M)$, we denote by $\nabla f$ the gradient and $\Hess(f) $ or $\Hess_f$  the Hessian. We write $\Ric$ for the Ricci curvature.
$\Hess$ and $\Ric$ are  symmetric 2-tensors: for vector fields $X,Y$, $\Hess_f(X,Y)(x)$ and $\Ric(X,Y) (x)$ depend only on $X(x)$ and $Y(x)$ for $x\in M$.

\subsubsection{Lie groups and algebras.}
Let $M_N$ be the space of real  $N\times N$ matrices
endowed with the Hilbert-Schmidt inner product
\begin{equ}[e:HS]
	\< X,Y\> =  \Tr (X Y^t)
	\qquad  \forall X,Y\in M_N.
\end{equ}
For any matrix $M\in M_N$ we write $M^t$ for the 
 transpose of $M$.
 Below we always choose
$$
G=SO(N),\qquad \mfg=\so(N)
$$
where $\so(N)$ is the Lie algebra of skew-symmetric matrices,
and write $d(\mfg)=\frac12 N(N-1)$ for its dimension.



Every $X\in \mfg$ induces a right-invariant vector field  $\widetilde X$
on $G$, and for each $Q\in G$, $\widetilde X(Q)$  is just given by $XQ$ since $G$ is a matrix Lie group. Indeed, given any $X \in \mfg$, the curve $t \mapsto e^{tX} Q = Q+tXQ + O(t^2)$
near $t=0$.

The inner product on $\mfg$ induces an inner product
on the tangent space at every $Q\in G$
via the right multiplication on $G$.
Hence, for $X,Y\in \mfg$, we have $XQ,YQ \in T_{Q}G$, and their inner product is given by $\Tr((XQ)(YQ)^t) = \Tr(XY^t)$.
This yields a bi-invariant Riemannian metric on $G$.

For  $f\in C^\infty (G)$ and $X \in \mfg$,
the right-invariant vector field $\widetilde X$ induced by $X$ acts on $f$ at $Q\in G$ by the right-invariant
derivative
$\widetilde X f (Q)=  \frac{\dif}{\dif t}|_{t=0} f(e^{tX} Q)$.
One also has
$\widetilde{[X,Y]} = [\widetilde X, \widetilde Y]$
where the LHS  is the Lie bracket on $\mfg$ and  the RHS is the vector fields commutator.
Also, for the Levi-Civita connection $\nabla$ we have
$\nabla_{\widetilde X} ( \widetilde Y) =\frac12 \widetilde{ [X,Y]  }$.
See e.g. \cite[Appendix~F]{AGZ}.

%

\medskip

\subsubsection{Brownian motions.}
Recall that Brownian motions are defined on arbitrary Riemannian manifolds, in particular Lie groups.
Let $\mathfrak{B}$ and $B$  be the Brownian motions on $G=SO(N)$ and its Lie algebra $\mfg$ respectively.
One has
$\E \big[  \<B(s),X \> \<B(t),Y \>  \big] = \min(s,t) \< X,Y \>$
for $X,Y \in \mfg$.
Writing $M^{ij}$ for the $(i,j)$th entry of a matrix $M$, one has
$\dif\< B^{ij}, B^{k\ell}\>=\frac12(\delta_{ik}\delta_{j\ell}-\delta_{i\ell}\delta_{jk})\dif t$  (c.f. \cite[(2.5)]{SSZloop}).
By  \cite[Sec.~1.4]{Levy11},  $\mathfrak{B}$ and $B$ are related through the following SDE:
\begin{equ}[e:dB]
	\dif \mathfrak B = \dif B \circ \mathfrak B = \dif B\, \mathfrak B
	+ \frac{c_\mfg}{2} \mathfrak B \dif t,
\end{equ}
where $\circ$ is the Stratonovich product, and $\dif B\, \mathfrak B$ is in the It\^o sense.
Here the constant $c_\mfg$ is determined by
$ \sum_{\alpha} v_\alpha^2  =c_\mfg I_N$ where
$(v_\alpha)_{\alpha=1}^{\dg}$ is an orthonormal basis of $\mfg$.
Moreover (c.f. \cite[Lem.~1.2]{Levy11}),
\footnote{Note that in \cite[Lem.~1.2]{Levy11}, the scalar product differs from \eqref{e:HS} by a scalar multiple depending on $N$ and $\mfg$, so we accounted for this in the expression for $c_\mfg$ above.}
\begin{equ}[e:c_g]
	c_{\so(N)} =  -\frac12(N-1) .
\end{equ}

Turning to the case of the sphere,
now let $\mathfrak{B}$ be the Brownian motion on unit sphere $\mS^{N-1}$.  We can choose suitable vector fields $\sigma_k$ with $\sum_k\sigma_k^2$ given by Beltrami--Laplace on $\mS^{N-1}$ and then  write
$\dif \mathfrak{B}= \sum_{k=1}^{N}\sigma_k(\mathfrak{B})\circ\dif B^k$ with  $B$ being $\mR^N$-valued Brownian motion, where $\circ$ is the Stratonovich product.
More concretely,
$\mathfrak{B}=(\mathfrak{B}^1,\cdots,\mathfrak{B}^N)$  is the solution to the SDEs
\begin{align}\label{bm:sp}
\dif \mathfrak{B}^i=	&\sum_{j=1}^N(\delta_{ij}-\mathfrak{B}^i\mathfrak{B}^j)\circ \dif B^j,\quad 1\leq i\leq N.
\end{align}
See \cite[page 83, Example 3.3.2]{Hsu}.

\subsubsection{Product manifolds and Lie groups}\label{sec:product}

The  configuration space of our fields will be a product of
Riemannian manifolds and Lie groups over the lattice vertices and edges. See \cite[Section~2.1]{SZZ22} for review of these definitions;
in particular recall that the tangent space (resp. Lie algebra) of such a product  is given by the direct sum of the tangent spaces (resp. Lie algebras).
For a finite collection of Riemannian manifolds $(M_e)_{e\in A}$ where $A$ is a finite set and all $M_e$ are the same Riemannian manifold $M$, we write  the product as $M^{A}$.

In this case, we will sometimes view a ``local'' tangent vector as a ``global'' tangent vector in the following way:
given a point $x= (x_e)_{e\in A} \in M^A$,
if $u_e \in T_{x_e} M_e$ for some $x_e \in M_e$,
we will sometimes view $u_e$ as a tangent vector in $T_x M^{A}$
which has zero components for all $\bar e \neq e$.
Continuing with this notation,
if $(v_e^i)^{i=1,...,d}$ is a basis (resp. orthonormal basis) of $T_{x_e} M_e$, then
$(v_e^i)_{e\in A}^{i=1,...,d}$ is a basis (resp. orthonormal basis) of $T_{x} M^A$.


With similar notation as above we write the product Lie group $G^A$ and its Lie algebra $\mfg^A$ for a finite set $A$.
Given $X \in \mfg^A$, the exponential map $t \mapsto \text{exp}(tX) \in G^A $ is also defined pointwise as
$\text{exp}(tX)_{e}\eqdef e^{tX_{e}}$ for each $e\in A$.

For a finite lattice $\Lambda$, we define the {\it configuration space}
for the YMH model as
$$
\cQ_\Lambda\eqdef
G^{E_{\Lambda}^+}\times M^{\Lambda},
\qquad
M\in \{\mR^N, \mS^{N-1},G\}
$$
consisting of all maps $(Q,\Phi)$ with
$$
Q: E_{\Lambda}^{+}  \ni e \mapsto Q_{e} \in G
\qquad
\mbox{and}
\qquad
\Phi: \Lambda \ni x \mapsto \Phi_x\in M.
$$
Write $d_M$ for the dimension of $M$, namely
$d_M=N$ for $M=\mR^N$, $d_M=N-1$ for $M=\mS^{N-1}$,
and  $d_M=\frac12 N(N-1)$ for $M=G$.
Let  $\q_\Lambda$ be the Lie algebra of the Lie group 
$G^{E_{\Lambda}^+}$
which is just the corresponding direct sum of $\mfg$.

For any matrix-valued functions $A,B$ on $E_{\Lambda}^+$, we denote by $AB $
the pointwise product $(A_e B_e)_{e\in E_{\Lambda}^+}$.


With the above notation,
one has
$$
T_{(Q,\Phi)} \cQ_\Lambda =
\Big\{
(XQ,v)=(X_eQ_e,v_x)_{e\in E_{\Lambda}^{+},x\in \Lambda}
\quad \mbox{where}\;\;
 X\in \q_\Lambda,v\in T_{\Phi}M^\Lambda
 \Big\}
$$
endowed with an inner product
$$
\Big\langle
(XQ,u),(YQ,v)
\Big\rangle_{T_{(Q,\Phi)} \cQ_\Lambda}
= \sum_{e \in E_{\Lambda}^{+} } \text{Tr}(X_{e}Y_{e}^{t})
+ \sum_{x\in \Lambda}\langle u_x,v_x\rangle_{T_{\Phi_x}M} \;.
$$
Recalling our above convention of viewing ``local'' tangent vectors as ``global'' ones,
we can find a basis of $T_{(Q,\Phi)} \cQ_\Lambda $:
\begin{equ}[e:TQ-basis]
T_{(Q,\Phi)} \cQ_\Lambda
=
\mbox{Span}
\Big\{
X_e^iQ ,v_x^j \; :\;
 e\in E^+_{\Lambda} , \;
 x\in \Lambda, \;
 1 \le i \le \dg, \;
 1\le j\le d_M
\Big\} \;,
\end{equ}
where for each $e$, $\{X_e^i\}_i$ is a basis of $\mfg$,
and for each $x$, $\{v_x^j\}_j$  is a basis of $T_{\Phi_x} M$.


Given $f \in C^{\infty}(\mathcal{Q}_\Lambda)$,
letting $v^i_e, v^j_x$ with $e,x,i,j$  as above
 be an orthonormal basis of  $T_{(Q,\Phi)}\cQ_\Lambda$,
we can write the gradient and Laplace--Beltrami operator at edges $e$ or vertices $x$ as
\begin{equs}[e:D_ex]
	\nabla_{e}f\eqdef \sum_{i=1}^{\dg}(v_e^if)v_e^i,
	\qquad
	\Delta_ef\eqdef  \div\nabla_ef=\sum_{i=1}^{\dg}\<\nabla_{v_e^i}\nabla_{e}f,v_e^i\>,
	\\
		\nabla_{x}f\eqdef \sum_{j=1}^{d_M}(v_x^jf)v_x^j,
	\qquad
	\Delta_xf\eqdef  \div\nabla_xf=\sum_{j=1}^{d_M}\<\nabla_{v_x^j}\nabla_{x}f,v_x^j\>.
\end{equs}
We then have
$\nabla f
=\sum_{e\in E_{\Lambda}^+}\nabla_e f
+\sum_{x\in\Lambda} \nabla_x f$.

\medskip

{\bf Notations.} Throughout the paper:
\begin{itemize}
\item
We always write $\mu_\Lambda$ for the YMH measure  \eqref{measure}
on a finite lattice $\Lambda$. We will hide its dependence on the Higgs target space $M$ in our notation since it will be clear which $M$ is under consideration from the context. We always write $\mu$ for an infinite volume limit of $\mu_\Lambda$.
\item
In the case $M=\R^N$, when we do disintegration for $\mu_\Lambda$, we will write $\mu_Q$ for the measure over the fields $\Phi$ with the field $Q$ fixed, the $\nu$ for the remaining measure, see Section~\ref{sec:pR}. For simplicity we hide their dependence on $\Lambda$ in the notation.
\item
Only in Section~\ref{sec:gauge fixing} we will use the notation $\widetilde\nu_\Lambda$ for a gauge fixed measure.
\item
We use the notation $a\lesssim b$ if there exists a constant $c>0$ such that $a\leq cb$, and we write $a\simeq b$ if $a\lesssim b$ and $b\lesssim a$.
\end{itemize}

\section{Yang--Mills--Higgs SDEs}
\label{sec:YMHSDEs}

In this section we consider the Langevin dynamics \eqref{e:LSDE} for the measure \eqref{measure}. More precisely, on the lattice $\Lambda$, consider
the following SDE on $\cQ_L$
\begin{equ}[eq:YM1]
	\dif (Q,\Phi) = \nabla \mathcal S_\YMH (Q,\Phi) \dif t + \sqrt 2\dif \mathfrak B\;,
\end{equ}
with $\mathfrak B=(\mathfrak B_e, \mathfrak B_x : e\in E^+_{\Lambda},x\in \Lambda)$ being a collection of  independent Brownian motions,
where each $\mathfrak B_e$ takes values in
 $G$ and each $\mathfrak B_x$ takes values in $M\in \{\mR^N, \mS^{N-1},G\}$.

We first derive the more explicit forms of the equations, and address the questions of well-posedness and invariant measures.
Compared to the pure Yang--Mills case in \cite{SSZloop} and \cite{SZZ22}, the main difficulty lies in the control
in the case    where $\Phi$ takes values in the unbounded (non-compact) space $\mR^N$, and we will provide more details in this case in Section~\ref{sec:SDE-RN}. For the bounded cases $M \in \{\mS^{N-1},G\}$ where the target spaces are also compact,  the proof follows similarly as in \cite{SZZ22}, see Section~\ref{sec:SDESG}.

Recall from \eqref{e:CS} in Section~\ref{sec:not} that (possibly up to an additive constant which is irrelevant for the measure) we write $\mathcal S_{\YMH} = \cS_1-\cS_2$ where
\begin{equ}
	\cS_1\eqdef N\beta \, \Re \sum_{p\in \CP^+_\Lambda} \Tr(Q_p),
	\qquad
	\cS_2\eqdef
	\begin{cases}
	\kappa N\sum_{e\in E^+_\Lambda }
	|Q_e \Phi_y - \Phi_x|^2  +mN\sum_{z\in \Lambda}|\Phi_z|^2 &\; (M=\mR^N)
	\\
	-2 \kappa N\sum_{e\in E^+_\Lambda }\Phi_x^tQ_e \Phi_y &\; (M=\mS^{N-1})
	\\
	-2 \kappa N\sum_{e\in E^+_\Lambda }\tr(\Phi_x^tQ_e \Phi_y) &\; (M=G)
	\end{cases}
\end{equ}
with $e=(x,y)$.

We will give the explicit expression for $\nabla \mathcal S_{\YMH}$.
To this end, we introduce the following notation:
\begin{enumerate}
\item
For a plaquette $p=e_1e_2e_3e_4 \in \CP$, we write $p\succ e_1$ to indicate that $p$ is a plaquette that starts from edge $e_{1}$.
Note that for each edge $e$, there are $2(d-1)$ plaquettes in $\CP$
such that $p\succ e$.
\item  For $x,y\in \mZ^d$, we write $x\sim y$ if there exists one edge $e$ such that $x,y\in e$. In particular $x\sim x$. For edges $e$ and $\bar e$,  we write $\bar e \sim e$ which means that there exists $p\in \cP$ with $p\succ e$ and $\{\bar e,\bar e^{-1}\}\cap p \neq \emptyset$.
In particular $e\sim e$.
\end{enumerate}
Recall the notation \eqref{e:D_ex}.

\bl\label{lem:nab}
For each $e \in E_{\Lambda}^{+}$, and $M\in \{\mR^N,\mS^{N-1}, G\}$ we have
\begin{equs}
	\nabla_e\cS_1(Q,\Phi)&= -\frac12 N\beta \sum_{p\in \cP_{\Lambda}, p\succ e} (Q_p- Q_{p}^{t}) \cdot Q_{e}\;,
	\\
	\nabla_e\cS_2(Q,\Phi)&= \kappa N(Q_e\Phi_y\Phi_x^t-\Phi_x\Phi_y^tQ_e^t)Q_e\;, \quad e=(x,y).
\end{equs}
\el

%


\begin{proof}
The first identity is  \cite[Lemma 3.1]{SSZloop}. The second identity follows  in a similar way. More precisely,
recalling \eqref{e:h}, it suffices to differentiate $\Tr(\Phi_x^t Q_e \Phi_y)$ for all the choices of $M$.
We denote by $\proj$ the orthogonal projection of $M_N$  onto the Lie algebra $\mfg$. Let $X_e\in \mfg$.
	We compute the derivative at $Q$ in the tangent direction $X_eQ_e$ (see Section~\ref{sec:mani} for background on derivatives)
	$$
	\frac{\dif}{\dif t} \Big|_{t=0} \Tr\Big(\Phi_x^t e^{tX_e}Q_e\Phi_y\Big)
	=\Tr \Big(\Phi_x^t X_eQ_e\Phi_y \Big)
	= \Tr\Big(X_eQ_e \Phi_y \Phi_x^t\Big)
	=\Tr\Big(X_e \proj(Q_e \Phi_y \Phi_x^t)\Big).
	$$
	For this to match with $\Tr (XQ_e D^t)$ where $D=\nabla_e \Tr(\Phi_x^t Q_e \Phi_y)$, we see that
	$$
	D=(Q_e^t \proj (Q_e \Phi_y \Phi_x^t))^t = -\frac12 (Q_e \Phi_y \Phi_x^t - \Phi_x \Phi_y^t Q_e^t )Q_e
	$$
	which yields the second identity in the lemma.
\end{proof}

Recall our notation $e=(x,y)\in E_\Lambda$ which means that
$x=u(e)$ and $y=v(e)$.

\bl\label{lem:nab_x}
For each 
$x\in \Lambda$,  we have
\begin{equ}
	\nabla_x \cS_{\YMH}(Q,\Phi)=
	\begin{cases}
	2\kappa N
	\sum_{e=(x,y)\in E_\Lambda}  
	(Q_e\Phi_{y}-\Phi_x)-2m N\Phi_x  &\qquad (M=\mR^N),
	\\
	2\kappa N\sum_{e=(x,y)\in E_\Lambda}\Big(Q_e\Phi_{y}- (\Phi_x^t Q_e \Phi_y)  \Phi_x\Big) &\qquad (M=\mS^{N-1}),
	\\
	\kappa N\sum_{e=(x,y)\in E_\Lambda}\big(Q_e\Phi_{y}-\Phi_x\Phi_y^tQ_e^t\Phi_x\big) &\qquad (M=G).
	\end{cases}
\end{equ}
\el
\begin{proof}
Consider the case $M=\mR^N$.
Since $\cS_1$ is independent of $\Phi$, we only need to calculate $\nabla_x\cS_2$. For fixed $x$,
	\begin{align*}
		\nabla_x\Big(m\sum_{z\in\Lambda}|\Phi_z|^2\Big)=2m\Phi_x.
	\end{align*}
Concerning  $\sum_{e=(z,y)\in E_\Lambda^+}|Q_e \Phi_y - \Phi_z|^2$, the terms which contribute to $\nabla_x\cS_2$ for fixed $x$
can be written as
$\sum_{e=(x,y)\in E_\Lambda}|Q_e \Phi_y - \Phi_x|^2$.
Here, we used the fact that  for an edge $e=(y,x)\in E_\Lambda^+$,
we can reverse its orientation so that $x$ becomes the beginning point, while the summand is invariant (see \eqref{e:rev-e}). 	
By \eqref{e:h}, we have
		\begin{align*}
			\nabla_x|Q_e \Phi_y - \Phi_x|^2=-2(Q_e\Phi_y-\Phi_x),
		\end{align*}
so the first identity follows.
%
For the case $M=\mS^{N-1}$, we first compute the gradient in $\R^{N}$ as above, which gives $2\kappa N\sum_{e=(x,y)} Q_e\Phi_{y}$.
We then project $Q_e\Phi_y$ onto $T_{\Phi_x}\mS^{N-1}$, which is given by
	$$
 Q_e \Phi_y - \frac{\<  Q_e \Phi_y, \Phi_x\>}{|\Phi_x|^2 } \Phi_x
	=
  Q_e \Phi_y -  (\Phi_x^t Q_e \Phi_y)  \Phi_x$$
so the second identity follows. The third identity follows in the same way.
\end{proof}

\br
For each $x$,
the term $\sum_{e=(x,y)} (Q_{e}\Phi_{y}-\Phi_{x})$
appearing in Lemma~\ref{lem:nab_x} can be then viewed as a discrete divergence of
the discrete covariant derivative \eqref{e:cov-der}, namely,  a ``discrete gauge covariant Laplacian'' at $x$. This is a standard discretization of the gauge covariant Laplacian, see also e.g.  \cite[Section~3.2]{MR2884879}, \cite[Eq.~(20)]{MR4408017}.
\er

%

\subsection{YMH SDE for $M=\mR^N$} 
\label{sec:SDE-RN}


By Lemmas~\ref{lem:nab} -- \ref{lem:nab_x}  and \eqref{e:c_g} we write  the SDE \eqref{eq:YM1} as follows: for $e\in E_{\Lambda}^+$, $x\in \Lambda$,
\begin{equation}\label{SDE:RN}
	\begin{dcases}
	\dif Q_e=-\frac12N\beta\sum_{p\in \cP_\Lambda,p\succ e}(Q_p-Q_p^t)Q_e\dif t
	-\kappa N\big(Q_e(\Phi_y\Phi_x^t)Q_e-\Phi_x\Phi_y^t\big)\dif t
	\\\qquad\qquad-\frac12(N-1)Q_e\dif t+\sqrt 2\dif B_eQ_e,
	\\
	\dif \Phi_x=2\kappa N\sum_{e=(x,y)\in E_\Lambda}(Q_e\Phi_{y}-\Phi_x)\dif t-2mN\Phi_x\dif t+\sqrt 2\dif B_x,
	\end{dcases}
\end{equation}
where  $(B_e)_{e\in E_\Lambda^+}$ is a collection of  $\mfg$-valued Brownian motions  and $\{B_x\}_{x\in \Lambda}$ is a collection of  $\mR^N$-valued Brownian motions, all independent.
Remark that we have the Stratonovich correction term of $Q_e$ and $B_e$, which is $\frac12\cdot \sqrt 2\<\dif Q_e,\dif B_e\>=\frac12\cdot  2\<\dif B_e Q_e,\dif B_e\>$ giving the term $\frac12(N-1)Q_e\dif t$.

We first prove global well-posedness
on a fixed finite lattice $\Lambda=\Lambda_L$,
given by the following lemma.
Remark that 
this lemma (as well as Proposition~\ref{lem:4.7} below)
will not be directly invoked in the rest of the paper,
but these global well-posedness results are of independent interest
and of theoretical importance, so we give the proofs.

\bl\label{lem:exist}
For fixed $N\in\mN, \beta, m, \kappa\in\mR$  and $\Lambda$, given any initial data $(Q(0),\Phi(0))\in \cQ_\Lambda$ 
there exists a  unique   solution $(Q,\Phi)\in C([0,\infty);\cQ_\Lambda)$
to \eqref{SDE:RN} such that for any $p\geq1, T>0$
	\begin{align}\label{bdf:Phi}
	\E\Big(	\sup_{s\in [0,T]}\sum_{x\in \Lambda}|\Phi_x(s)|^{2p}\Big)\leq C_{N,m,p,d,T,\kappa,\Lambda}\Big(1+\sum_{x\in\Lambda}\|\Phi_x(0)\|^{2p}\Big)
\end{align}
for a constant $C_{N,m,p,d,T,\kappa,\Lambda}$ depending on $N, m, p, d, T, \kappa$ and $|\Lambda|$.
\el
\begin{proof}
	For fixed $N$ and $\Lambda$,  \eqref{SDE:RN} is a finite dimensional SDE with locally Lipschitz coefficients. We introduce  stopping time $\tau_R=\tau_0\wedge \tau_R^0$ for $R>0$ with
	\begin{align*}
		\tau_0:=&\inf \{t\geq 0: \|Q_e(t)\|_{\infty}>2, \, \, \text{for at least one }  \, e\in E_{\Lambda}\}\wedge T,
		\\\tau_R^0:=&\inf \{t\geq 0: \|\Phi_x(t)\|_{\infty}>R, \, \, \text{for at least one }  \, x\in \Lambda\}\wedge T,
	\end{align*}
where $\|Q_e\|_\infty\eqdef \max_{i,j}|Q_e^{ij}|$ and $\|\Phi_x\|_{\infty}\eqdef \max_{i}|\Phi_x^i|$.
We then have unique local solution $(Q,\Phi)=((Q_e)_{e\in E_\Lambda^+},(\Phi_x)_{x\in \Lambda})$ with $Q_e\in C([0,\tau_R];M_N)$ for $e\in E_\Lambda^+$ and $\Phi_x\in C([0,\tau_R];\mR^N)$ for $x\in \Lambda$,  which satisfies \eqref{SDE:RN} before $\tau_R$, by a usual cut-off technique (see e.g. \cite[Theorem 3.2]{MR3170231}).

	Since for  $e\in E_\Lambda^+$,
	$\nabla \cS_{\YMH}(Q)_e$ belongs to $T_{Q_e}G$, 
	exactly the same argument as in the proof of \cite[Lemma 3.2]{SSZloop} implies that $Q_e(t)\in G$, $\forall t\geq0$,  and $\tau_0=T$ a.s..
	
	We then turn to global bound on $\Phi$. Applying It\^o's formula for $|\Phi_x|^{2p}$, $p\geq1$ we obtain
	\begin{equs}[ep:PhiL]
		\dif |\Phi_x|^{2p} 
		&\leq 
		2p(2p-1)N |\Phi_x|^{2p-2}\dif t
		+4p|\kappa| N\sum_{e=(x,y)}
		|\Phi_x|^{2p-2}\big(|\Phi_y||\Phi_x|+|\Phi_x|^2\big)\dif t
		\\&\quad
		+4|m|p N|\Phi_x|^{2p}\dif t
		+2\sqrt 2p\,\big\<|\Phi_x|^{2p-2}\Phi_x,\dif B_x\big\>.
	\end{equs}
Regarding  the martingale part 
$$
M_t:=2\sqrt 2p\int_0^t \big\<|\Phi_x|^{2p-2}\Phi_x,\dif B_x\big\>
$$
by  Burkholder--Davis--Gundy inequality,
followed by Cauchy-Schwarz, we obtain
\begin{equs}[eq:BDG]
\E &\sup_{s\in [0,t\wedge \tau_R]}|M_s|
\leq 
C\E\Big(\int_0^{t\wedge \tau_R}|\Phi_x|^{4p-2}\dif s\Big)^{1/2}
\\
&\leq 
C\E \Big[\sup_{s\in [0,t\wedge \tau_R]}|\Phi_x(s)|^{p}
\Big(\int_0^{t\wedge \tau_R}|\Phi_x|^{2p-2}\dif s\Big)^{1/2}\Big]
\\
&\leq
\frac12\E\Big(	\sup_{s\in [0,t\wedge\tau_R]}|\Phi_x(s)|^{2p}\Big)+C\E\Big(\int_0^{t\wedge \tau_R}|\Phi_x|^{2p-2}\dif s\Big).
\end{equs}
Summing over $x$ for \eqref{ep:PhiL}, 
noting that $|\Phi_x|^{2p-1} |\Phi_y| \leq C_p( |\Phi_x|^{2p} + |\Phi_y|^{2p})$
and $\sum_x  \sum_{e=(x,y)}|\Phi_y|^{2p}=2d \sum_x |\Phi_x|^{2p}$,
we obtain
	\begin{align*}
\E\Big(	\sup_{s\in [0,t\wedge\tau_R]}\sum_{x\in \Lambda}|\Phi_x(s)|^{2p}\Big)\leq C_{N,m,p,d,T,\kappa,\Lambda}+C\sum_{x\in \Lambda}|\Phi_x(0)|^{2p}+C_{m,p}\E\int_0^{t\wedge \tau_R}\sum_{x\in \Lambda}|\Phi_x(s)|^{2p}\dif s.
\end{align*}
Using Gronwall's inequality we find
	\begin{align*}
	\E\Big(	\sup_{s\in [0,t\wedge\tau_R]}\sum_{x\in \Lambda}|\Phi_x(s)|^{2p}\Big)\leq C_{N,m,p,d,T,\kappa,\Lambda}\Big(1+\sum_{x\in\Lambda}\|\Phi_x(0)\|^{2p}\Big),
\end{align*}
with $C_{N,m,p,d,T,\kappa,\Lambda}$  independent of $R$. Letting $R\to \infty$, we have $\tau_R\to T$ a.s.. Hence, we can extend the local solutions to global solutions on $[0,T]$ and have \eqref{bdf:Phi}.
Since $T$ is arbitrary the result follows.
\end{proof}

It is obvious that the measure \eqref{measure} is well-defined for $m>0$ and $\kappa>0$.
When we discuss the measure below, we always assume $m>0$ and $\kappa>0$.
The same argument as in \cite[Lemma 3.3]{SSZloop} implies the following result.

\begin{lemma}\label{lem:inv}
For $m>0$, $\kappa>0$,	\eqref{measure} is invariant under the SDE system \eqref{SDE:RN}.
\end{lemma}

By Lemma~\ref{lem:exist}, the solutions to \eqref{SDE:RN}
form a Markov process in $\cQ_\Lambda$.
We write $(P_t^\Lambda)_{t\geq0}$ for the associated semigroup,
i.e. for $f\in C^\infty(\cQ_\Lambda)$, $(P_t^\Lambda f)(z)=\E f(Q^z(t),\Phi^z(t))$ for $z\in \cQ_\Lambda$, where $(Q^z(t),\Phi^z(t))$ denotes the solution at time $t$ to \eqref{SDE:RN} starting from $z\in \cQ_\Lambda$.
We can also write down the Dirichlet form associated with $(P_t^\Lambda)_{t\geq0}$ which will be used later when we study functional inequalities. More precisely, for $F\in C^\infty(\cQ_\Lambda)$ we consider the following  symmetric quadratic form
\begin{equs}[def:EL]
	\cE^\Lambda(F,F)
	&\eqdef\int \<\nabla F,\nabla F\>_{T_{(Q,\Phi)}\cQ_\Lambda}\dif \mu_{\Lambda}.
\end{equs}
Then it is easy to see that
\begin{align*}
	\cE^\Lambda(F,F)
	&=\sum_{e\in E_{\Lambda}^+}\int\<\nabla_{e}F,\nabla_{e} F\>\dif \mu_{\Lambda}+\sum_{x\in \Lambda}\int\<\nabla_{x}F,\nabla_{x} F\>\dif \mu_{\Lambda}
	\\&=\sum_{e\in E_{\Lambda}^+}\int\tr(\nabla_{e}F(\nabla_{e} F)^t)\dif \mu_{\Lambda}+\sum_{x\in \Lambda}\int|\nabla_{x}F|^2\dif \mu_{\Lambda}.
\end{align*}
Using integration by parts formula for the Haar measure,
we have that $(\cE^\Lambda, C^\infty(\cQ_\Lambda))$ is closable,  and its closure $(\cE^\Lambda,D(\cE^\Lambda))$ is a regular Dirichlet form on $L^2(\cQ_\Lambda,\mu_{\Lambda})$. (c.f. \cite{Fukushima}.)

\medskip

In the following we extend the SDE \eqref{SDE:RN} from $\Lambda$ to the entire space $\mZ^d$.

Consider the SDE \eqref{SDE:RN} on the entire space, namely
 for $e\in E^+$, $x\in \mZ^d$,
\begin{equation}\label{SDE:RNi}
	\begin{dcases}
		\dif Q_e=-\frac12N\beta\sum_{p\in \cP,p\succ e}(Q_p-Q_p^t)Q_e\dif t
		-\kappa N\big(Q_e(\Phi_y\Phi_x^t)Q_e-\Phi_x\Phi_y^t\big)\dif t
		\\\qquad\qquad-\frac12(N-1)Q_e\dif t+\sqrt 2\dif B_eQ_e,
		\\
		\dif \Phi_x=2\kappa N\sum_{e=(x,y)}(Q_e\Phi_{y}-\Phi_x)\dif t-2m N\Phi_x\dif t+\sqrt2\dif B_x,
	\end{dcases}
\end{equation}
where  $(B_e)_{e\in E^+}$ are  $\mfg$-valued Brownian motions and $(B_x)_{x\in \mZ^d}$ are  $\mR^N$-valued Brownian motions, all independent.

Define configuration space  $\cQ=G^{E^+}\times \ell^p_k(\mZ^d;\mR^N)$ for some even $p>4$ and $k>d$, where
\begin{equ}[e:lpk]
\ell^p_k\eqdef
\ell^p_k(\mZ^d;\mR^N)\eqdef
\Big\{f:\mZ^d\to \mR^N, \quad \sum_x \frac1{|x|^k+1}|f_x|^p<\infty\Big\}
\end{equ}
is a weighted $\ell^p$ space with polynomial weight of degree $k$.
It is easy to check that $\ell^p_k$ is compactly embedded into $\ell^p_a$ for $a>k$.
Write $M_N^{E^+} = \prod_{e\in E^+} M_N$ for the direct product of (infinitely many) vector spaces $M_N$.
Define  a norm  on $M_N^{E^+}\times\ell^p_k(\mZ^d;\mR^N) $  by
\begin{equs}[norm]
	\|(Q,\Phi)\| &\eqdef \|Q\|+\|\Phi\|,\\ \|Q\|^2\eqdef\sum_{e\in E^+}\frac1{|e|^k+1}|Q_e|^2
	&,\quad \|\Phi\|^p\eqdef \sum_{x\in \mZ^d}\frac1{|x|^k+1}|\Phi_x|^p,
\end{equs}
with $|Q_e|^2=\<Q_e,Q_e\>=\tr(Q_eQ_e^t)$   and $|e|$ (resp. $|x|$) given by the distance from $0$ to $e$ (resp. $x$) in $\mZ^d$, respectively.
(More precisely, $|e|$ is the minimum of the distances from the two vertices of $e$ to $0$.) Now we give existence and uniqueness of solutions to \eqref{SDE:RNi}.

\bp\label{lem:4.7}
Fix $N\in\mN, \beta\in\mR, m\in \mR$. For any initial data $(Q^0,\Phi^0)\in \cQ$, there exists a unique probabilistically strong solution $(Q,\Phi)$ to \eqref{SDE:RNi} in $C([0,\infty);\cQ)$
satisfying for any $T>0$
\begin{equ}[e:mom-Phi]
	\E\sup_{s\in[0,T]}\|\Phi(s)\|^p+	\E\Big(\sup_{s\neq t\in [0,T]}\frac{\|\Phi(t)-\Phi(s)\|}{|t-s|^\alpha}\Big)<\infty.
\end{equ}
\ep

\begin{proof} 
For the pure Yang--Mills model on the entire lattice, such a well-posedness result was given in \cite[Proposition 3.4]{SZZ22}. The main difference here is that $\cQ$ is not compact and we need uniform bounds for the component $\Phi$.

We start with solutions on $\Lambda_L$.
More precisely, for every initial data $X^0\eqdef (Q^0,\Phi^0)\in \cQ$   we can easily find
$X^L(0)\in \cQ_{\Lambda_L}$,
with  periodic extension to the entire lattice still denoted by
$X^L(0)$,
such that $\|X^L(0)-X^0\|\to0$ as $L\to\infty$.
By Lemma \ref{lem:exist} we obtain a unique solution
$X^L\eqdef (Q^L,\Phi^L)\in C([0,\infty); \cQ_{\Lambda_L})$
to \eqref{SDE:RN} from $X^L(0)$. 
We can also extend $X^L$ to $\cQ$ by periodic extension.

We first show tightness.	Applying It\^o's formula to $|\Phi_x^L|^{p}$, we have \eqref{ep:PhiL} as before with $2p$ there replaced by $p$.
Multiplying both sides of \eqref{ep:PhiL} by $\frac1{|x|^k+1}$ with $k>d$, summing over $x$, applying Young's inequality and the Burkholder--Davis--Gundy inequality for the martingale part as in \eqref{eq:BDG}, and using $\frac1{|x|^k+1}\leq C\frac1{|y|^k+1}$ for $x\sim y$ with $C$ independent of $x,y$,  we obtain
\begin{align*}
		\E\sup_{s\in[0,t]}\|\Phi^L(s)\|^p\leq C\|\Phi^L(0)\|^{p}+ C_{p,N}\int_0^t\E\Big(1+\|\Phi^L(s)\|^{p}\Big)\dif s,
\end{align*}
which by Gronwall's inequality implies 
\begin{align}\label{PhiL}
	\E\sup_{s\in[0,T]}\|\Phi^L(s)\|^{p}\leq C\|\Phi^L(0)\|^{p}+ C_{p,N,T}.
\end{align}
Here all the constants are independent of $L$.
Furthermore, using SDEs \eqref{SDE:RN}  we also have
\begin{align*}
	\E|\Phi^L_x(t)-\Phi^L_x(s)|^p &\leq C_p\E\Big(\int_s^t \big(\sum_{y\sim x}|\Phi_y^L| \big)\dif r\Big)^p+C_p\E|B_x(t)-B_x(s)|^p
	\\&\leq C_p|t-s|^{p-1}\E\int_s^t \big(\sum_{y\sim x}|\Phi_y^L|^p \big)\dif r+C_p|t-s|^{p/2}.
\end{align*}
Hence, multiplying both sides by $\frac1{|x|^k+1}$ with $k>d$ and summing over $x$, we have, by \eqref{PhiL},
\begin{align*}
	\E\|\Phi^L(t)-\Phi^L(s)\|^p\leq C_{p}(|t-s|^p+|t-s|^{p/2}).
\end{align*}
As a consequence, by Kolmogorov's criterion we have for $\alpha<1/2$
\begin{align}\label{holderphi}
	\sup_L\E\Big(\sup_{s\neq t\in [0,T]}\frac{\|\Phi^L(t)-\Phi^L(s)\|}{|t-s|^\alpha}\Big)<\infty,
\end{align}
which combined with \eqref{PhiL} implies that $\{\Phi^L\}_L$ is tight in $C([0,T];\ell^p_a)$ with $a>k$.
Since $Q^L(t)$ takes values in a compact space  $G^{E^+}$,
 the marginal laws of $\{Q^L\}$ at each $t\ge 0$ form a tight set in $G^{E^+}$. Using \eqref{PhiL} and $\frac1{|e|^k+1}\simeq \frac1{|x|^k+1}$ for $x\in e$ we get
\begin{align*}
		\E\|Q^L(t)-Q^L(s)\|^{p/2}\leq C_{p}(|t-s|^{p/2}+|t-s|^{p/4}),
\end{align*}
which  again
by Kolmogorov criterion implies \eqref{holderphi}
with $\Phi^L$ replaced by $Q^L$.
Hence, $X^L=(Q^L,\Phi^L)$ is tight in $C([0,\infty);G^{E^+}\times \ell^p_a)$ equipped with the distance
\begin{align*}
	\widetilde\rho(X,X')\eqdef
	\sum_{n=0}^\infty 2^{-n}\Big(1\wedge \sup_{t\in [n,n+1]}(\|X(t)-X'(t)\|\Big),
\end{align*}
for  $X,X'\in C([0,\infty);G^{E^+}\times \ell^p_a)$.
Since $(C([0,\infty);G^{E^+}\times \ell^p_a), 	\widetilde\rho)$ is a Polish space, existence of probabilistically weak solutions $(Q,\Phi)$ follows from the usual Skorokhod Theorem and
taking limit on both sides of the equation. We refer to the proof of \cite[Proposition 3.4]{SZZ22} for more details. Using \eqref{PhiL} and \eqref{holderphi} and by lower semicontinuity,
  we conclude  that \eqref{e:mom-Phi} holds for every $T>0$, $\alpha\in (0,1/2)$
which implies that $(Q,\Phi)\in C([0,\infty),\cQ)$.

Now we prove pathwise uniqueness.
Consider two solutions  $(Q,\Phi), (Q',\Phi')\in C([0,T];\cQ)$ starting from the same initial data $(Q(0),\Phi(0))\in \cQ$. Set
$\delta Q_e:=Q_e-Q'_e$ and $\delta \Phi_x:=\Phi_x-\Phi'_x$, and choose the following weight
$$\rho_x(t)=e^{-t(|x|^k+1)}.$$
Define stopping time
\begin{align*}
	\tau_R\eqdef \inf\{t\geq0; \|\Phi(t)\|+\|\Phi'(t)\|\geq R\}\wedge T.
\end{align*}
Since $Q_e, Q_e'\in G$, 
by It\^o's formula for $\rho_x(t)|\delta Q_e(t)|^2$, $\rho_x(t)|\delta \Phi_x(t)|^2$ and the Burkholder--Davis--Gundy inequality  for stochastic integrals, we obtain for each $e=(x,y)$
\begin{align*}
&\E  \sup_{t\in [0,T\wedge \tau_R]}\rho_x|\delta Q_e|^2
\\
&\leq C_{N}\E\int_0^{T\wedge \tau_R}(1+|\Phi_x||\Phi_y|)|\delta Q_e|^2\rho_x\dif s
+C_N \sum_{\bar e\sim e}   
\E\int_0^{T\wedge \tau_R}|\delta Q_e||\delta Q_{\bar e}|\rho_x\dif s
	\\&\quad +C_{\kappa,N}\E\int_0^{T\wedge\tau_R} (|\Phi_y||\delta\Phi_x|+|\Phi_x'||\delta \Phi_y|)|\delta Q_e|\rho_x\dif s-\E\int_0^{T\wedge \tau_R}(1+|x|^k)|\delta Q_e|^2\rho_x\dif s,
\end{align*}
and for each $x$
\begin{align*}
	\E\sup_{t\in[0,T\wedge \tau_R]}\rho_x|\delta \Phi_x|^2
	&\leq C_{\kappa,N}\sum_{e=(x,y)}\E\int_0^{T\wedge \tau_R}(|\Phi_y||\delta Q_e|+|\delta \Phi_y|+|\delta \Phi_x|)|\delta \Phi_x|\rho_x\dif s
	\\&
	\quad-(1+|x|^k)\E
	\int_0^{T\wedge \tau_R}|\delta \Phi_x|^2\rho_x\dif s.
\end{align*}
Here the last terms with the factor $(1+|x|^k)$ arise from
 the time derivative of the weight $\rho$.
Since $|\Phi_x|+|\Phi_y|\leq C_R (1+|x|^{k/p})$ for $x\sim y$ before $\tau_R$ and $|x|^{2k/p}\leq \eps |x|^k+C_\eps$ for $\eps>0$ and $p>4$, we obtain
$|\Phi_x|^2+|\Phi_y|^2\leq C_{\eps,R} +\eps|x|^{k}$ before $\tau_R$.
Thus for $\eps>0$ small we can use the last terms with the factor $(1+|x|^k)$
 in the above two inequalities to absorb the terms involving $|\Phi_x|, |\Phi_y|, |\Phi_x'|$ and obtain
\begin{align*}
	\E  \sup_{t\in [0,T\wedge \tau_R]}\rho_x|\delta Q_e|^2
	&
	\leq C_{N,R}\E\int_0^{T\wedge \tau_R}|\delta Q_e|^2\rho_x\dif s
	+C_N \sum_{\bar e\sim e}
	\E\int_0^{T\wedge \tau_R}|\delta Q_e||\delta Q_{\bar e}|\rho_x\dif s
	\\&\quad+C_{N}\E\int_0^{T\wedge \tau_R} (|\delta\Phi_x|^2+|\delta \Phi_y|^2)\rho_x\dif s,
\end{align*}
and
\begin{align*}
\E\sup_{t\in[0,T\wedge \tau_R]}\rho_x|\delta \Phi_x|^2
\leq C_N\sum_{e=(x,y)}
\E\int_0^{T\wedge \tau_R}\Big[|\delta Q_e|^2+(|\delta \Phi_y|+|\delta \Phi_x|)|\delta \Phi_x|\Big]\rho_x\dif s.
\end{align*}
We write $\rho_e:=\rho_x$ for $x=u(e)$.
Since $|\rho_x/\rho_y|\lesssim 1$ and $|\rho_e/\rho_{\bar e}|\lesssim 1$ for $x\sim y$ and $e\sim \bar e$, we can change the weight $\rho_x$, $\rho_e$ to $\rho_y$, $\rho_{\bar e}$, respectively and
 obtain
\begin{equs}
\E\sup_{t\in [0,T\wedge \tau_R]}|\delta Q_e|^2\rho_e
&\leq
C_{N,R}
\E\int_0^{T\wedge \tau_R}
\Big(\sum_{x\in e} |\delta \Phi_x|^2\rho_x
+
 \sum_{\bar e\sim e} 
 |\delta Q_{\bar e}|^2\rho_{\bar e}\Big)\dif s,
\\
	\E\sup_{t\in [0,T\wedge \tau_R]}|\delta \Phi_x|^2\rho_x
	&\leq
	C_{N,R}
	\E\int_0^{T\wedge \tau_R}
	\Big(\sum_{e\ni x}|\delta Q_e|^2\rho_e
	+\sum_{y\sim x}|\delta \Phi_y|^2\rho_{y}\Big)\dif s.
\end{equs}
Summing over $e$ and $x$ we get
\begin{align*}
	\E\sup_{t\in [0,T\wedge \tau_R]}\|\delta Q\|^2_\rho+\E\sup_{t\in [0,T\wedge \tau_R]}\|\delta \Phi\|^2_\rho\leq C_{N,R} \E\int_0^{T\wedge \tau_R}(\|\delta Q\|_\rho^2+\|\delta \Phi\|_\rho^2)\dif s,
\end{align*}
where
\begin{align*}
 \|Q\|^2_\rho\eqdef\sum_{e\in E^+}\rho_e|Q_e|^2
	&,\quad \|\Phi\|^2_\rho\eqdef \sum_{x\in \mZ^d}\rho_x|\Phi_x|^p.
\end{align*}
Hence, pathwise uniqueness follows by Gronwall's lemma and sending $R\to\infty$.
By Yamada--Watanabe Theorem \cite{Kurtz},
weak existence and pathwise uniqueness gives us  existence and uniqueness of probabilistically strong solution.
\end{proof}

By Proposition \ref{lem:4.7}, the solutions to \eqref{SDE:RNi} form a Markov process in $\cQ$.
We denote by $(P_t)_{t\geq0}$ the associated semigroup. Using the dynamics, we also have the following tightness result.

\bl
\label{lem:tightR}
For $m>0$ and $\kappa>0$, the fields  $(\mu_{\Lambda_{L}})_L$ in $\cQ$ form a tight set with respect to the topology induced by the norm in \eqref{norm}.
\el
\begin{proof}
	Consider the solutions $(Q^L,\Phi^L)$ to SDEs \eqref{SDE:RN} starting from $(\mu_{\Lambda_{L}})_L$. By Lemma \ref{lem:inv} $(Q^L,\Phi^L)$ is stationary solutions to \eqref{SDE:RN}. Furthermore, applying It\^o's formula to $|\Phi_x^L(t)|^{p}$ as in  \eqref{ep:PhiL}, we obtain
	\begin{equs}\label{bdpPhi}
	&\E |\Phi_x^L(t)|^{p}-\E|\Phi_x^L(0)|^{p}
	\\
	&\leq 2p\kappa N\sum_{e=(x,y)} t\, \E|\Phi_x^L|^{p-2}\big(|\Phi_y^L||\Phi_x^L|-|\Phi_x^L|^2\big)
	-2mNpt\,\E|\Phi_x^L|^{p}
		+p(p-1)Nt\,\E|\Phi_x^L|^{p-2}.
	\end{equs}
Since $\Phi_x^L$ is stationary, the LHS is zero. By translation invariance  $\E|\Phi_x|^p=\E|\Phi_y|^p$, which substituted into \eqref{bdpPhi} implies
	\begin{equs}
	&\E|\Phi_x^L|^{p}
	\leq C_{p,m}\E|\Phi_x^L|^{p-2}\leq \eps \E|\Phi_x^L|^p+C_\eps.
\end{equs}
Hence, $\E|\Phi_x^L|^{p}\lesssim 1$
uniformly in $x$ and $L$.
Therefore $\E\|\Phi^L\|_{\ell^p_{k}}^p\lesssim 1$
for $k>d$.
 By compact embedding $\ell^p_k\subset \ell^p_{k'}$ for $k<k'$, the result follows.
\end{proof}

 Since by Lemma \ref{lem:inv} $\mu_{\Lambda_{L}}$ is an invariant measure for \eqref{SDE:RN}, we then obtain the following result by the same argument as in the proof of \cite[Theorem 3.5]{SZZ22}.

 \bt\label{th:in1} Every tight limit $\mu$ of $\{\mu_{\Lambda_{L}}\}$ is an invariant measure for \eqref{SDE:RNi}.
 \et

The dynamics  \eqref{SDE:RNi} is gauge covariant.
Namely, for every $G$-valued function $g$ on $\mZ^d$,
define $g\circ (Q,\Phi)$ by \eqref{e:gauge}.
%
If $(Q,\Phi)$ is a solution to \eqref{SDE:RNi}, it is easy to check that 
$g\circ(Q,\Phi)$
also satisfies \eqref{SDE:RNi} with $B_e$ and $B_x$ replaced by $g_xB_eg_x^{-1}$ and $g_xB_x$, which are still  Brownian motion in $\mfg$ and $\mR^N$, respectively.  By Proposition \ref{lem:4.7}, the uniqueness in law to SDE \eqref{SDE:RNi} holds. Hence, we obtain the following result.

\bp
Fix $N\in\mN, \beta, m\in \mR$. Let $(Q,\Phi)$ and $(\bar Q,\bar \Phi)$  be the unique solutions to \eqref{SDE:RNi}
with initial conditions $(Q^0,\Phi^0)$ and $(\bar Q^0,\bar \Phi^0)$ in $\cQ$ respectively.
If $(\bar Q^0,\bar \Phi^0) = g\circ (Q^0,\Phi^0)$
 for some $G$ valued function $g$
on $\mZ^d$, then, $(\bar Q(t), \bar \Phi(t))$ and  $g\circ(Q(t),\Phi(t))$ are equal in law
for all $t\geq0$.
\ep

We introduce the Dirichlet form  associated with \eqref{SDE:RNi}.
Define the space of cylinder functions 
\begin{equs}[e:cylQ]
	C^\infty_{cyl}(\cQ)
	=\Big\{
	F: F&=f(Q_{e_1},\dots, Q_{e_n},\Phi_{x_1},\dots,\Phi_{x_k}),
	\\&n,k\in \mN, e_i\in E^+, x_j\in \mZ^d, f\in C^\infty(G^n\times \mR^{Nk})\Big\}.
\end{equs}
For every tight limit $\mu$ and $F\in C^\infty_{cyl}(\cQ)$ we define
the following  symmetric quadratic form
\begin{equs}[e:EFF]
	\cE^{\mu}(F,F)&\eqdef\sum_{e\in E^+}\int\<\nabla_{e}F,\nabla_{e} F\>\dif \mu+\sum_{x\in \mZ^d}\int\<\nabla_{x}F,\nabla_{x} F\>\dif \mu.
\end{equs}
The same argument as in \cite[Proposition 3.7]{SZZ22} implies that
$(\cE^{\mu},C^\infty_{cyl}(\cQ))$ is closable and its closure $(\cE^{\mu},D(\cE^{\mu}))$ is a  Dirichlet form on $L^2(\cQ,\mu)$.
Furthermore, since $\cQ$ is a Polish space,
$(\cE^{\mu},D(\cE^{\mu}))$ is a quasi-regular Dirichlet form  
(c.f.  \cite[IV. Sec. 4.b]{MaRockner}).

\subsection{YMH SDE for  $M\in \{ \mS^{N-1}, G\}$} 
\label{sec:SDESG}


%

In the case where $\Phi_x$ is $\mS^{N-1}$ valued,
the equation for $\Phi_x$ is driven by $\sqrt 2 \mathfrak{B}_x$
where $\mathfrak{B}_x$ is the Brownian motion on the sphere.  By \eqref{bm:sp}, we write the $i$-th component of $\sqrt 2\dif \mathfrak{B}_x$ as
	\begin{align*}
		&\sqrt2\sum_{j=1}^N(\delta_{ij}-\Phi_x^i\Phi_x^j)\circ \dif B_x^j=\sqrt2\sum_{j=1}^N(\delta_{ij}-\Phi_x^i\Phi_x^j)\dif B_x^j-\frac{\sqrt2}{2}\sum_{j=1}^N\<\dif (\Phi_x^i\Phi_x^j),\dif B^j_x\>
	\end{align*}
and the It\^o correction term (e.g. the last term)
is equal to
$$
-\sum_{j=1}^N(\delta_{ij}-\Phi_x^i\Phi_x^j)\Phi_x^j\dif t-\Phi_x^i\sum_{j=1}^N(\delta_{jj}-\Phi_x^j\Phi_x^j)\dif t
=-(N-1)\Phi_x^i\dif t,
$$
where we used $|\Phi_x|^2=1$.
%
Hence by  Lemmas~\ref{lem:nab} -- \ref{lem:nab_x}  and \eqref{e:c_g}  we write SDEs \eqref{eq:YM1} as follows: for $e\in E_{\Lambda}^+, x\in \Lambda$,
	\begin{equs}[SDES]
		\begin{dcases}
		\dif Q_e=-\frac12N\beta\sum_{p\in \cP_\Lambda,p\succ e}(Q_p-Q_p^t)Q_e\dif t-\kappa N (Q_e(\Phi_y\Phi_x^t)Q_e-\Phi_x\Phi_y^t)\dif t\\\qquad\qquad-\frac12(N-1)Q_e\dif t+\sqrt2\dif B_eQ_e,\\
		\dif \Phi_x
		=2\kappa N\sum_{e=(x,y)\in E_\Lambda}(Q_e\Phi_{y}-(\Phi_x^tQ_e\Phi_y)\Phi_x)\dif t +\sqrt 2(I-\Phi_x\otimes\Phi_x)\dif B_x
		 -(N-1)\Phi_x\dif t,
		\end{dcases}
	\end{equs}
with $B_e$ being $\mfg$-valued Brownian motion and $B_x$ being $\mR^N$-valued Brownian motion.


%
%

In the case where $\Phi$ is $G$ valued,
we write SDEs \eqref{eq:YM1} as follows: for $e\in E_{\Lambda}^+, x\in \Lambda$
\begin{equs}[SDEG]
	\begin{dcases}
	\dif Q_e=-\frac12N\beta\sum_{p\in \cP_\Lambda,p\succ e}(Q_p-Q_p^t)Q_e\dif t-\kappa N(Q_e(\Phi_y\Phi_x^t)Q_e-\Phi_x\Phi_y^t)\dif t\\\qquad\qquad-\frac12(N-1)Q_e\dif t+\sqrt2\dif B_eQ_e,\\
	\dif \Phi_x=\kappa N\sum_{e=(x,y)\in E_\Lambda}(Q_e\Phi_y-(\Phi_x\Phi_y^tQ_e^t)\Phi_x)\dif t+\sqrt2\dif B_x\Phi_x-\frac{N-1}{2}\Phi_x\dif t,
	\end{dcases}
\end{equs}
with $B_e$, $B_x$ being independent Brownian motions in $\mfg$.

\vspace{3ex}
Since $\nabla_x\cS$ and $\nabla_e\cS$ stay in the related tangent space and the noise part is given by the corresponding Brownian motions, Lemma \ref{lem:exist} and Lemma \ref{lem:inv} still hold for the above two cases by exactly the same argument as in \cite[Lemma 3.2, Lemma 3.3]{SSZloop}. Hence, the solutions form a Markov process in $\cQ_\Lambda$ and we also use $(P_t^\Lambda)_{t\geq0}$ for the associated semigroup. Furthermore, we could also extend SDEs \eqref{SDES} and \eqref{SDEG} to $\mZ^d$ and construct a Markov process in $\cQ$ as in Proposition \ref{lem:4.7}. Furthermore, every tight limit $\mu$ of $(\mu_{\Lambda_L})$ is an invariant measure of the infinite volume dynamics i.e. SDEs \eqref{SDES} and \eqref{SDEG} with  $e\in E_{\Lambda}^+, x\in \Lambda$ replaced by $e\in E^+, x\in \mZ^d$. We can also write the related Dirichlet form as \eqref{e:EFF}.

\section{Functional inequalities and exponential ergodicity}\label{sec:ergo}

In this section we investigate the long time behavior of the Markov dynamics constructed in Section~\ref{sec:YMHSDEs}. 
Our primary approach is to utilize functional inequalities, i.e.  the log-Sobolev and the Poincar\'e inequalities.

To state these functional inequalities in our setting we first define some notation.
For a probability measure $\mu$ on $\cQ_\Lambda$ we define the entropy functional for positive $f\in C^\infty(\cQ_\Lambda)$
 $$\text{ent}_\mu(f)=\int f\log f\dif \mu-\int f\dif \mu\log \int f\dif \mu.$$
The entropy functional for a measure on $\cQ$ and $f\in C^\infty_{cyl}(\cQ)$ is defined in the same way. We will write $\E_\mu$ for the expectation w.r.t. $\mu$.

As alluded around \eqref{e:strategy},
a natural approach to prove  functional inequalities would be
to verify the well-known Bakry--\'Emery condition.
Recall \cite{BE,MR3155209,GZ} that for a probability measure
$\dif \mu\propto \exp(\cS)\dif \sigma$  on  a Riemannian manifold with $\sigma$ being the volume measure,
if there exists a positive constant $K>0$ such that
\begin{equs}[con:B-E]\Ric-\Hess_{\cS}\geq K,\end{equs}
then $\mu$ satisfies the log-Sobolev inequality with constant $K>0$, i.e. $K\text{ent}_\mu(F^2)\leq \mu(|\nabla F|^2)$.
This approach
was demonstrated in \cite{SZZ22} for the pure Yang--Mills, to derive the log-Sobolev  and thus the Poincar\'e inequality.

However,
when $M=\mR^N$,
the unbounded nature of $\Phi$ poses challenges in controlling the Hessian of the term $\Phi_x^tQ_e\Phi_y$ showing up in $\cS_{\YMH}$ as required for the Bakry--\'Emery condition.
For example we  have a term of the form $\tr(\Phi_x^tX_e^2Q_e\Phi_y)$ when calculating $\Hess_{\cS}$, and we cannot bound it by the Ricci curvature. This motivates the strategy of  first integrating out $\Phi$. In light of this, in Section~\ref{sec:pR}, our strategy will be to employ a disintegration, and establish the functional inequalities for the conditional probability given $Q$ first, from which we deduce  the mass gap of the conditional probability. Subsequently, we  combine it with further decomposition of the domain $\Lambda$ to give a delicate control of the correlation functions, appearing in the proof of the Poincar\'e inequality.
For the cases $M\in \{ \mS^{N-1}, G\}$ which we consider in
Section~\ref{sec:LSB}, we can directly  verify  the Bakry--\'Emery condition, with some careful calculations of Hessians for the Higgs terms on
$\mS^{N-1}$ or  $G$. We also refer to the introduction for more explanation on the idea of the proof.

\subsection{Poincar\'e inequalities for $M=\mR^N$} 
\label{sec:pR}

In this section we suppose that $m,\kappa>0$.
Recall that
\begin{align*}
	\dif\mu_{\Lambda}(Q,\Phi)
	:= Z_{\Lambda}^{-1}
	\exp\Big(\mathcal S_{\YMH} \Big)
	\prod_{e\in E^+_\Lambda} \dif\sigma_N(Q_e)
	\prod_{z\in \Lambda} \dif\Phi_z\,
\end{align*}
with  $\mathcal S_{\YMH}=\cS_1-\cS_2$
for
$\cS_1=N\beta \, \Re \sum_{p\in \CP^+_\Lambda} \Tr(Q_p)$ and
\begin{align*}
	\cS_2
	=\kappa N\sum_{e\in E_\Lambda^+}(\Phi_x^t\Phi_x-2\Phi_x^tQ_e\Phi_y+\Phi_y^t\Phi_y)
	+mN\sum_x\Phi_x^t\Phi_x.
\end{align*}
Our strategy is to employ disintegration to express $\mu_\Lambda$ as
\begin{equ}[e:disinteg]
	\mu_\Lambda(F)=\int F(Q,\Phi)\mu_\Lambda(\dif \Phi,\dif Q)=\int \int F(Q,\Phi)\mu_Q(\dif \Phi)\nu(\dif Q)=\E_{\nu}(\E_{\mu_Q}(F)),
\end{equ}
with $\mu_Q$ being the regular conditional probability given $Q$.  We write $\mu_Q(\dif \Phi)$ as probability measure on $\mR^{N\Lambda}=\mR^{N|\Lambda|}$
\begin{align}\label{def:muQ}
	\mu_Q(\dif \Phi)=\frac1{Z_Q}	\exp(-\cS_2(Q,\Phi))\prod_{x}\dif \Phi_x,
\end{align}
where $Z_Q=\int\exp(-\cS_2(Q,\Phi))\prod_{x}\dif \Phi_x$.
We write $\nu$ for the following probability measure on $G^{E^+_{\Lambda}}$
\begin{align}\label{def:nu}
	\nu(\dif Q)=\frac1Z\exp\Big(\cS_1-V(Q)\Big)\prod_{e}\dif\sigma_N(Q_e),
\end{align}
where $\cS_1$ is the pure Yang--Mills action which is independent of $\Phi$, and
\begin{align}\label{def:VQ}
V(Q)=-\log\int\exp(-\cS_2(Q,\Phi))\prod_{x}\dif \Phi_x.
\end{align}
It is easy to find that $V$ is well-defined and is a smooth function on $G^{E_\Lambda^+}$.

\br\label{rem:strat1}
Our disintegration may have some similarity with 
the recent works \cite{MR3926125,RolandIsing,RolandPhi4},
but the difference is that we have two fields $Q,\Phi$ while the above works in some sense decompose the field into two scales.
\er

Before proceeding we first give a uniform estimate under the measure $\mu_Q$,
which we prove using stochastic differential equation techniques i.e. by exploiting the Langevin dynamics.
This estimate plays an important role in the sequel, such as in the proof of functional inequalities
(Lemma~\ref{log:1}, Theorem~\ref{th:poin}) and proving mass gap (Section~\ref{sec:mass-RN}). In fact, since $\mu_Q$ is stationary measure of the Langevin dynamics,
 and we can use moment estimates of the dynamics via It\^o's calculus to derive suitable moments bounds of $\mu_Q$.

\bl
\label{lem:b:phi}
For any $x\in \Lambda$, 
\begin{align}\label{b:phi}
	\E_{\mu_Q}(|\Phi_x|^2)\leq \frac{1}{2m}.
\end{align}
\el
\begin{proof}
		For fixed $Q$, consider the Langevin dynamics for
	the measure $\mu_Q$:
	\begin{equ}[SDE:Q]
		\dif \Phi_x
		=2\kappa N\sum_{e=(x,y)}(Q_e\Phi_{y}-\Phi_x)\dif t
		-2m N\Phi_x
		\dif t
		+\sqrt2\dif B_x.
	\end{equ}
	Similarly as in the proof of  Lemma \ref{lem:exist} and Lemma \ref{lem:inv} (it is even simpler here since \eqref{SDE:Q}
	is exactly the 2nd equation in \eqref{SDE:RN} and  with fixed $Q$
	it is linear), we easily see the global well-posedness of SDE \eqref{SDE:Q} and $\mu_Q$ is an invariant measure of it.  We then take $\Phi=(\Phi_x)$ as the stationary solutions starting from $\mu_Q$.
	Applying It\^o's formula to $|\Phi_x|^2$ we get
	\begin{align}
		\dif |\Phi_x|^2
		=2N\dif t
		+4\kappa N\sum_e(Q_e\Phi_y\cdot \Phi_x-|\Phi_x|^2)\dif t
		-4mN|\Phi_x|^2
		\dif t
		+\dif M.
	\end{align}
	Here $M$ is a martingale and the term $N\dif t$ arises from the quadratic variation of the $\R^N$ valued Brownian motion $B_x$.

	In stationarity setting, taking expectation $\E$ w.r.t. $\mu_Q$ to eliminate the martingale term, and then using translation invariance $\E|\Phi_x|^2=\E|\Phi_y|^2$ and $|Q_e\Phi_y\Phi_x|\leq |\Phi_y||\Phi_x|$, we obtain \eqref{b:phi}.
\end{proof}

We then prove the log-Sobolev and the Poincar\'e inequalities for $\mu_Q$ and $\nu$,
which is the content of Lemma~\ref{log:1} below.
To state the result, we introduce the following  symmetric quadratic form for $F\in C^\infty(\mR^{N\Lambda})$
\begin{align}\label{dir:muQ}
	\cE^{\mu_Q}(F,F)\eqdef\sum_{x\in \Lambda}\int\<\nabla_{x}F,\nabla_{x} F\>\dif \mu_{Q}.
\end{align}
Using integration by parts formula for $\mu_Q$ (which obviously holds since it is the Lebesgue measure multiplied with a smooth density),
we have that $(\cE^{\mu_Q}, C^\infty(\mR^{N\Lambda}))$ is closable,  and its closure $(\cE^{\mu_Q},D(\cE^{\mu_Q}))$ is a regular Dirichlet form on $L^2(\mR^{N\Lambda},\mu_{Q})$ (c.f. \cite{Fukushima}).
We also
consider the following  symmetric quadratic form for $F\in C^\infty(G^{E_\Lambda^+})$
\begin{align}\label{dir:nu}
	\cE^{\nu}(F,F)&\eqdef\sum_{e\in E_\Lambda^+}\int\<\nabla_{e}F,\nabla_{e} F\>\dif \nu.
\end{align}
By \eqref{def:VQ} and Lemma \ref{lem:nab},
\begin{align}\label{def:neV}
\nabla_eV=\kappa N\E_{\mu_Q}(Q_e\Phi_y\Phi_x^tQ_e-\Phi_x\Phi_y^t).
\end{align}
This together with Lemma \ref{lem:nab} and the bound \eqref{b:phi} implies that
  $\nabla_e\cS_1-\nabla_eV \in L^2(\nu)$.
Therefore, $(\cE^{\nu}, C^\infty(G^{E^+_\Lambda}))$ is closable (c.f. \cite[Proposition 3.7]{SZZ22}, \cite[Chap. I, Prop. 3.3]{MaRockner}),  and its closure $(\cE^{\nu},D(\cE^{\nu}))$ is a regular Dirichlet form on $L^2(G^{E^+_\Lambda},\nu)$.

In Lemma~\ref{log:1} below,
the proof of \eqref{log:muQ} will basically follow the strategy \eqref{e:strategy},
whereas the proof of \eqref{log:nu} will follow the lines of  \eqref{e:PfLem43};
both involving verifying  the Bakry--\'Emery condition.
To verify  the B\'E condition for $\mu_Q$ and $\nu$ we  need to calculate
$\Hess_{\cS_2}(v^\Phi,v^\Phi)$, $\Hess_{\cS_1}(v^Q,v^Q)$ and $\Hess_{V}(v^Q,v^Q)$ for  $v^\Phi\in \mR^{N\Lambda}$ and $v^Q\in T_QG^{E^+_\Lambda}$.
Following the convention in Section~\ref{sec:product},
we write
\begin{equ}[e:v]
v^\Phi=(v_x)_{x\in \Lambda},\qquad 	v^Q=(v_e)_{e\in E^+_\Lambda}=\sum_{e\in E_\Lambda^+}X_eQ_e
\end{equ}
with $X_e \in \mfg^{E_\Lambda^+}$  being zero for 
all except the component indexed by $e$.

To state the lemma, we also
define a constant $K_{\R^N}^\nu$
  by
\begin{equ}[e:wPcon]
	K_{\R^N}^\nu
	\eqdef\frac14(N-2)
	-\frac{\kappa N}{m}
	-2N \frac{\kappa^2}{m^2}-8(d-1)N|\beta|,
\end{equ}
which will be the lower bound
 of the B\'E condition \eqref{con:B-E} for $\nu$ and will be required to be positive.
The strict positivity requirement of $K_{\R^N}^\nu$ means that $\kappa/m$ and $\beta$ are required to be small
(this turns out to be  precisely the condition in Theorem~\ref{main:2}).

\bl\label{log:1}
For $m, \kappa>0$,
the log-Sobolev and the Poincar\'e inequalities hold for $\mu_Q$,
namely,
 for $F\in C^\infty(\mR^{N\Lambda})$
\begin{align}\label{log:muQ}
	\text{ent}_{\mu_Q}(F^2)\leq \frac1{mN}\cE^{\mu_Q}(F,F)
	\qquad \text{ and } \qquad \var_{\mu_Q}(F)\leq \frac1{2mN}\cE^{\mu_Q}(F,F).
\end{align}
Note that the constants in \eqref{log:muQ} are independent of $Q$.
If moreover one has $K_{\R^N}^\nu>0$, 
then the log-Sobolev and the Poincar\'e inequalities hold for $\nu$, 
namely,
 for $F\in C^\infty(G^{E_\Lambda^+})$,
\begin{align}\label{log:nu}
	\text{ent}_{\nu}(F^2)\leq \frac2{K_{\R^N}^\nu}\cE^{\nu}(F,F)
\qquad \text{ and } \qquad \var_{\nu}(F)\leq \frac1{K_{\R^N}^\nu}\cE^{\nu}(F,F).
\end{align}
\el
\begin{proof}
	We have for $v^\Phi=(v_x)$
	\begin{equation}\label{eq:Hs-Rn}
		\aligned
		\Hess_{\cS_2}(v^\Phi,v^\Phi)
		&= \kappa N \sum_{e}\Big[2|v_x|^2+2|v_y|^2-4v_x^tQ_ev_y\Big]+2mN\sum_x|v_x|^2
		\geq 2mN|v|^2.
		\endaligned
	\end{equation}
	Hence, by the Bakry--\'Emery condition \eqref{con:B-E}, for $m,\kappa>0$, the log-Sobolev and the Poincar\'e inequalities \eqref{log:muQ} hold for $\mu_Q$.
	
	Now we turn to $\nu$.
	 Recall the definition of $V$ in \eqref{def:VQ}. By direct calculation we  have, for $v^Q = (X_e Q_e)_{e\in E^+}$,
	\begin{equ}
		\Hess_V(v^Q,v^Q)
		=-2\kappa N\sum_e\Big[\E_{\mu_Q}(\Phi_x^tX_e^2Q_e\Phi_y)+2\kappa N \,\var_{\mu_Q}(\Phi_x^tX_eQ_e\Phi_y)\Big]
	\end{equ}
	where we write $e=(x,y)$ as before.
	Hence, we have
	\begin{equ}\label{est:h1}
		|\Hess_V(v^Q,v^Q)|
		\leq
		2\kappa N\sum_e\Big[|X_e|^2\E_{\mu_Q}(|\Phi_x|\cdot|\Phi_y|)+2\kappa N\,\var_{\mu_Q}(\Phi_x^tX_eQ_e\Phi_y)\Big].
	\end{equ}
Our next step is then to upper-bound the right-hand side of \eqref{est:h1}.
We will use \eqref{b:phi} in Lemma~\ref{lem:b:phi} to bound
the first term
on the RHS of \eqref{est:h1},
and invoke the Poincar\'e inequality \eqref{log:muQ} for $\mu_Q$ already proved above  to bound
the second term in the RHS of \eqref{est:h1}.

By \eqref{b:phi} the first term in \eqref{est:h1} can be bounded by $2\kappa N|v_Q|^2\frac{1}{2m}$.
Moreover, by the Poincar\'e inequality \eqref{log:muQ} for $\mu_Q$, we have
	\begin{align*}
		\Big|\sum_e\var_{\mu_Q}(\Phi_x^tX_eQ_e\Phi_y)\Big|\leq \frac1{2mN} \sum_e\sum_z\E_{\mu_Q}(|\nabla_z F_e|^2)\leq \frac{1}{2m^2 N}|v_Q|^2,
	\end{align*}
	where $F_e=\Phi_x^tX_eQ_e\Phi_y$ and we used that $\nabla_zF_e\neq0$ iff $z=x, y$ and $\E_{\mu_Q}(|\nabla_xF_e|^2)\leq \frac{1}{2m}|v_e|^2$ by direct calculation.
	
	Combining the above calculations we obtain
	\begin{align*}
		|\Hess_V(v^Q,v^Q)|\leq\Big(\frac{\kappa N}{m}+2\kappa^2\frac{N}{m^2}\Big)|v^Q|^2.
	\end{align*}
	By the calculation in the pure Yang--Mills case  (\cite[Lemma~4.1]{SZZ22}), we have
	\begin{equ}[e:YMHess]
		|\Hess_{\cS_1}(v^Q,v^Q)|\leq 8(d-1)N|\beta||v^Q|^2.
	\end{equ}
	Recall also that
	\begin{equ}[e:RicciG]
	\Ric(v^Q,v^Q)
	=\frac14(N-2) |v^Q|^2
	\end{equ}
	 which follows from \cite[(F.6)]{AGZ} or \cite[(4.8)]{SZZ22}.
	Hence,   combining the last three equations,
	we see that for $K_{\R^N}^\nu>0$,
	we have verified the Bakry--\'Emery condition \eqref{con:B-E} for $\nu$. Therefore, the log-Sobolev and the Poincar\'e inequalities \eqref{log:nu} hold. 
\end{proof}

Even if we have the log-Sobolev inequality for $\mu_Q$ and $\nu$, the log-Sobolev or the Poincar\'e inequality does not automatically hold for $\mu_\Lambda$. To derive the Poincar\'e inequality for $\mu_\Lambda$, we need to control the correlation function under $\mu_Q$ (see \eqref{nab:e} below). Hence, we  first establish the mass gap property for $\mu_Q$ with fixed $Q$ to achieve better control of the correlation.
To this end, we  consider the generator $(\cL^Q,D(\cL^Q))$ corresponding to $\cE^{\mu_Q}$ in the sense that $\cE^{\mu_Q}(F,G)=-\int \cL^Q FG\dif \mu_Q$ for $F, G\in D(\cL^Q)$. For smooth  functions $F\in C^\infty(\mR^{N\Lambda})$,
\begin{align*}
	\cL^QF=\sum_{x\in \Lambda}\Delta_x F-\sum_{x\in \Lambda}\<\nabla_x\cS_2,\nabla_xF\>.
\end{align*}
We write $(P_t^Q)_{t\geq0}$ for  the related Markov semigroup (c.f. \cite{Fukushima}), which can also be constructed from dynamics \eqref{SDE:Q}.
In the rest of this section,
for $F\in C^\infty(\mR^{N\Lambda})$ we write $\Lambda_F\subset \Lambda$ for
\begin{equ}[def:LambdaF0]
\Lambda_F =\text{ the set of the  points $F$ depends on.}
\end{equ}
Let $|\Lambda_F|$ denote the
cardinality of $\Lambda_F$.

Note that the following lemma holds without smallness assumptions on $\beta$ and $\kappa$.

\bl\label{lem:co}
Let $m,\kappa>0$. For  $F, H \in C^\infty(\mR^{N\Lambda})$, assuming $\Lambda_F\cap \Lambda_H=\emptyset$ with $\Lambda_F,\Lambda_H$ defined by \eqref{def:LambdaF0},
one has
\begin{align}\label{eq:mass-gap}
	|\cov_{\mu_Q}(F,H)|\leq C_1e^{-C_Nd(\Lambda_F,\Lambda_H)}\$F\$_{2,Q}\$H\$_{2,Q},
\end{align}
with  $C_1=C|\Lambda_F||\Lambda_H|$ and $C, C_N$  independent of $Q$, $F$ and $H$.
Here $C_N$ depends on $K_\cS$,  $N$ and $d$,  and $\$F\$_{2,Q}\eqdef \sum_x\|\nabla_x F\|_{L^2(\mu_Q)}$.
\el

\begin{proof}
	Since we have the log-Sobolev and the Poincar\'e inequalities for $\mu_Q$ by Lemma \ref{log:1},  we can use similar arguments as in the proof of \cite[Corollary 4.11]{SZZ22} to conclude the results.
	More precisely, since
	$P_t^Q$ leaves $\mu_Q$ as an invariant measure, we have
	\begin{align}\label{coQ}
		|\cov_{\mu_Q}(F,H)|&=|\E_{\mu_Q}(P_t^{Q}(FH)-P_t^{Q}FP_t^Q H)+\cov_{\mu_Q}(P_t^Q F,P_t^Q H)|\no
		\\&\leq |\E_{\mu_Q}(P_t^{Q}(FH)-P_t^{Q}FP_t^Q H)|+\var_{\mu_Q}(P_t^Q F)^{1/2}\var_{\mu_Q}(P_t^Q H)^{1/2}.
	\end{align}
	Since we have the Poincar\'e inequality for $\mu_Q$, we use \cite[Theorem 5.6.1]{Wang} (see also (4.9) in \cite{SZZ22}) to have
	\begin{align}
		&	\var_{\mu_Q}(P_t^QF)\leq\frac1{2mN} \sum_x\E_{\mu_Q}(|\nabla_x P_t^QF|^2)\no
		\\&\leq\frac1{2mN} e^{-4tmN}\sum_x\E_{\mu_Q}((P_t^Q|\nabla_x F|)^2)\leq\frac1{2mN} e^{-4tmN}\sum_x\E_{\mu_Q}(|\nabla_x F|^2).\label{coQ1}
	\end{align}
Here, we refrain from employing the bound $\var_{\mu_Q}(P_t^QF)\leq e^{-4tmN} \sum_x\E_{\mu_Q}(|F|^2)$ through the Poincar\'e inequality, as utilized in \cite{SZZ22}. Our objective is to retain $\sum_x\E_{\mu_Q}(|\nabla_x F|^2)$ on the RHS. This is essential for establishing the Poincar\'e inequality for $\mu_\Lambda$ in Theorem \ref{th:poin}.

	We also have as in the proof of \cite[Corollary 4.11]{SZZ22}
	\begin{align*}
		P_t^{Q}(FH)-P_t^{Q}FP_t^Q H=2\sum_x\int_0^tP_s^Q\<\nabla_x P_{t-s}^Q F,\nabla_x P_{t-s}^Q H\>\dif s.
	\end{align*}
	Recall that
	\begin{align*}
		\nabla_x P_{t}^Q F-P_t^Q \nabla_x F=\int_0^tP_{t-s}^Q[\nabla_x,\cL^Q]P_s^Q F\dif s.
	\end{align*}
	By direct calculation we have
	$$|[\nabla_x,\cL^Q]F|\leq \sum_{y\sim x}a_{xy}|\nabla_y F|,$$
	for $a_{xy}= 2\kappa N^{5/2}$ for $x\neq y$ and $a_{xx}=2(2d\kappa+m)N^{5/2}$,  since for a basis $\{v_x^i\}$ in $\mR^N$ $|v_x^iv_y^j\cS_2|\leq 2\kappa N$ for $x\sim y$ and $y\neq x$ and $|v_x^iv_x^j\cS_2|\leq 2N(2d\kappa+m)$, where we used for each $x$ there are at most $2d$ edges containing $x$.
	Hence, we derive
	\begin{align*}
		\|\nabla_x P_{t}^Q F\|_{L^2(\mu_Q)}\leq \|\nabla_x F\|_{L^2(\mu_Q)}+\int_0^t\sum_{y\sim x}a_{xy}\|\nabla_y P_{s}^Q F\|_{L^2(\mu_Q)}\dif s.
	\end{align*}
		The same arguments as in the proof of  \cite[Corollary 4.11]{SZZ22} implies that for any $c>0$ there exists $B>0$ such that for $d(x,\Lambda_F)\geq Bt$
		\begin{align}\label{nab:x}
			\|\nabla_x P_{t}^Q F\|_{L^2(\mu_Q)}\leq e^{-2cd(x,\Lambda_F)}\$F\$_{2,Q}.
		\end{align}
		where we change the corresponding $L^\infty$-norm to $L^2$-norm compared to \cite{SZZ22}.
		We then choose $t\sim d(\Lambda_F,\Lambda_H)/B$ below and consider $\sum_x\<\nabla_xP_{t-s}^Q F,\nabla_xP_{t-s}^Q H\>$. For $x\in \Lambda_H$ we use \eqref{nab:x} for $\|\nabla_xP_{t-s}^Q F\|_{L^2(\mu_Q)}$ and use \eqref{coQ1} for $\|\nabla_x P_{t-s}^QH\|_{L^2(\mu_Q)}$ to have
		\begin{align*}
			\|\<\nabla_xP_{t-s}^Q F,\nabla_xP_{t-s}^Q H\>\|_{L^1(\mu_Q)} &\leq \|\nabla_xP_{t-s}^Q F\|_{L^2(\mu_Q)}\|\nabla_xP_{t-s}^Q H\|_{L^2(\mu_Q)}
			\\&\leq e^{-c\,d(\Lambda_F,\Lambda_H)}e^{-c\,d(x,\Lambda_F)}\$F\$_{2,Q} \$ H\$_{2,Q}.
		\end{align*}
		The same bound also holds for $x\in \Lambda_F$. For $x\notin \Lambda_F\cup \Lambda_H$, we have $d(x,\Lambda_F)\vee d(x,\Lambda_H)\geq d(\Lambda_F,\Lambda_H)/2$.
		Hence, we use \eqref{nab:x} to have
		\begin{align*}
			\|\<\nabla_xP_{t-s}^Q F,\nabla_xP_{t-s}^Q H\>\|_{L^1(\mu_Q)}\leq e^{-c\, d(\Lambda_F,\Lambda_H)/2-c\,(d(x,\Lambda_F)\wedge d(x,\Lambda_H))/2}\$F\$_{2,Q} \$H\$_{2,Q}.
		\end{align*}
		Combining the above calculations and taking sum over $x$ we obtain
		\begin{align}\label{coQ2}
			\sum_x\|\<\nabla_xP_{t-s}^Q F,\nabla_xP_{t-s}^Q H\>\|_{L^1(\mu_Q)}\leq C_1e^{-c \,d(\Lambda_F,\Lambda_H)/4}\$F\$_{2,Q} \$H\$_{2,Q}.
		\end{align}
		Substituting \eqref{coQ1} and \eqref{coQ2} into \eqref{coQ} the results follow.
	\end{proof}

	%
	
	\br\label{rem:V} As mentioned in Section \ref{sec:es-1} (see  Remark \ref{rem:Other-M} below), it is straightforward to extend Lemma \ref{lem:co} to the case that $V(\Phi_x)=m|\Phi_x|^2+\lambda|\Phi_x|^4$ with $m>0, \lambda\geq0$.
	More precisely, we have
	$$|\nabla_x P_t^QF|^2-|P_t^Q\nabla_xF|^2\leq2\int_0^tP_{t-s}^Q\<\nabla_xP_s^QF,[\nabla_x,\cL^Q]P_s^QF\>\dif s, $$
	(c.f. \cite[Lemma 1.1]{Zegarlinski1996}) and we
	need to use $\lambda\geq0$ to have $$2\<\nabla_xP_s^QF,[\nabla_x,\cL^Q]P_s^QF\>\leq \sum_{y\sim x}a_{xy}|\nabla_yP_s^QF|^2$$ for suitable constants $a_{xy}$, which implies $$\|\nabla_xP_t^QF\|_{L^2(\mu_Q)}^2\leq \|\nabla_xF\|_{L^2(\mu_Q)}^2+\int_0^t\sum_{y\sim x}a_{xy}\|\nabla_yP_s^Q F\|_{L^2(\mu_Q)}^2\dif s.$$
	By iteration we obtain \eqref{nab:x} and mass gap of $\mu_Q$, i.e. \eqref{eq:mass-gap}.
	\er
	
\br\label{rem:co} In the proof of Theorem~\ref{th:poin} below, we will also use  mass gap for $\mu_Q^{\Lambda_0}$, where $\Lambda_0 \subset \Lambda$ and   $\mu_Q^{\Lambda_0}$ is the regular conditional probability of $\mu_Q$ given $\Phi_y$ for $y\in \Lambda\backslash \Lambda_0$. More precisely,
\begin{align*}
	\mu_Q^{\Lambda_0}\propto\exp\Big(&-\kappa N\sum_{e\in E_{\Lambda_0}^+}|Q_e\Phi_y-\Phi_x|^2-\kappa N\sum_{x\in \Lambda_0, y\in \Lambda\backslash\Lambda_0}(|\Phi_x|^2-2\Phi_x^tQ_e\Phi_y)\\&\qquad\qquad\qquad-mN\sum_{z\in \Lambda_0}|\Phi_z|^2\Big)\prod_{x\in \Lambda_0}\dif \Phi_{x}.
\end{align*}
We can easily verify the Bakry--\'Emery condition as in Lemma \ref{log:1} with a uniform lower bound independent of $\Lambda_0$ and $\Lambda$. Hence, the same argument as in Lemma \ref{lem:co} implies that for smooth functions $F, H$ on $\mR^{N\Lambda}$
\begin{align}\label{co:cond}
	\cov_{\mu_Q^{\Lambda_0}}(F,H)\leq C_1e^{-C_Nd(\Lambda_F^0,\Lambda_H^0)}\Big(\sum_{x\in\Lambda_0}\|\nabla_xF\|_{L^2(\mu_Q^{\Lambda_0})}\Big)\Big(\sum_{x\in\Lambda_0}\|\nabla_xH\|_{L^2(\mu_Q^{\Lambda_0})}\Big),
\end{align}
where $C_1=C|\Lambda_F^0||\Lambda_H^0|$ with  $\Lambda_F^0=\Lambda_F\cap \Lambda_0$, and $C$, $C_N$ can be chosen the same as in Lemma \ref{lem:co}.
\er
	
	Now we are in a position to prove the Poincar\'e inequality for $\mu_\Lambda$.
	Recall $\cE^\Lambda$ and $\cE^\mu$  introduced in \eqref{def:EL} and \eqref{e:EFF}.

	\bt\label{th:poin}
	For $m, \kappa>0$ and $	K_{\R^N}^\nu>0$ with $	K_{\R^N}^\nu$ defined in \eqref{e:wPcon},
	the  Poincar\'e inequality holds for $\mu_{\Lambda}$,
	 i.e. there exists a  universal constant $C>0$ depending on $N, d$ but  independent of $\Lambda$ such that
	 \begin{equ}[e:mu-poin]
	 	\var_{\mu_{\Lambda}}(F)\leq C\cE^{\Lambda}(F,F).
	 \end{equ}
	Moreover,
	 for any infinite volume limit $\mu$  of $\mu_\Lambda$,
	\begin{align*}
		\var_\mu(F)\leq C\cE^\mu(F,F)
	\end{align*}
	for any  $F\in C^\infty_{cyl}(\cQ)$.
	\et
	\begin{proof}
It suffices to prove the claim for $\mu_\Lambda$ uniform in $\Lambda$ and the result for $\mu$ immediately follows.
By the disintegration \eqref{e:disinteg}, we have
		\begin{align*}
			\var_{\mu_\Lambda}(F)=\E_\nu(\var_{\mu_Q}(F))+\var_\nu(H),
		\end{align*}
		with 
		$H(Q)=\E_{\mu_Q}(F)$.
		By Lemma \ref{log:1},
		i.e. the Poincar\'e inequalities \eqref{log:muQ} \eqref{log:nu} for $\mu_Q$ and $\nu$,
		\begin{align}\label{poin}
			\var_{\mu_\Lambda}(F)
			\leq
			\frac1{2mN}\sum_x\E_{\mu_{\Lambda}}(|\nabla_xF|^2)
			+\frac1{K_{\mR^N}^\nu}\sum_e\E_\nu(|\nabla_e H|^2).
		\end{align}
		The term $\sum_x\E_{\mu_{\Lambda}}(|\nabla_xF|^2)$
		is part of the desired terms in \eqref{def:EL}.
		In the rest of the proof  we estimate the second term $\sum_e\E_\nu(|\nabla_e H|^2)$.
		
		To this end we calculate $\nabla_eH$ and by Lemma \ref{lem:nab} to have
		\begin{align}\label{nab:e}
			\nabla_e H=\E_{\mu_Q}(\nabla_e F)-\kappa N\cov_{ \mu_Q}(F,Q_e(\Phi_y\Phi_x^t)Q_e-\Phi_x\Phi_y^t).
		\end{align}
	The first term on the RHS of \eqref{nab:e} can be easily bounded.
	For most of the proof below, we estimate the following correlation function
	\begin{equ}[e:covFg]
	\cov_{ \mu_Q}(F,g_e) ,\qquad  \mbox{with}\quad g_e \eqdef Q_e(\Phi_y\Phi_x^t)Q_e-\Phi_x\Phi_y^t.
	\end{equ}

A naive bound for $\cov_{ \mu_Q}(F,g_e)$ is to apply the Cauchy--Schwarz inequality and the Poincar\'e inequality to have
		\begin{align}\label{nai:bd}
			|\cov_{ \mu_Q}(F,g_e)|
			\leq \var_{\mu_Q}(F)^{1/2}\var_{\mu_Q}(g_e)^{1/2}
			\lesssim \cE^{\mu_Q}(F,F)^{1/2}
		\end{align}
		where the second step follows by Lemma~\ref{lem:b:phi} and boundedness of $Q$,
and this substituted into \eqref{poin} implies
		$$	\var_{\mu_\Lambda}(F)\lesssim \sum_x\E_{\mu_{\Lambda}}(|\nabla_xF|^2)+\sum_e\E_{\mu_{\Lambda}}(|\nabla_eF|^2)+\sum_e\E_\nu(\cE^{\mu_Q}(F,F)).$$
		The last term does not give the correct Dirichlet form and might be infinity.
	Even if we use Lemma~\ref{lem:co} instead of \eqref{nai:bd}, we would get a final bound depending on $\Lambda_F$, which is not  desired.
(In fact, one of the main points of the argument below is to ensure that the final bound does not depend on $\Lambda_F$.)
	
	Instead, we will integrate ``vertex by vertex''.
To this end, we write  $\Lambda=\{x_1,\dots,x_p\}$ for some $p\in \mN$, and $\Lambda_i=\{x_1,\dots,x_i\}$\footnote{Throughout this proof, we abuse the notation $\Lambda_i$, and it means a subset of $\Lambda$ and is different from   $\Lambda_L=\mZ^d\cap L\mT^d$ defined in Section \ref{sec:not}. } and set $\mu_Q^{\Lambda_i}$ as the regular conditional probability of $\mu_Q$ given $\Phi_y$ for $y\in \Lambda\backslash \Lambda_i$. 
For $i\ge 0$ we write 
$$
\cov_Q^{\Lambda_i}=\cov_{ \mu_Q^{\Lambda_i}},\qquad
\var_Q^{\Lambda_i}=\var_{ \mu_Q^{\Lambda_i}},\qquad
\E_Q^{\Lambda_i}=\E_{ \mu_Q^{\Lambda_i}},
\qquad
 F_i\eqdef \E_Q^{\Lambda_i}(F).
$$
The superscript $\Lambda_i$ in the above notation indicates that integration is only over 
the field $\Phi$ on $\Lambda_i$. In particular we have 
$F_0 =F$, and $\cov_Q^{\Lambda_0}$ and $\var_Q^{\Lambda_0}$ are just zero
since the entire field is fixed.

As in \eqref{def:LambdaF0} we write $\Lambda_{F_i}\subset \Lambda$
 for the set of the  points $F_i$ depends on.
(Of course $F_i$ also depends on the edge variables $Q$ but $Q$ is fixed here.)
By definition, $\Lambda_{F_i}=\Lambda_F \backslash \Lambda_i$.
	
\begin{figure}[h]
\begin{tikzpicture}
\draw (0,0) -- (5,0) -- (5,4) -- (0,4) -- (0,0);
\draw  (2,0) -- (2,2) -- (2.5,2) -- (2.5,4);
\node at (1,2) {$\Lambda_i$};
   \fill (0.2,3.8) circle (1pt); \node at (0.5,3.8) {$x_1$};
    \fill (2.2,2.2) circle (1pt); \node at (2.2,2.4) {$x_i$};
     \fill (2.2,1.8) circle (1pt); \node at (2.7,1.8) {$x_{i+1}$};
     \fill (4.8,0.2) circle (1pt); \node at (4.5,0.2) {$x_p$};
     \node at (3.7,2.8) {Fix $\Phi$};
     \node at (1.2,2.8) {Integrate $\Phi$};
\end{tikzpicture}
\caption{Illustration of conditioning in Step 1.}
\end{figure}
		
		{\bf Step 1. Calculation of $\cov_{ \mu_Q}(F,g_e)$.}
		
In this step we  prove  that there exists $c_N>0$ such that
\begin{equs}\label{cov:muq}
	|\cov_{ \mu_Q}(F,g_e)|
		 &\lesssim
		\sum_{i=1}^{p}\Big(\E_{\mu_Q}(|\nabla_{x_i}F|^2)^{1/2}+\E_Q^{\Lambda}(|\cov_Q^{\Lambda_{i-1}}(F,f_{x_i})|^2)^{1/2}\Big)e^{-c_Nd(e,x_i)},
	 \\
	 f_{x_{i}}
	 &\eqdef \sum_{y\sim x_i}Q_{(x_i,y)}\Phi_y 1_{y\in \Lambda_{i-1}}
\end{equs}
where $d(e,x_i)$ denotes the distance from $e$ to $x_i$ in $\mZ^d$.
	The proportional constant in \eqref{cov:muq} is independent of $F$ and $p$.
		
By telescoping,
definition of $F_i$, and
 property of conditional expectations, we can write
\begin{equs}[dec:cov]
	\cov_{ \mu_Q}(F,g_e)
	&=\E_Q^\Lambda\Big((F-\E^\Lambda_Q(F))g_e\Big)
	=\sum_{i=0}^{p-1}\E_Q^\Lambda\Big(\big(F_i-F_{i+1}\big)g_e\Big)
	\\
	&=\sum_{i=0}^{p-1}\E_Q^\Lambda\E_Q^{\Lambda_{i+1}}\Big(\big(F_i-\E^{\Lambda_{i+1}}_Q(F_i)\big)g_e\Big)
	=\sum_{i=0}^{p-1}\E_Q^\Lambda\cov_Q^{\Lambda_{i+1}}(F_i,g_e).
\end{equs}
	We then estimate each term on the RHS of \eqref{dec:cov}. 
	%
		We can  apply Remark \ref{rem:co}  for $\mu_Q^{\Lambda_{i+1}}$. Consequently, we use \eqref{co:cond} to obtain
		\begin{align}\label{cov:i1}
			|\cov_{Q}^{\Lambda_{i+1}}(F_{i},g_e)|
			\lesssim e^{-c_Nd(e,x_{i+1})} \E_Q^{\Lambda_{i+1}}(|\nabla_{x_{i+1}}F_i|^2)^{1/2}\sum_{y\in e}\E_Q^{\Lambda_{i+1}}(|\nabla_{y}g_e|^2)^{1/2}
		\end{align}
	where the implicit constant is independent of $i$ and we used $\Lambda_{i+1}\cap \Lambda_{F_i}\subset \{x_{i+1}\}$.
Considering the RHS of \eqref{cov:i1}, we have
		\begin{align}\label{nab:Fi}
			\nabla_{x_{i+1}}F_i=\E_Q^{\Lambda_{i}}(\nabla_{x_{i+1}}F)+2\kappa N\cov_Q^{\Lambda_{i}}(F,f_{x_{i+1}}),
		\end{align}
		where $f_{x_{i+1}}=\sum_{ y\sim x_{i+1}}Q_{(x_{i+1},y)}\Phi_y1_{ y\in\Lambda_{i}}$.
		Substituting \eqref{nab:Fi} into \eqref{cov:i1}, taking expectation $\E_Q^\Lambda$, and applying  H\"older's inequality and Lemma \ref{lem:b:phi}  we obtain
		\begin{align}\label{cov:i}
			\E_Q^\Lambda|\cov_{Q}^{\Lambda_{i+1}}(F_{i},g_e)|
			\lesssim \Big(\E_Q^\Lambda(|\nabla_{x_{i+1}}F|^2)^{1/2}+\E_Q^{\Lambda}(|\cov_Q^{\Lambda_{i}}(F,f_{x_{i+1}})|^2)^{1/2}\Big)e^{-c_Nd(e,x_{i+1})},
		\end{align}
			where we used Lemma \ref{lem:b:phi} to have $\sum_{y\in e}\E_{\mu_Q}(|\nabla_{y}g_e|^2)\lesssim1$.
	Here the implicit constant is independent of $i$.
		Substituting \eqref{cov:i} 
		 into \eqref{dec:cov} leads to \eqref{cov:muq}.

		{\bf Step 2. Calculation of $\cov_Q^{\Lambda_{i}}(F,f_{x_{i+1}})$.}
		
		In this step we
		prove
\begin{align}\label{cov:fx}
|\cov^{\Lambda_{i}}_Q(F,f_{x_{i+1}})|
\lesssim\sum_{k=0}^\infty e^{-c_Nk}
	\E_Q^{\Lambda_i}(\cE_Q^{\Lambda_{i,k+1}}(F)^{1/2}),
\end{align}
	where for $i=1,\dots, p-1$ and $k \ge 0$,
	\begin{equ}
	\Lambda_{i,k} 
	 \eqdef \Lambda_i\cap \{x \, : \,d(x,x_{i+1})\leq k\},
	\qquad
		\cE_Q^{\Lambda_{i,k+1}}(F) 
		\eqdef
		\sum_{y\in\Lambda_{i,k+1} }
		\E_Q^{\Lambda_{i,k+1}}
	(|\nabla_y F|^2).	\label{de:Eik}
	\end{equ}
Here as above $\E_Q^{\Lambda_{i,k}}$
is the integration over $\Phi$ on $\Lambda_{i,k}$ only.
	
\begin{figure}[h]
\begin{tikzpicture}
\draw (0,0) -- (5,0) -- (5,4) -- (0,4) -- (0,0);
\draw  (2,0) -- (2,2) -- (2.5,2) -- (2.5,4);
\node at (1,3) {$\Lambda_i$};
   \fill (0.2,3.8) circle (1pt); \node at (0.5,3.8) {$x_1$};
     \fill (2.2,1.8) circle (1pt); \node at (2.7,1.8) {$x_{i+1}$};
     \fill (4.8,0.2) circle (1pt); \node at (4.5,0.2) {$x_p$};
  \draw[dashed] (2.2,1.8) circle (1.2cm);
  \node at (1.8,2.2) {$\Lambda_{i,k}$};
\end{tikzpicture}
\caption{Illustration of $\Lambda_{i,k}$. $F_{i,k}$ depends on $\Phi$ outside $\Lambda_{i,k}$, and $f_{x_{i+1}}$ depends on $\Phi$ near $x_{i+1}$, which leads to a decay $e^{-c_Nk}$ in the correlation analysis}
\end{figure}
		A naive upper bound for $\cov_Q^{\Lambda_{i}}(F,f_{x_{i+1}})$ is obtained by applying the Cauchy-Schwarz inequality and the Poincar\'e inequality for $\mu_Q^{\Lambda_i}$, resulting in a bound of $C\cE^{\Lambda_i}_Q(F)^{1/2}$. However, substituting this bound into \eqref{cov:muq} and \eqref{poin} leads to a bound involving  $\sum_{e,i}\E_{\mu_\Lambda}(|\nabla_{x_1}F|^2)e^{-c_Nd(e,x_i)}$, which  depends on $p$. To address this issue, we further use idea from \cite[Section 4]{Yosida}  to decompose $\Lambda_i$ and apply the exponential decay property as stated in Lemma \ref{lem:co}, leading to a more refined bound \eqref{cov:fx}.
		
We set 
$$
F_{i,k} \eqdef \E^{\Lambda_{i,k}}_QF  \quad \mbox{for } k\geq1,
\qquad
\mbox{and}
\qquad
F_{i,0} \eqdef F.
$$
		Since $\E_Q^{\Lambda_i}\E_Q^{\Lambda_{i,k+1}}=\E_Q^{\Lambda_i}$,
		 we have
		\begin{equation}\label{cov:de}
			\aligned
			\cov^{\Lambda_{i}}_Q(F,f_{x_{i+1}})&=\sum_{k=0}^\infty	\E^{\Lambda_{i}}_Q[(F_{i,k}-F_{i,k+1})f_{x_{i+1}}]\\
			&=\sum_{k=0}^\infty	\E_Q^{\Lambda_i}\E_Q^{\Lambda_{i,k+1}}[(F_{i,k}-F_{i,k+1})f_{x_{i+1}}]
			\\&=\sum_{k=0}^\infty	\E_Q^{\Lambda_i}\cov_Q^{\Lambda_{i,k+1}}(F_{i,k},f_{x_{i+1}}).
			\endaligned
		\end{equation}
By Remark \ref{rem:co}, \eqref{co:cond}  holds for $\mu_Q^{\Lambda_{i,k+1}}$. Hence, we use  $\Lambda_{i,k+1}\cap \Lambda_{F_{i,k}}=\Lambda_{i,k+1}\backslash \Lambda_{i,k}$ to obtain
\begin{align}\label{cov:ik}
		&|\cov^{\Lambda_{i,k+1}}_Q(F_{i,k},f_{x_{i+1}})|\no
			\\
		&\lesssim (|k|^{d}+1)e^{-c_Nk}\Big(\sum_{y\in \Lambda_{i,k+1}\backslash \Lambda_{i,k}}\|\nabla_yF_{i,k}\|_{L^2(\mu_Q^{\Lambda_{i,k+1}})}\Big)\Big(\sum_{y\in \Lambda_{i,k+1}}\|\nabla_yf_{x_{i+1}}\|_{L^2(\mu_Q^{\Lambda_{i,k+1}})}\Big)\no
			\\
		&\lesssim (|k|^{3d/2}+1)e^{-c_Nk}\cE^{\Lambda_{i,k+1}}_Q(F_{i,k})^{1/2}\cE^{\Lambda_{i,k+1}}_Q(f_{x_{i+1}})^{1/2},
		\end{align}
		where  the proportional constant is independent of $k$ and we used that
		$|\Lambda_{F_{i,k}}\cap \Lambda_{i,k+1}|\lesssim |k|^d$, $d(\Lambda_{F_{i,k}},\Lambda_{f_{x_{i+1}}})\geq k$,  and
		$$\sum_{y\in \Lambda_{i,k+1}\backslash \Lambda_{i,k}}\|\nabla_y F_{i,k}\|_{L^2(\mu_Q^{\Lambda_{i,k+1}})}\lesssim (|k|^{d/2}+1) \cE_Q^{\Lambda_{i,k+1}}(F_{i,k})^{1/2}.$$
			We then estimate $\cE_Q^{\Lambda_{i,k+1}}(F_{i,k})$ from the RHS of \eqref{cov:ik}. To this end, we  have for $y\in \Lambda_{i,k+1}\backslash\Lambda_{i,k}$
		\begin{align*}
			\nabla_y F_{i,k}=\E_Q^{\Lambda_{i,k}}(\nabla_y F)+2\kappa N\cov_Q^{\Lambda_{i,k}}(F,f_y),
		\end{align*}
		with $f_y=\sum_{ x\sim y}Q_{(y,x)}\Phi_{x}1_ {x\in\Lambda_{i,k}}$.
		Note that we also have the Poincar\'e inequality for $\mu_Q^{\Lambda_{i,k}}$ with uniform constant, which implies that
		\begin{align*}
			|\cov_Q^{\Lambda_{i,k}}(F,f_y)|\lesssim \cE_Q^{\Lambda_{i,k}}(F)^{1/2}\cE_Q^{\Lambda_{i,k}}(f_y)^{1/2},
		\end{align*}
	where $\cE_Q^{\Lambda_{i,k}}(F)$ is defined as in \eqref{de:Eik}.
			Hence, we  derive
		\begin{align}
			\cE_Q^{\Lambda_{i,k+1}}(F_{i,k})&=\sum_{y\in\Lambda_{i,k+1}\backslash \Lambda_{i,k} }\E_Q^{\Lambda_{i,k+1}}
			(|\nabla_y F_{i,k}|^2)\no\\
			&\lesssim \sum_{y\in  \Lambda_{i,k+1}\backslash\Lambda_{i,k}}\E_Q^{\Lambda_{i,k+1}}(|\nabla_y F|^2)+\sum_{y\in  \Lambda_{i,k+1}\backslash\Lambda_{i,k}}\E_Q^{\Lambda_{i,k+1}}\Big(\cE_Q^{\Lambda_{i,k}}(F)\cE_Q^{\Lambda_{i,k}}(f_y)\Big)\no
			\\
			&\lesssim \sum_{y\in  \Lambda_{j,k+1}\backslash\Lambda_{i,k}}\E_Q^{\Lambda_{i,k+1}}(|\nabla_y F|^2)+(|k|^d+1)\cE_Q^{\Lambda_{i,k+1}}(F)\no
			\\
			&\lesssim (|k|^d+1)\cE_Q^{\Lambda_{i,k+1}}(F),\label{bd:E}
		\end{align}
		where the proportional constant is independent of $k$ and we used $|\Lambda_{i,k+1}\backslash\Lambda_{i,k}|\lesssim |k|^d+1$ and $\cE_Q^{\Lambda_{i,k}}(f_y)\lesssim 1$, since $|\nabla_x f_y|\lesssim 1$ for $|x-y|=1$, $\nabla_x f_y=0$ for $|x-y|>1$.
		
		We then substitute \eqref{bd:E} and the bound $\cE_Q^{\Lambda_{i,k+1}}(f_{x_{i+1}})\lesssim1$ into \eqref{cov:ik} to  get
		\begin{align}\label{cov:k}
			|\cov_Q^{\Lambda_{i,k+1}}(F_{i,k},f_{x_{i+1}})|\lesssim (|k|^{2d}+1)e^{-c_Nk}\cE_Q^{\Lambda_{i,k+1}}(F)^{1/2}.
		\end{align}
		Substituting \eqref{cov:k} into \eqref{cov:de} and using $(|k|^{2d}+1)e^{-\frac{c_N}2k}\lesssim 1$ uniformly in $k$, we arrive at \eqref{cov:fx}.

		{\bf Step 3. Proof of the Poincar\'e inequality}
		
		In this step we combine \eqref{cov:fx} and \eqref{cov:muq} to conclude the results. More precisely,
		combining \eqref{cov:fx} with \eqref{cov:muq} we derive
		\begin{align*}
			|\cov_{ \mu_Q}(F,g_e)|
			&\lesssim\E_{\mu_Q}(|\nabla_{x_1}F|^2)^{1/2}e^{-c_Nd(e,x_1)}+ \sum_{i=2}^{p}\Big(\E_{\mu_Q}(|\nabla_{x_i}F|^2)^{1/2}
			\\&\quad+\sum_{k=0}^\infty e^{-c_Nk}\E_{\mu_Q}(\cE_Q^{\Lambda_{i-1,k+1}}(F))^{1/2}\Big)e^{-c_Nd(e,x_i)},
		\end{align*}
	where $\cE_Q^{\Lambda_{i-1,k+1}}(F)$ is defined as in \eqref{de:Eik}.
		Hence, we use H\"older's inequality w.r.t. $i$ and $\sum_{i}e^{-c_Nd(e,x_i)}<\infty$ to obtain
		\begin{align*}
			|\cov_{ \mu_Q}(F,g_e)|^2
			&\lesssim\sum_{i=1}^p\E_{\mu_Q}(|\nabla_{x_i}F|^2)e^{-c_Nd(e,x_i)}\\&\quad+ \sum_{i=2}^{p}\sum_{k=0}^\infty e^{-c_Nk}\E_{\mu_Q}(\cE_Q^{\Lambda_{i-1,k+1}}(F))e^{-c_Nd(e,x_i)},
		\end{align*}
		where the proportional constant is independent of $e, i$ and $p$.
		We then use \eqref{nab:e} to have
		\begin{align}
			\sum_e\E_\nu(|\nabla_e H|^2) &\lesssim \sum_e\E_{\mu_\Lambda}(|\nabla_e F|^2)+\sum_{i=1}^p\E_{\mu_\Lambda}(|\nabla_{x_i}F|^2)\no
			\\&\quad + \sum_{i=2}^{p}\sum_{k=0}^\infty e^{-c_Nk}\E_\nu(\E_{\mu_Q}(\cE_Q^{\Lambda_{i-1,k+1}}(F))),\label{sum:nu}
		\end{align}
		where we used $\sum_ee^{-c_Nd(e,x_i)}<\infty$. Recall that $$\E_\nu(\E_{\mu_Q}(\cE_Q^{\Lambda_{i-1,k+1}}(F)))=\sum_{x\in \Lambda_{i-1,k+1}}\E_{\mu_\Lambda}(|\nabla_x F|^2)\;.$$
		Hence, for the last term in \eqref{sum:nu} we change the order of summation and write it as
		\begin{align*}
			&\sum_{k=0}^{\infty}\sum_{x\in \Lambda}\E_{\mu_\Lambda}(|\nabla_xF|^2)e^{-c_Nk}\sum_{i:x\in \Lambda_{i-1,k+1}}1
			\\&\lesssim \sum_{x\in \Lambda}\E_{\mu_\Lambda}(|\nabla_xF|^2)\sum_{k=0}^{\infty}	(|k|^{d}+1)e^{-c_Nk}\lesssim \sum_{x\in \Lambda}\E_{\mu_\Lambda}(|\nabla_xF|^2),
		\end{align*}
		where we bound $\sum_{i:x\in \Lambda_{i-1,k+1}}1$ by $(|k|^d+1)$  because, for a fixed $x$ there are at most $C(|k|^d+1)$ points in $\Lambda$ such that the distance to $x$ are smaller than $k+1$, with $C$ independent of $k$.
		Note that the above bound is independent of $p$ and we get
		\begin{align*}
			\sum_e\E_\nu(|\nabla_e H|^2)\lesssim \sum_e\E_{\mu_\Lambda}(|\nabla_e F|^2)+\sum_{x\in \Lambda}\E_{\mu_\Lambda}(|\nabla_{x}F|^2),
		\end{align*}
		which combined with \eqref{poin} implies the result \eqref{e:mu-poin}.
	\end{proof}
	
	\bc\label{co:ergodi}
	Under the same condition as in Theorem \ref{th:poin}, it holds that there exists a universal constant $C>0$ such that for any infinite volume limit $\mu$
	\begin{align*}
		\|P_tf-\mu(f)\|_{L^2(\mu)}\leq e^{-Ct}\|f\|_{L^2(\mu)},
	\end{align*}
where $(P_t)_{t\geq0}$ is the Markov semigroup associated with the dynamics \eqref{SDE:RNi}.
	\ec
	\begin{proof}
	This is a standard consequence of the
	Poincar\'e inequality (Theorem~\ref{th:poin}),
	as we alluded above \eqref{e:strategy},
	c.f. \cite[Section 1.1]{Wang}.
	\end{proof}

\br
For the pure Yang--Mills, 
 \cite[Section 5]{SZZ22}
proved exponential mixing under a Wasserstein distance (instead of $L^2$ as in Corollary~\ref{co:ergodi}) by a coupling argument,
and the positive Ricci curvature of the Lie group played a key role in the argument.
Here since $\Phi$ takes values in an unbounded space $\R^N$, constructing such a coupling for the Markov process poses a considerable challenge: we cannot simply get uniform control on the nonlinearity of the dynamics by exploiting  the Ricci curvature of the underlying manifold.
Nevertheless, we can still prove exponential mixing in $L^2$ in Corollary~\ref{co:ergodi}
by functional inequality approach, i.e. using Theorem~\ref{th:poin}.
\er	
	
\br\label{rem:strat2}
We remark on some works which may relate with our strategies,
continuing with Remark~\ref{rem:strat1}.
The bounds \eqref{e:covFg} and \eqref{nai:bd} would be similar with the methods in for instance \cite{RolandIsing},
but as explained above in our setting \eqref{nai:bd}  would not yield the desired result.
We need more calculation; in particular we  further decompose the domain.
The idea of integrating point by point to prove the log-Sobolev, as far as we know, goes back to \cite{LuYau1993} (see also the review by Ledoux \cite{MR1837286}).
\er

\subsection{Log-Sobolev inequalities for $M\in \{ \mS^{N-1}, G\}$}
\label{sec:LSB}

The following two  lemmas give  the log-Sobolev inequalities in the two cases where the Higgs fields are bounded.
In these cases we can apply the Bakry--\'Emery condition directly and we also do not need $m,\kappa>0$.
Recall $\mathcal S_{\YMH} = \cS_1-\cS_2$ as in Section~\ref{sec:YMHSDEs}.


To state Lemma~\ref{lem:ba-S} below,
we define a constant $K_{\mS^{N-1}}$,
which will be
the lower bound of the B\'E condition \eqref{con:B-E} for $\mu_\Lambda$,
as
\begin{equ}[e:KS-S]
	K_{\mS^{N-1}}
	\eqdef	
	\max_{\delta>0}\min\Big\{ K_{\YM}  -2(2\delta+1)|\kappa|N , \;\;
	(N-2)-4|\kappa|N d(2+1/\delta)\Big\}		 \tag{B\'E-$\mS$}
\end{equ}
where $K_{\YM}$ is defined in \eqref{e:K_YM},
and the maximum can be attained due to monotonicity of the two terms.
In fact, it turns out that \eqref{e:K_YM}+\eqref{e:small-kappa} is equivalent with $K_{\mS^{N-1}}>0$.
Indeed, a necessary condition for $K_{\mS^{N-1}}>0$ is  obviously
$$
  K_{\YM}>0,
\qquad
B_\kappa\eqdef (N-2)-8|\kappa|N d >0.
$$
Assume that these two inequalities hold. 
The mininum  in \eqref{e:KS-S} being positive is equivalent with
$$
\frac{
  K_{\YM} -2 |\kappa|N }{ 4 |\kappa| N}>\delta > 4|\kappa|N d /B_\kappa
  =\frac{4|\kappa|N d}{(N-2)-8|\kappa|N d},
$$
so to obtain an equivalent condition for $K_{\mS^{N-1}}>0$
we only need that the LHS is indeed larger than the RHS,
which is precisely \eqref{e:small-kappa}.

We then have the following result, under  small $(\beta,\kappa)$ condition.

\bl\label{lem:ba-S}
For $M= \mS^{N-1}$,
suppose that   
\begin{itemize}
\item
$K_{\YM}>0$ and
$ 2N (4d+ (N-2)K_{\YM}^{-1})|\kappa|< N-2$,
 with $K_{\YM}$ defined in
  \eqref{e:K_YM},
\item
or equivalently 
$K_{\mS^{N-1}}>0$, with $K_{\mS^{N-1}}$ defined in \eqref{e:KS-S}.
\end{itemize}
Then, the
log-Sobolev and the Poincar\'e inequalities hold for $\mu_{\Lambda}$, i.e. for $F\in C^\infty(\cQ_\Lambda)$
\begin{align*}
	\text{ent}_{\mu_{\Lambda}}(F^2)
	\leq \frac2{K_{\mS^{N-1}}}\cE^{\Lambda}(F,F),
	\qquad
	\text{var}_{\mu_{\Lambda}}(F)\leq \frac1{K_{\mS^{N-1}}}\cE^{\Lambda}(F,F).
\end{align*}

\el

\begin{proof}
	Consider a tangent vector $v=(v^Q,v^\Phi)\in T_QG^{E^+_\Lambda}\times T_\Phi M^{\Lambda}$ with $M=\mS^{N-1}$.
	We have
	\begin{equation} \label{eq:ric}
		\aligned
		\Ric(v,v)
		&=\Ric(v^Q,v^Q)+\Ric(v^\Phi,v^\Phi)
		\\	&=\frac14(N-2)|v^Q|^2+(N-2)|v^\Phi|^2.
		\endaligned
	\end{equation}
	Here we used \eqref{e:RicciG}, and the well-known
	the formula $\Ric(u,u)=\frac{N-2}{r^2}|u|^2$ for tangent vector $u$ of a sphere of radius $r$ in $\R^N$,
	and standard fact for Ricci curvatures of product manifolds.
	
	Write $v^Q=(v_e)_{e\in E^+_\Lambda}$, $v^\Phi=(v_x)_{x\in \Lambda}$ and $|v^Q|^2=\sum_e |v_e|^2$,  $|v^\Phi|^2=\sum_x |v_x|^2$.
One has (c.f. \cite[Sec. 1.10]{FigalliVillani})
	\begin{equ}[e:Hess-vv]
		\Hess_{\cS(Q,\Phi)}(v,v)=\frac{\dif^2}{\dif t^2}\Big|_{t=0}\cS(\Gamma(t)),
	\end{equ}
	with $\Gamma(t)=(\exp_Q(tv^Q), \exp_\Phi(tv^\Phi))$. On the product manifold we have
	\begin{equ}[e:Ga-t]
		\Gamma(t)=((\exp_{Q_e}(tv_e))_{e\in E^+},(\exp_{\Phi_x}(tv_x))_{x\in \Lambda_L}),
	\end{equ}
	with $v_e=X_eQ_e$
	\begin{equs}
		\exp_{Q_e}(tv_e) & =\exp(tX_e)Q_e,
		\\
		\exp_{\Phi_x}(tv_x) & =	\cos(|v_x|t)\Phi_x+\sin(|v_x|t)\frac{ v_x}{|v_x|}.
	\end{equs}
	By the calculation in the pure Yang-Mills case (\cite[Lemma~4.1]{SZZ22}),
	\begin{equ}[e:HessYM]
		|\Hess_{\cS_1}(v^Q,v^Q)|\leq 8(d-1)N|\beta||v^Q|^2.
	\end{equ}
Recall that
$\cS_2= - 2\kappa N\sum_{e\in E_\Lambda^+}\Phi_x^tQ_e\Phi_y$.
Writing $v_e=X_e Q_e$ with $X_e\in \mfg$ we get
\begin{equs}
		\Hess_{\cS_2}(v,v)
		&=-2\kappa N\sum_{e\in E_\Lambda^+}\frac{\dif^2}{\dif t^2} \Big|_{t=0}\Big[\Big(\cos(|v_x|t)\Phi_x+\sin(|v_x|t)\frac{ v_x}{|v_x|}\Big)^t
		\\
		&\qquad\qquad\qquad\exp\Big(tX_e\Big)Q_e \cdot \Big(\cos(|v_y|t)\Phi_y+\sin(|v_y|t)\frac{v_y}{|v_y|}\Big)\Big]
		\\
		&= -2\kappa N\sum_{e\in E_\Lambda^+}\Big[2v_x^t Q_e v_y-(|v_x|^2+|v_y|^2)\Phi_x^tQ_e \Phi_y\\&\qquad\qquad\qquad+\Phi_x^tX_e^2Q_e\Phi_y+ 2v_x^tX_eQ_e\Phi_y+2\Phi_x^tX_eQ_ev_y\Big]
\end{equs}
where $e=(x,y)$. 
Bounding $|Q_e|$ and $|\Phi_x|$ by $1$, by the Cauchy--Schwarz inequality we have
\begin{equs}
|\Hess_{\cS_2}(v,v)|
		&\leq 2|\kappa| N\sum_{e\in E_\Lambda^+}
		\Big[|v_x|^2+ |v_y|^2+(|v_x|^2+|v_y|^2)
		+|X_e|^2
		+2|v_x||X_e|
		+2|X_e||v_y|\Big]
		\\
		&\leq 2|\kappa|N \Big((4d+2d/\delta)|v^\Phi|^2+(1+2\delta )|v^Q|^2 \Big),
		\label{eq:Hs}
\end{equs}
	for any $\delta>0$.
Here we used $\sum_e|v_x|^2\leq d|v^\Phi|^2$ and $\sum_e|v_y|^2\leq d|v^\Phi|^2$,  since there are at most $d$ edges in $E^+_\Lambda$ with $x=u(e)$ and $y=v(e)$,
and Young's inequality
		\begin{align}\label{eq:young}
			2ab\leq a^2 /\delta+b^2\delta
		\end{align} 
		for  $|v_x||X_e|$ and $|X_e||v_y|$. 
		
		Combining \eqref{eq:ric}\eqref{e:HessYM}\eqref{eq:Hs},
		and recalling $K_{\YM}$ defined in \eqref{e:K_YM},
		we see that for $K_{\mS^{N-1}}>0$ the
		 B\'E condition \eqref{con:B-E}  is verified which proves the desired result.
%
\end{proof}



For the purpose of Lemma~\ref{lem:log-G}, we
define a constant $K_G$, which will be the lower bound of the B\'E condition \eqref{con:B-E} for the YMH measure $\mu_\Lambda$ with $M=G$, as
\begin{equ}[e:KS-G]
	K_G\eqdef	
	\max_{\delta>0}\min\Big\{
	K_{\YM} 
	-2|\kappa|N(1+2\delta), \;\; 
	\frac14 (N-2)
	-4|\kappa|Nd(2+1/\delta)\Big\}  .		\tag{B\'E-G}
\end{equ}
Similarly as above,
\eqref{e:K_YM}+\eqref{e:small-kappa} is equivalent with $K_G>0$.

\bl\label{lem:log-G}
For $M= G$,  suppose that 
\begin{itemize}
\item
$K_{\YM}>0$ and
$ 2N (16d+ (N-2)K_{\YM}^{-1})|\kappa|<  N-2$,
 with $K_{\YM}$ defined in
  \eqref{e:K_YM},
\item
or equivalently 
$K_{G}>0$ with $K_{G}$ defined in \eqref{e:KS-G}.
\end{itemize}
Then,
the
log-Sobolev and the Poincar\'e inequalities hold for $\mu_{\Lambda}$, i.e. for $F\in C^\infty(\cQ_\Lambda)$
\begin{align*}
	\text{ent}_{\mu_{\Lambda}}(F^2)\leq
	\frac2{K_G}\cE^{\Lambda}(F,F), \qquad
	\text{var}_{\mu_{\Lambda}}(F)\leq \frac1{K_G}\cE^{\Lambda}(F,F).
\end{align*}
\el
\begin{proof}
Recall that
$
\cS_2=-2\kappa N\sum_{e\in E_\Lambda^+}\tr(Q_e\Phi_y\Phi_x^t)$.
	For $v=(v^Q,v^\Phi)   \in T_{(Q,\Phi)}\cQ_L$.
	Instead of \eqref{eq:ric} we now have
	\begin{equ}
		\Ric(v,v)
		=\frac14(N-2) \left(|v^Q|^2
		+|v^\Phi|^2 \right).
	\end{equ}
	
	We write $v^Q=(v_e)_{e\in E^+_\Lambda}$, $v^\Phi=(v_x)_{x\in \Lambda}$ and proceed as in
	\eqref{e:Hess-vv} - \eqref{e:Ga-t} in the previous lemma.
	%
	Now with $v_e=X_eQ_e$ and $v_x=Y_x\Phi_x$
	\begin{align*}
		\exp_{Q_e}(tv_e)=\exp(tX_e)Q_e,\qquad \exp_{\Phi_x}(tv_x)=\exp(tY_x)\Phi_x.
	\end{align*}
	For the pure Yang--Mills part, we again have \eqref{e:HessYM}.
	For $\cS_2$ we have 
	\begin{equation}\label{eq:Hs-1}
		\aligned
		\Hess_{\cS_2}(v,v)
		&=-2N\kappa\sum_{e\in E_\Lambda^+}\frac{\dif^2}{\dif t^2} \Big|_{t=0}\tr\Big(\exp(tX_e)Q_e \cdot \exp(tY_y)\Phi_y \Big(\exp(tY_x)\Phi_x\Big)^t  \Big)
		\\&=-2N\kappa\sum_{e\in E_\Lambda^+}\Big(\tr(X_e^2Q_e\Phi_y\Phi_x^t)+\tr(Q_eY_y^2\Phi_y\Phi_x^t)+\tr(Q_e\Phi_y(Y_x^2\Phi_x)^t)
		\\&\qquad+2\tr(X_eQ_eY_y\Phi_y\Phi_x^t)+2\tr(X_eQ_e\Phi_y(Y_x\Phi_x)^t)+2\tr(Q_eY_y\Phi_y(Y_x\Phi_x)^t)\Big)
		\\&\leq2N|\kappa|\sum_{e\in E_\Lambda^+}(|X_e|^2+|Y_y|^2+|Y_x|^2+2|X_e||Y_y|+2|X_e||Y_x|+2|Y_x||Y_y|)
		\\&\leq 2N|\kappa|((1+2\delta)|v^Q|^2+(4d+2d/\delta)|v^\Phi|^2),
		\endaligned
	\end{equation}
where $\delta>0$ comes from the application of  Young's inequality \eqref{eq:young} and we use a  similar argument as in the last step of \eqref{eq:Hs}.	The lemma again follows from the Bakry--\'Emery criteria \eqref{con:B-E}  as in the proof of Lemma \ref{lem:ba-S}.
\end{proof}

Using \cite[Theorem 1.2]{MR1361304} (which states that  if spins take values in a compact manifold, the log-Sobolev inequality implies exponential ergodicity w.r.t. $L^\infty$-norm and uniqueness of the invariant measure).  Hence by Lemma \ref{lem:ba-S} or Lemma \ref{lem:log-G} we can deduce the following uniqueness and exponential ergodicity result.

\bc
\label{cor:SG-erg}
Assume the condition in Lemma \ref{lem:ba-S}  for $M= \mS^{N-1}$ or the condition in Lemma \ref{lem:log-G}  for $M= G$.   The invariant measure $\mu$ for the infinite volume dynamics \eqref{SDES} is unique  and satisfies the log-Sobolev and the Poincar\'e inequalities, i.e. for $F\in C^\infty_{cyl}(\cQ)$
\begin{align*}
	\text{ent}_{\mu}(F^2)\leq \frac2{K_M}\cE^\mu(F,F),\quad \quad \text{var}_{\mu}(F)\leq \frac1{K_M}\cE^{\mu}(F,F),
\end{align*}
where $K_M\in \{K_{ \mS^{N-1}}, K_G\}$ is as in   \eqref{e:KS-S} or   \eqref{e:KS-G}.
In particular, any infinite volume limit of $\mu_\Lambda$ is $\mu$.
 The associated Markov semigroup $(P_t)_{t\geq0}$ satisfies
\begin{align*}
	\|P_tF-\mu(F)\|_{L^\infty(\mu)}\leq C_1e^{-Ct}\Big(\sum_x\|\nabla_x F\|_{L^\infty}+\sum_e\|\nabla_e F\|_{L^\infty}\Big),
\end{align*}
for some universal constants $C_1, C>0$.
\ec

\br Since $\Phi_x$ is bounded, we could also construct a coupling as in \cite[Section 5]{SZZ22} to prove uniqueness of $\mu$ and its exponential ergodicity w.r.t. a suitable Wasserstein distance.
\er


\section{Mass gap}\label{sec:mass}
	
In this section we prove mass gap for the lattice Yang--Mills--Higgs measure under the assumptions on the model parameters as in Theorem~\ref{main4}. We 
will apply the functional inequalities obtained in Section \ref{sec:ergo}.

Many calculations below, such as the calculations for the commutators as alluded in Section~\ref{sec:intro}, will rely on a convenient choice
of the basis vectors in the tangent space of the configuration space.
Recall \eqref{e:TQ-basis} for the generic form of such a basis,
which consists of a basis of $T_QG^{E^+_\Lambda}$ and a basis of $T_\Phi M^{\Lambda}$.

To this end, we will take	
 an orthonormal  basis $\{v_e^i, i=1,\dots,d(\mfg), e\in E^+_\Lambda\}$ of $T_QG^{E^+_\Lambda}$, in the same way
  as in \cite[Section 4.3]{SZZ22}:
they are right-invariant vector fields on the Lie group $G^{E^+_\Lambda}$, obtained by a standard explicit basis of the Lie algebra $\mfg=\so(N)$, namely the skew-symmetric matrices $(\frac{1}{\sqrt 2} (e_{mn}-e_{nm}))_{1\le m,n\le N}$ where $e_{mn}$ is the $N\times N$ matrix
with its $(m,n)$-th entry being 1 and
all the other entries being $0$.

When $M\in \{G,\mS^{N-1}\}$
we will also choose a convenient basis $\{v_x^i, i=1,\dots, d_M, x\in \Lambda\}$ in $T_\Phi M^{\Lambda}$.
In particular, when $M=G$,  for each fixed $x\in \Lambda$,
we will take  $\{v_x^i, i=1,\dots, d_M\}$ (where $d_M=d(\mfg)$)
to be exactly the same  right-invariant vector fields as above.
The case $M=\mS^{N-1}$ is more subtle and will be described in Section~\ref{sec:massGS}.
	
\subsection{Mass gap for  $M=\mR^N$} 
\label{sec:mass-RN}

In this section we prove  the existence of a mass gap for the case where the Higgs field $\Phi$ is  $\mR^N$-valued. Given that $\Phi$ is unbounded, the basic method of calculating the commutator to establish the mass gap, 
as discussed above \eqref{e:idea-comm},
may lead to unbounded terms. To overcome this challenge, 
we proceed as in \eqref{eq:cor-i} \eqref{eq:cor-i1},
namely we turn to the Poincar\'e inequalities 
for $\mu_Q$ and $\nu$
 introduced in \eqref{def:muQ} and \eqref{def:nu}
which arise from the disintegration \eqref{e:disinteg}. 
As explained above \eqref{e:PfThm51},
we will need to estimate  
the commutators involving  the following generator  associated with 
the Dirichlet form $(\cE^\nu,D(\cE^\nu))$ defined in \eqref{dir:nu} on $C^\infty(G^{E^+_{\Lambda}})$
\begin{equ}[e:defLnu]
\cL^\nu F=\sum_e\Delta_eF +\sum_e\<\nabla_e\cS_1-\nabla_eV,\nabla_e F\>,
\end{equ}
with $\nabla_eV$ given in \eqref{def:neV},
where $V$  is the potential defined in \eqref{def:VQ} appearing in the measure $\nu$,
and $\cS_1$ is the Yang--Mills action as in above  \eqref{e:disinteg}.
Since $V$ is a nonlocal function on $Q$,
 we will leverage the mass gap of $\mu_Q$ (see Lemma \ref{lem:co}) to control $\nabla_e V$. 
	Denote by $(P_t^\nu)_{t\geq0}$  the related Markov semigroup (c.f. \cite{Fukushima}).
	In this section, for $F\in C^\infty_{cyl}(\cQ)$ we write
	 \begin{align}\label{def:LambdaF}
		\Lambda_F =\text{ the set of the edges and points $F$ depends on. }
	\end{align}
	Let $|\Lambda_F|$ denote the
	cardinality of $\Lambda_F$.
Recall the constant $K_{\mR^N}^\nu$ defined in \eqref{e:wPcon},
and that  $K_{\mR^N}^\nu >0$ is equivalent with 
the small $(\beta,\kappa/m)$ condition \eqref{e:small-bkm}.

We then have the following mass gap result 
under small $(\beta,\kappa/m)$ condition, which is the precise statement of 
Theorem~\ref{main4} in the case $M = \mR^N$.

	\bt\label{co:mass}
	Let $m,\kappa>0$. Assume that $K_{\mR^N}^\nu>0$, or equivalently, assume \eqref{e:small-bkm}. 
	
	For  $F, H\in C^\infty_{cyl}(\cQ)$, suppose that $\Lambda_F\cap \Lambda_H=\emptyset$ with $\Lambda_F, \Lambda_H$ defined in \eqref{def:LambdaF}. Let $\mu$ be an infinite volume limit of $\mu_\Lambda$.
	Then one has
	\begin{align*}
		|\cov_\mu(F,H)|\leq c_{1}  e^{-c_N d(\Lambda_F,\Lambda_H)}(\$F\$_2 \$H\$_2 +\|F\|_{L^2}\|H\|_{L^2}),
	\end{align*}
	where $c_{1}$  depends on $|\Lambda_F|$, $|\Lambda_H|$, $N$,
	and $c_N$ depends on $K_{\mR^N}^\nu$,  $N$ and $d$.
	Here 
	$$
	\$F\$_2\eqdef\sum_e\|\nabla_e F\|_{L^2(\mu)}+\sum_x\|\nabla_x F\|_{L^2(\mu)}.
	$$
	\et

	\begin{proof}
	We prove the result for $\mu_\Lambda$ with the constants $c_1, c_N$ independent of $\Lambda$. The result then follows by letting $\Lambda\to\mZ^d$.  We write $\E_\Lambda=\E_{\mu_\Lambda}$ and $\cov_\Lambda=\cov_{\mu_\Lambda}$.
Recall that the proof strategy was hinted in \eqref{e:PfThm51}.
		
		\medskip
	{\bf Step 1. Apply disintegration \eqref{e:disint} and mass gap for $\mu_Q$.}
		
	For cylinder functions $F$ and $H$ with $\Lambda_F\cap \Lambda_H=\emptyset$
		\begin{align}
			\cov_\Lambda(F,H)&=\E_\Lambda\Big((F-\E_\Lambda F)H\Big)\no
			\\&=\E_\nu\E_{\mu_Q}\Big((F-\E_{\mu_Q}F+\E_{\mu_Q}F-\E_\Lambda F)H\Big)\no
			\\&=\E_\nu\cov_{\mu_Q}(F,H)+\cov_{\nu}(\E_{\mu_Q}F,\E_{\mu_Q}H).\label{cov1}
		\end{align}
		Using Lemma \ref{lem:co} we have
		\begin{align*}
			|\cov_{\mu_Q}(F,H)|\leq C_1e^{-C_Nd(\Lambda_F,\Lambda_H)}\$F\$_{2,Q}\$H\$_{2,Q},
		\end{align*}
		with constants $C_1, C_N$, where $\$F\$_{2,Q}$ is defined in Lemma \ref{lem:co}.   Hence, by H\"older's inequality for $\nu$ we have
		\begin{align}\label{cov2}
			\E_\nu|\cov_{\mu_Q}(F,H)|\leq C_1e^{-C_Nd(\Lambda_F,\Lambda_H)}\$F\$_{2,\Lambda}\$H\$_{2,\Lambda},
		\end{align}
		with $\$F\$_{2,\Lambda}\eqdef \sum_e\|\nabla_e F\|_{L^2(\mu_\Lambda)}+\sum_x\|\nabla_x F\|_{L^2(\mu_\Lambda)}$.

				\medskip
	{\bf Step 2. Prove exponential decay of $\cov_{\nu}(\E_{\mu_Q}F,\E_{\mu_Q}H)$.}
	
		In the following we consider the second term  on the RHS of \eqref{cov1}, i.e.
		$
		\cov_{\nu}(F_1,H_1)
		$
		with $F_1=\E_{\mu_Q}F$ and $H_1=\E_{\mu_Q}H$. Recall that
		$P_t^\nu$ leaves $\nu$ as an invariant measure.
		We then have
		\begin{align}
			|\cov_{\nu}(F_1,H_1)|&=|\E_\nu(P_t^{\nu}(F_1H_1)-P_t^{\nu}F_1P_t^\nu H_1)+\cov_{\nu}(P_t^\nu F_1,P_t^\nu H_1)|\no
			\\&\leq|\E_\nu(P_t^{\nu}(F_1H_1)-P_t^{\nu}F_1P_t^\nu H_1)|+\var_{\nu}(P_t^\nu F_1)^{1/2}\var_\nu(P_t^\nu H_1)^{1/2}.\label{cov4}
		\end{align}
Since  the log-Sobolev inequality holds for $\nu$ (Lemma \ref{log:1}), we have
	\begin{align}\label{cov3}
			|\var_{\nu}(P_t^\nu F_1)|
			\leq e^{-2 K_{\mR^N}^\nu t}\|F_1\|_{L^2(\nu)}^2\leq e^{-2 K_{\mR^N}^\nu t}\|F\|_{L^2(\mu_\Lambda)}^2.
	\end{align}
		We also have as in \cite[Corollary 4.11]{SZZ22} (straightforward calculation using \eqref{e:defLnu})
		\begin{align*}
		P_t^{\nu}(F_1H_1)-P_t^{\nu}F_1P_t^\nu H_1
		=2\sum_e\int_0^tP_s^\nu
		\<\nabla_e P_{t-s}^\nu F_1,\nabla_e P_{t-s}^\nu H_1\>\dif s.
		\end{align*}
		Suppose for the moment that we can prove the following: for every edge $e$
		\begin{align}\label{nabF}
			\|\nabla_e P_{t}^\nu F_1\|_{L^2(\nu)}
			\leq C_1e^{C_0t-c_Nd(e,\Lambda_F)}\$F\$_{2,\Lambda},
		\end{align}
		for some universal constants $C_0, c_N>0$ independent of $F$, and $C_1>0$, which  depends on $\Lambda_F$ and $N$.
	We then choose
	$t = c_Nd(\Lambda_F,\Lambda_H) / (8C_0)$.
	 
Using  \eqref{nabF}, we then have
\begin{equs}[e:three-exp]
{} \|\< & \nabla_eP_{t-s}^\nu F_1,  \nabla_eP_{t-s}^\nu H_1\>\|_{L^1(\nu)}
\leq 
\|\nabla_eP_{t-s}^\nu F_1\|_{L^2(\nu)}\|\nabla_eP_{t-s}^\nu H_1\|_{L^2(\nu)}
\\&\leq 
C_1^2
\$F\$_{2,\Lambda} \$H\$_{2,\Lambda}\,e^{-c_N d(\Lambda_F,\Lambda_H)/4}
 \times
\begin{cases}
e^{-c_Nd(e,\Lambda_F)/4} 
& \text{if } e\in \Lambda_H 
\\
e^{-c_Nd(e,\Lambda_H)/4}  
& \text{if } e\in \Lambda_F
\\
e^{-c_N(d(e,\Lambda_F)\wedge d(e,\Lambda_H))/2}
& \text{if } e\notin \Lambda_F\cup \Lambda_H.
\end{cases}
\end{equs}
Here $e\notin \Lambda_F\cup \Lambda_H$, we used $d(e,\Lambda_F)\vee d(e,\Lambda_H)\geq d(\Lambda_F,\Lambda_H)/2$.
		Summing over $e$, since $\sum_{e\notin \Lambda_F}e^{-c_Nd(e,\Lambda_F)}$ is finite  we obtain 
		\begin{align}\label{cov5}
		|\E_\nu(P_t^{\nu}(F_1H_1)-P_t^{\nu}F_1P_t^\nu H_1)|
		&\leq 2	t\sum_e\|\<\nabla_eP_{t-s}^\nu F_1,\nabla_eP_{t-s}^\nu H_1\>\|_{L^1(\nu)}\no
		\\
		&\leq c_1e^{-c_N d(\Lambda_F,\Lambda_H)/4}\$F\$_{2,\Lambda} \$H\$_{2,\Lambda}.
		\end{align}
		Substituting \eqref{cov5} and \eqref{cov3} into \eqref{cov4} we obtain for $t= c_Nd(\Lambda_F,\Lambda_H)/(8C_0)$
		\begin{align*}
			|\cov_\nu(F_1,H_1)|\leq c_1e^{-c_N d(\Lambda_F,\Lambda_H)/8}\Big(\$F\$_{2,\Lambda} \$H\$_{2,\Lambda}+\|F\|_{L^2(\Lambda)}\|H\|_{L^2(\Lambda)}\Big)
		\end{align*}
		which combined with \eqref{cov2} and \eqref{cov1} the result follows. 
		Here  $c_N$ may differ from line to line.
		
		\medskip
	{\bf Step 3. Commutator estimates (which invoke the mass gap for $\mu_Q$).}
		
It remains to check the claimed bound \eqref{nabF}.
Recall the basis $\{v_e^i\}$ of the tangent space $G^{E^+_\Lambda}$ introduced in the beginning of this section. Consider
		\begin{align*}
			v_e^i P_{t}^\nu F_1-P_t^\nu v_e^i F_1=\int_0^t\frac{\dif}{\dif s}(P_{t-s}^\nu v_e^iP_s^\nu F)\dif s=\int_0^tP_{t-s}^\nu[v_e^i,\cL^\nu]P_s^\nu F_1\dif s.
		\end{align*}
		We first show that for smooth function $\widetilde F$ only depending on $Q$, and for every edge $e$,
	\begin{equs}\label{com:H}
	\left\|[v_e^i,\cL^\nu] \widetilde F \right\|_{L^2(\nu)}
	&\leq 
	\sum_{\bar e\in E_\Lambda^+}c_{e\bar e}\,
	e^{-2c_Nd(e,\bar e)}\|\nabla_{\bar e}\widetilde F\|_{L^2(\nu)}
	\\
	\max\{c_{ee},& \; c_{e\bar e} \; : \;  e,\bar e\in E_\Lambda^+\}\leq C_N
	\end{equs}
	for some coefficients $c_{e\bar e}>0$, and universal constants $c_N>0, C_N>0$.
It is important to note that 
unlike the simple situation discussed in \eqref{e:idea-comm},
here $\bar e$ is summed over {\it all} the edges,
and we will see that  this is due to non-locality of $V$. 

In fact, to prove \eqref{com:H}, 
we should start by calculating $[v_e^i,\cL^\nu] \widetilde F$ using the form of $\cL^\nu$ in \eqref{e:defLnu}. 
Such a calculation is done for the pure Yang--Mills case in  \cite[Lemma 4.10]{SZZ22},
with exactly the same choice of the basis $v_e^i$ and the pure Yang--Mills generator.
The main ingredients in this calculation are:

(1) The Beltrami--Laplacian in  \eqref{e:defLnu} does not contribute since 
it commutes with $v_e^i$ (this is because $\{v_e^i\}$ is chosen as right-invariant vector fields and the metric on $G$ is bi-invariant);

(2) $[\nabla_{v_e^i},\nabla_{\bar e}]=0$ when $e\neq \bar e$;

(3) $[v_e^i,v_e^j]$ can be calculated exactly since they are chosen as explicit matrices,
and can be easily bounded  which is  the content of  \cite[Lemma 4.9]{SZZ22}.

One then obtains the same result expressing $[v_e^i,\cL^\nu] \widetilde F$
as three terms called $I_1,I_2,I_3$ as in \cite[Lemma 4.10]{SZZ22},
with the only difference being that
 $\nabla \cS$  therein is now replaced by $\nabla_e(\cS_1-V)$,
the edge $\bar e$  is now summed over {\it all} edges.
Proceeding by the Cauchy--Schwarz inequality as in \cite[Lemma 4.10]{SZZ22} we  have
\begin{equ}\label{com-1}
	|[v_e^i,\cL^\nu] \widetilde F|
	\leq 
	\sum_{\bar e\in E_\Lambda^+}\Big(\sum_j|v_e^iv_{\bar e}^j(\cS_1-V)|^2\Big)^{1/2}|\nabla_{\bar e}\widetilde F|
	+\sqrt2|\nabla_e(\cS_1-V)||\nabla_e\widetilde F| \;.
\end{equ}
The terms  $v_e^iv_{\bar e}^j\cS_1$  and $\nabla_e \cS_1$ have been calculated in \cite{SZZ22};
in particular, since $\cS_1$ only involves local interaction,
 the derivative $\nabla_e \cS_1$ only depends on the variables $Q_{\bar e}$ with $\bar e\sim e$,
and  $v_{\bar e}v_e^i\cS_1 \neq 0$ only if $e\sim \bar e$.
This is {\it not} the case for  $v_e^iv_{\bar e}^jV$  and $\nabla_e V$,
since $V$ is nonlocal.
		Recall from \eqref{def:neV} that $\nabla_eV=\kappa N\E_{\mu_Q}(Q_e\Phi_y\Phi_x^tQ_e-\Phi_x\Phi_y^t)$,
		 and we use Lemma \ref{lem:b:phi} to have
		\begin{equ}[e:nabla_eV]
			|\nabla_eV|\lesssim \kappa N.
		\end{equ}
Since  $v_{\bar e}v_e^iV\neq 0$
for generic $e, \bar e$, 
we will use the mass gap of $\mu_Q$ from Lemma \ref{lem:co} to obtain better control on $v_e^iv_{\bar e}^jV$ below.

By direct calculation we have
\begin{equ}
	v_e^iv_e^jV=
	-2\kappa N\E_{\mu_Q}(\Phi_x^tX_e^jX_e^iQ_e\Phi_y)+4\kappa^2N^2\cov_{\mu_Q}(\Phi_x^tX_e^iQ_e\Phi_y,\Phi_x^tX_e^jQ_e\Phi_y)
\end{equ}
and for any $e\neq\bar e$, we have 
\begin{equ}
	v_e^iv_{\bar e}^jV=4\kappa^2N^2\cov_{\mu_Q}(\Phi_x^tX_e^iQ_e\Phi_y,\Phi_{\bar x}^tX_{\bar e}^jQ_{\bar e}\Phi_{\bar y}),
	\quad \mbox{where} \quad
	e=(x,y),\quad \bar e=(\bar x,\bar y)
\end{equ}
which may be non-zero for {\it all} (not just neighboring) edges  $e, \bar e$.
We then apply the Cauchy--Schwarz inequality and the Poincar\'e inequality \eqref{log:muQ} for $\mu_Q$ to have
		\begin{equs}[e:veveV]
			| v_e^iv_e^jV|
			&\lesssim \kappa N+\kappa^2 N^2\, \var_{\mu_Q}(\Phi_x^tX_e^iQ_e\Phi_y)^{1/2}\,\var_{\mu_Q}(\Phi_{ x}^tX_{ e}^jQ_{ e}\Phi_{ y})^{1/2}
			\\
			&\lesssim
			\kappa N
			+\kappa^2 N\, 
			\cE^{\mu_Q}(\Phi_x^tX_e^iQ_e\Phi_y)^{1/2}\,
			\cE^{\mu_Q}(\Phi_{ x}^tX_{ e}^jQ_{ e}\Phi_{ y})^{1/2}
			\\
			&\lesssim (\kappa +\kappa^2)N,
		\end{equs}
		where the last step follows from Lemma~\ref{lem:b:phi}. Here and in the sequel we write $\cE^{\mu_Q}(F)=\cE^{\mu_Q}(F,F)$. 
		For $e\neq \bar e$, we use mass gap of $\mu_Q$ (Lemma \ref{lem:co}) and Lemma \ref{lem:b:phi} to get
		\begin{equ}[e:vevebarV]
			|	v_e^iv_{\bar e}^jV|\lesssim \kappa^2N^2 e^{-2c_Nd(e,\bar e)}\cE^{\mu_Q}(\Phi_x^tX_e^iQ_e\Phi_y)^{1/2}\cE^{\mu_Q}(\Phi_{\bar x}^tX_{\bar e}^iQ_{\bar e}\Phi_{\bar y})^{1/2}\lesssim  e^{-2c_Nd(e,\bar e)}\kappa^2 N^2.
		\end{equ}
Hence, substituting the bounds
\eqref{e:nabla_eV}\eqref{e:veveV}\eqref{e:vevebarV}
 into \eqref{com-1},  we obtain \eqref{com:H}.
		
				\medskip
	{\bf Step 4. Iteration and bound $\|\nabla_e P_{t}^\nu F_1\|_{L^2(\nu)}$ (again using mass gap for $\mu_Q$).}
		
We then apply \eqref{com:H} to obtain
		\begin{align*}
			\|v_e^i P_{t}^\nu F_1\|_{L^2(\nu)}\leq \| v_e^i F_1\|_{L^2(\nu)}+\int_0^t\sum_{\bar e}c_{e\bar e}\,e^{-2c_Nd(e,\bar e)}\|\nabla_{\bar e}P_s^\nu F_1\|_{L^2(\nu)}\dif s,
		\end{align*}
		which implies that
		\begin{align}\label{nab:1}
			\|\nabla_e  P_{t}^\nu F_1\|_{L^2(\nu)}\leq d(\mfg)\|\nabla_e F_1\|_{L^2(\nu)}+\int_0^t\sum_{\bar e}D_{e\bar e}e^{-2c_Nd(e,\bar e)}\|\nabla_{\bar e}P_s^\nu F_1\|_{L^2(\nu)}\dif s,
		\end{align}
		for $D_{e\bar e}=d(\mfg)c_{e\bar e}$.	
This is almost the desired bound of the type \eqref{e:idea-comm}, except 
that we eventually require a bound in terms of the norm of $F$, not $F_1$,
and also we want the first term on the RHS to be small when $e$ is far from $\Lambda_F$.

To this end,		
by direct calculation as in \eqref{nab:e}, one has
		\begin{align*}
			\nabla_e F_1
			=\E_{\mu_Q}(\nabla_eF)-\kappa N\cov_{ \mu_Q}(F,Q_e(\Phi_y\Phi_x^t)Q_e-\Phi_x\Phi_y^t).
		\end{align*}
By  Lemma \ref{lem:co} we have
		\begin{align*}
			|\cov_{ \mu_Q}(F,Q_e(\Phi_y\Phi_x^t)Q_e-\Phi_x\Phi_y^t)|\leq C_1e^{-2c_Nd(\Lambda_F,e)}\$F\$_{2,Q},
		\end{align*}
		for $C_1$ depending on $\Lambda_F$ and universal $c_N$ independent of $F$.
		Hence, we have
		\begin{align*}
			\|\nabla_e F_1\|_{L^2(\nu)}\leq \|\nabla_eF\|_{L^2(\mu_\Lambda)}+C_1e^{-2c_Nd(\Lambda_F,e)}\$F\$_{2,\Lambda}\leq C_1e^{-2c_Nd(\Lambda_F,e)}\$F\$_{2,\Lambda},
		\end{align*}
		where we  adjust $C_1$ to have the second inequality.
		Substituting this into \eqref{nab:1}  we get
		\begin{align}\label{ena}
			\|\nabla_e  P_{t}^\nu F_1\|_{L^2(\nu)} &\leq C_1e^{-2c_Nd(\Lambda_F,e)}\$F\$_{2,\Lambda}+\int_0^t\sum_{\bar e\in E_\Lambda^+}D_{e\bar e}e^{-2c_Nd(e,\bar e)}\|\nabla_{\bar e}P_s^\nu F_1\|_{L^2(\nu)}\dif s.
		\end{align}
	We then use \eqref{ena} again to bound  $\|\nabla_{\bar e}P_s^\nu F_1\|_{L^2(\nu)}$ and obtain
	\begin{align*}
		\|\nabla_e  P_{t}^\nu F_1\|_{L^2(\nu)}&\leq C_1e^{-2c_Nd(\Lambda_F,e)}\$F\$_{2,\Lambda}+\sum_{\bar e}D_{e\bar e}e^{-2c_Nd(e,\bar e)}C_1\int_0^te^{-2c_Nd(\bar e,\Lambda_F)}\$F\$_{2,\Lambda}\dif s+J
		\\&\leq C_1e^{-c_Nd(\Lambda_F,e)}(1+C_0t)\$F\$_{2,\Lambda}+J,
	\end{align*}
where $C_0\geq\sum_{\bar e}e^{-c_Nd(e,\bar e)}D_{e\bar e}$ is finite and can be chosen to be independent of $e$,
$$J\eqdef\int_0^t\sum_{\bar e\in E_\Lambda^+}D_{e\bar e}e^{-2c_Nd(e,\bar e)}\int_0^s\sum_{\tilde e}D_{\bar e\tilde e}e^{-2c_Nd(\bar e,\tilde e)}\|\nabla_{\tilde e}P_r^\nu F_1\|_{L^2(\nu)}\dif s\dif r,$$
 and we used $d(\Lambda_F,\bar e)+d(\bar e,e)\geq d(\Lambda_F,e)$ in the last step.

	Iterating again we finally obtain
		\begin{align*}
		\|\nabla_e  P_{t}^\nu F_1\|_{L^2(\nu)}
		&\leq C_1e^{-c_Nd(\Lambda_F,e)}\sum_{n=0}^\infty \frac{(C_0t)^n}{n!}\$F\$_{2,\Lambda}=C_1e^{C_0t-c_Nd(\Lambda_F,e)}\$F\$_{2,\Lambda}.
	\end{align*}
		 Hence, \eqref{nabF} and the theorem follow.
	\end{proof}

\br\label{rem:1-dep-Phi}
Without the condition $K_{\mR^N}^\nu>0$, 
 the above theorem still holds if  
one of the two observables, say $F$, only depends on the Higgs field $\Phi$.
In this case   \eqref{cov1} becomes  $\cov_{\Lambda}(F,H)=\E_\nu\cov_{\mu_Q}(F,H)$ since $\E_{\mu_Q}F$ is a constant, and Lemma \ref{lem:co} suffices to prove the result.
\er

\subsection{Mass gap for $M\in \{ \mS^{N-1}, G\}$}
\label{sec:massGS}
	
In this section we prove the existence of a mass gap for the two cases where the Higgs fields take values in $M\in \{G,\mS^{N-1}\}$.
Since  $M$ are now compact, i.e. $\Phi$ is bounded, the general proof strategy is similar as \cite[Section 4.3]{SZZ22}, which was discussed below \eqref{e:gen-mass};
however, the estimates in the $\mS^{N-1}$ case are harder for geometric reasons (e.g.  $\mS^{N-1}$ is only a manifold, not a Lie group).

As alluded above \eqref{e:idea-comm} the key of the proof is to estimate 
the commutators between the derivatives and the generator. 
Such commutator estimates will be proved in Lemma~\ref{lem:cG} and \eqref{com:1},
and the desired type of bound \eqref{e:idea-comm}
will be realized in \eqref{e:DxPF} \eqref{e:DePF}.
	
Writing $\cS=\cS_{\YMH}$ for short,	
the generator has the following form
	\begin{align}\label{eq:L}
		\cL_\Lambda F=	\sum_{e\in E_{\Lambda}^+}\Delta_{e}F+\sum_{x\in {\Lambda}}\Delta_{x}F+\sum_{e\in E_{\Lambda}^+}\<\nabla_e \cS(Q),\nabla_{e} F\>+\sum_{x\in {\Lambda}}\<\nabla_x \cS(\Phi),\nabla_{x} F\>.
	\end{align}
Recall the notations such as $\bar e\sim e$, $x\sim y$ in
Section~\ref{sec:YMHSDEs}.
We will also sometimes write
$x\sim e$  if $x\in e$.

 Recall the basis vectors $(v_e^i)$ introduced above, as well as $(v_x^i)$ in the case $M=G$ (in which case it is the same as $(v_e^i)$).
As in the proof of Theorem~\ref{co:mass},
the key step is to estimate the commutator $[v_e^i,\cL_\Lambda]$, $[v_x^i,\cL_\Lambda]$
with $\cL_\Lambda $ as in \eqref{eq:L}.
The computation for $v_x^i$ in the case $M=\mS^{N-1}$ is more involved so we postpone it to the end of this subsection.

%

\bl\label{lem:cG}
	There exist positive constants $a_{e\bar e}$, $a_{ex}$, $a_{xe}, a_{xy}$ such that 
	for $M\in \{G, \mS^{N-1}\}$,
	\begin{align}\label{e:1}
		|[v_e^i,\cL_\Lambda]F|\leq \sum_{\bar e\sim e}a_{e\bar e}|\nabla_{\bar e}F|+\sum_{x\sim e}a_{ex}|\nabla_{x}F|,\end{align}
	and for $M= G$
	\begin{align}\label{e:2}
		|[v_x^i,\cL_\Lambda]F|\leq \sum_{e\sim x}a_{xe}|\nabla_{e}F|+\sum_{x\sim y}a_{xy}|\nabla_{y}F|.\end{align}
	 	Moreover, in \eqref{e:1}, $a_{e\bar e}= N|\beta|\sqrt{d(\mfg)}$ for $e\neq \bar e$,
	$$
	a_{ee}=((d-1)|\beta|+|\kappa|) \Big(2 N\sqrt{d(\mfg)}+2\sqrt2N^{\frac32}\Big),
	$$
	and
	$$
	a_{ex}= \begin{cases}
		 2|\kappa| N\sqrt{d(\mfg)}& (M=G)\\
		 2|\kappa| N\sqrt{N-1} & (M= \mS^{N-1}).
		 \end{cases}
	$$
	Also, in \eqref{e:2}, $a_{xe}= 2|\kappa| N\sqrt{d(\mfg)}$ and
	\begin{equ}
	a_{xy}=2|\kappa| N\sqrt{d(\mfg)}\quad (x\neq y),\qquad
	a_{xx}= 4|\kappa| dN\sqrt{d(\mfg)}+4\sqrt2|\kappa|N^{\frac32}d.
	\end{equ}
	\el
	
	\begin{proof}
	We first prove \eqref{e:1}.
	As discussed above \eqref{com-1}, analogous calculations were performed in 
	\cite[Lemma 4.10]{SZZ22} for pure Yang--Mills.
In our case we have
		\begin{align*}
			[v_e^i,\cL_\Lambda]F
			&=v_e^i \cL_\Lambda F-\cL_\Lambda v_e^i F
			\\
			&=\sum_{\bar e\sim e}  \big\<\nabla_{v_e^i} \nabla_{\bar e}\cS,\nabla_{\bar e} F \big\>
			+\big \<\nabla_e\cS \;, \; \nabla_{v_e^i}\nabla_{e}F-\nabla_{e}v_e^iF \big\> \;
			\\&\qquad+\sum_{x\sim e}  \big\<\nabla_{v_e^i} \nabla_{x}\cS,\nabla_{x} F\big\>\eqdef \sum_{i=1}^3I_i \;,
		\end{align*}
where the last term $I_3$ arises 
due to the Higgs field on the vertices $x$.
The first two terms ($I_1$ and $I_2$) are the same 
as the first equation in the proof of  \cite[Lemma 4.10]{SZZ22}, except that $\cS=\cS_1-\cS_2$ now also has a Higgs term $\cS_2$,
and  the corresponding part for $\cS_1$ has been estimated as in \cite[Lemma 4.10]{SZZ22}.  Recall that each term in $\cS_2$  only depends on one edge. Hence, 	to estimate $I_1+I_2$,  we only need to calculate $|v_e^i v_{e}^j \cS|$ and $|\nabla_e \cS|$. Other terms  satisfy the same bounds as in \cite[Lemma 4.10]{SZZ22}.  By direct calculation we have
		\begin{align*}
			|v_e^i v_{e}^j \cS|\leq 2(d-1)N|\beta|+2|\kappa| N,
		\end{align*}
		and by Lemma \ref{lem:nab}
		\begin{align*}
			\nabla_e \cS=
				-\frac12N\beta\sum_{p\in \cP_{\Lambda}, p\succ e}(Q_p-Q_p^t)Q_e-\kappa N (Q_e\Phi_y\Phi_x^tQ_e-\Phi_x\Phi_y^t).
		\end{align*}
		Hence, we have
		$$|\nabla_e\cS|\leq 2(d-1)N^{3/2}|\beta|+2N^{3/2}|\kappa|.$$
		Consequently, by a similar argument as in \cite[Lemma 4.10]{SZZ22} we have
		\begin{equs}\label{bd:I1I2}
			|I_1+I_2| &\leq N|\beta| \sqrt{d(\mfg)} \sum_{\bar e\sim e,\bar e\neq e}|\nabla_{\bar e}F|
			+2\Big((d-1)N|\beta|+N|\kappa|\Big) \sqrt{d(\mfg)} |\nabla_{ e}F|\\
			&\qquad + 2\sqrt2\Big((d-1)N^{3/2}|\beta|+N^{3/2}|\kappa|\Big)|\nabla_eF|\;.
		\end{equs}
		  
Now we consider the term $I_3$.  We have
$|I_3|\le \sum_{x\sim e}  |\nabla_{v_e^i} \nabla_{x}\cS | \, |\nabla_{x} F|$.
To compute the gradient $\nabla_{x}\cS$, let $v$ be a unit vector 
in $T_{\Phi_x} M$.
		For $M= \mS^{N-1}$ we have
		\begin{align*}
		|v_e^i v \cS|
		=2|\kappa| N\, \Big|v^t X_e^iQ_e\Phi_y \Big|\leq 2|\kappa| N|\Phi_y|\leq 2|\kappa| N.
		\end{align*}
		For $M= G$ we have
		\begin{align}\label{exS}
			|v_e^i v \cS|=2|\kappa| N \,
			\Big|\tr(X_e^iQ_e\Phi_y( v \Phi_x)^t) \Big|
			\leq 2|\kappa| N.
		\end{align}
So $ |\nabla_{v_e^i} \nabla_{x}\cS | \le 2|\kappa| N \dim(M)$ and we get
		\begin{align}\label{bd:I3}
			|I_3|\leq \sum_{x\sim e}a_{ex}|\nabla_xF |.
		\end{align}
		Combining \eqref{bd:I1I2} \eqref{bd:I3} we obtain \eqref{e:1}.

\bigskip

	Now we prove \eqref{e:2}. 	Here $M= G$.	Again $v_x^i$ are right-invariant fields and thus commutes with the Beltrami-Laplacian $\Delta_x$. So similarly as before we have
		\begin{equs}[e:J123]
			{}[v_x^i,\cL_\Lambda]F
			&=v_x^i \cL_\Lambda F-\cL_\Lambda v_x^i F
			\\
			&=\sum_{x\sim y}  \big\<\nabla_{v_x^i} \nabla_{ y}\cS,\nabla_{ y} F \big\>
			+\big \<\nabla_x\cS \;, \; \nabla_{v_x^i}\nabla_{x}F-\nabla_{x}v_x^iF \big\> \;
			\\&\qquad+\sum_{x\sim e}  \big\<\nabla_{v_x^i} \nabla_{e}\cS,\nabla_{e} F \big\>\eqdef \sum_{i=1}^3J_i \;,
		\end{equs}
			where we used $v_x^i\<\nabla_y\cS,\nabla_yF\>=\<\nabla_{v_x^i}\nabla_y\cS,\nabla_yF\>+\<\nabla_y\cS,\nabla_{v_x^i}\nabla_yF\>$ and $\nabla_{v_x^i}\nabla_yF=\nabla_yv_x^iF$ for $x\neq y$.
			
To estimate $J_1+J_2$,
since for each $x$ there are $2d$ edges containing $x$, we have
		\begin{align*}
			|v_x^iv_x^j\cS| & =2|\kappa| N\Big|\sum_{e=(x,y)}\tr((X_x^iX_x^j\Phi_x)^tQ_e\Phi_y)\Big|\leq 4|\kappa| dN,
\\
			|v_x^iv_y^j\cS| & =2|\kappa| N\Big|\tr((X_x^i\Phi_x)^tQ_eX_y^j\Phi_y)\Big|\leq 2|\kappa| N,
\\			
|\nabla_x\cS| & \leq 4d|\kappa| N^{3/2},
		\end{align*}
where in the last bound we used Lemma \ref{lem:nab_x}. 
	
		Combining these estimates and by a similar argument as in \cite[Lemma 4.10]{SZZ22}, we have
		\begin{align}\label{bd:J1}
			|J_1+J_2|\leq 2|\kappa| N\sqrt{d(\mfg)}\sum_{y\sim x,y\neq x}|\nabla_yF|+4d|\kappa| N\sqrt{d(\mfg)}|\nabla_xF|+4\sqrt2|\kappa| d N^{3/2}|\nabla_xF|.
		\end{align}
		We also have
$J_3 
=\sum_{x\sim e}\sum_j  ( v_x^iv_e^j\cS) v_e^j F$.
		Hence, similarly to \eqref{exS} we get
		\begin{align}\label{bd:J3}
			|J_3|\leq \sum_{x\sim e}\Big(\sum_j  | v_x^iv_e^j\cS|^2\Big)^{1/2}|\nabla_eF |\leq \sum_{x\sim e}2|\kappa| N\sqrt{d(\mfg)}|\nabla_eF |.
		\end{align}
		Combining \eqref{bd:J1} and \eqref{bd:J3} leads to \eqref{e:2}.
	\end{proof}
	
Based on \eqref{e:1} and \eqref{e:2} the proof of mass gap for $M= G$ follows similarly as in \cite{SZZ22}. 
Namely, these bounds lead to an inequality of the form \eqref{e:idea-comm}, with uniformly bounded constants $D_{e,\bar e}$, so the mass gap follows by iteration, as explained below  \eqref{e:idea-comm}. 

The case $M= \mS^{N-1}$ requires more calculations with  a suitable choice of basis $\{v_x^i\}$. Below we focus on this case.
Recall the constants $K_{\mS^{N-1}}$ defined  in \eqref{e:KS-S}
and  $K_G$ defined  in \eqref{e:KS-G}.
These are the constants appearing in the functional inequalities in Lemmas~\ref{lem:ba-S}
and \ref{lem:log-G}.
Recall that $K_{\mS^{N-1}}>0$ or $K_G>0$
is equivalent with the small $(\beta,\kappa)$ condition \eqref{e:K_YM}+\eqref{e:small-kappa}.

Also, recall Corollary~\ref{cor:SG-erg}
that infinite volume limit $\mu$ of $\mu_\Lambda$ is unique.
The following is the more precise statement of Theorem~\ref{main4} for $M\in \{G,\mS^{N-1}\}$.
	
	\bc\label{c:ex}
	Suppose  $K_{\mS^{N-1}}>0$ if $M= \mS^{N-1}$ 
	or  $K_G>0$ if $M= G$. 
	For $F, H\in C^\infty_{cyl}(\cQ)$, suppose that $\Lambda_F\cap \Lambda_H=\emptyset$ with $\Lambda_F, \Lambda_H$ defined in \eqref{def:LambdaF}.
	Then 
	\begin{align*}
		|\cov_\mu(F,H)|\leq c_{1} d(\mfg) e^{-c_N d(\Lambda_F,\Lambda_H)}(\$ F\$_{\infty}\$ H \$_{\infty}+\|F\|_{L^2}\|H\|_{L^2}),
	\end{align*}
	where $c_{1}$ depends on $|\Lambda_F|$, $|\Lambda_H|$,
	and $c_N$ depends on $N$, $d$, $K_{\mS^{N-1}}$ or $K_G$. 
	Here
	$$\$F\$_\infty\eqdef \sum_{e}\|\nabla_eF\|_{L^\infty}+\sum_{x}\|\nabla_xF\|_{L^\infty}.$$
	\ec

	\begin{proof}
		The proof for the case $M=G$ follows by the same arguments as in \cite[Corollary 4.11]{SZZ22} using Lemma \ref{lem:cG}. In the following we focus on the case
		$M= \mS^{N-1}$.
		As before (e.g. \eqref{e:gen-mass}, \eqref{cov4}) we have
		\begin{equ}
			|\cov_{\mu_\Lambda}(F,H)|
			\leq
			|\E_{\mu_\Lambda}(P_t^\Lambda(FH)-P_t^\Lambda F P_t^\Lambda H)|+\var_{\mu_\Lambda}(P_t^\Lambda F)^{1/2}\var_{\mu_\Lambda}(P_t^\Lambda H)^{1/2},
		\end{equ}
	with $(P_t^\Lambda )$ the semigroup of the dynamics as in Section \ref{sec:SDESG}.
	By the Poincar\'e inequality
		$$\var_{\mu_\Lambda}(P_t^\Lambda F)\leq e^{-2tK_{\mS^{N-1}}}\|F\|_{L^2(\mu_{\Lambda})}^2.$$
Therefore, as discussed below \eqref{e:gen-mass} (and recalling the form of the generator \eqref{eq:L}), it remains to estimate 
	\begin{equs}
		P_t(FH)-P_tF P_tH
		&=\int_0^t\frac{\dif}{\dif s}[P_s(P_{t-s}FP_{t-s}H)]\dif s\\
		&=2\int_0^tP_s\Big(\sum_x\<\nabla_x P_{t-s}F,\nabla_x P_{t-s}H\>+\sum_e\<\nabla_e P_{t-s}F,\nabla_e P_{t-s}H\>\Big)\dif s.
	\end{equs}
		Here and in the sequel we omit the superscript $\Lambda$ for notation's simplicity. 
		In the following we will prove $L^\infty$ bounds for $|\<\nabla_x P_{t-s}F,\nabla_x P_{t-s}H\>|$ and $|\<\nabla_e P_{t-s}F,\nabla_e P_{t-s}H\>|$.
		
		To this end,  for each {\it fixed} $\Phi_x\in \mS^{N-1}$  we choose $N-1$ vector fields $\{v_x^1,\cdots, v_x^{N-1}\}$, which form an orthonormal basis in tangent space at $\Phi_x$ and each $v_x^i$ is a Killing vector field.
		Note that  $\{v_x^i\}$ does not need to form an orthonormal basis at points of $\mS^{N-1}$ other than $\Phi_x$.
	Here the advantage for choosing Killing vector fields is that they commute with the Beltrami--Laplace operator \cite[Exercise~8.4.14]{petersen}, which makes it easier  to calculate $[v_x^i,\cL]$. 
		Without loss of generality we  assume that $\Phi_x=(1,0,\dots,0)\in \mS^{N-1}\subset \mR^N$ and
		we choose the  $N-1$ Killing vector fields
		$\{v_x^i : i=1,\cdots,N-1\}$
		at  a generic point $\Phi=(\phi_1,\dots,\phi_N)\in \mS^{N-1}\subset \mR^N$ as
		$$
		(-\phi_2,\phi_1,0,\cdots,0), \quad (-\phi_3,0,\phi_1,0,\cdots,0),\quad \cdots, \quad (-\phi_N,0,\cdots, 0,\phi_1).
		$$
		At the point $\Phi_x=(1,0,\cdots,0)$, they are precisely
		$$
		(0,1,0,\cdots,0), \quad \cdots,  \quad (0,\cdots,0,1),
		$$
		which form an  orthonormal basis at $\Phi_x$ as desired
		and we have
		\begin{align}\label{killvec}
			\sup_{\Phi\in \mS^{N-1}}\sup_{i}|v_x^i(\Phi)|\leq 1.
		\end{align}
		
		For the rest of the proof, we suppose that  for any  $F\in C^\infty_{cyl}(\cQ)$, there exist positive constants $a_{xy}, a_{xe}$, which are independent of $x,y, e$, such that
		\begin{align}\label{com:1}
			\|[v_x^i,\cL]F\|_{L^\infty}\leq \sum_{y\sim x}a_{xy}\|\nabla_yF\|_{L^\infty}+\sum_{e\sim x}a_{xe}\|\nabla_eF\|_{L^\infty}.
		\end{align}
		This will be justified later in Lemma~\ref{lem:com}.
		
		We have
		\begin{align}\label{zmm1}
			v_x^iP_tF-P_tv_x^iF=\int_0^t\frac{\dif}{\dif s} \left( P_{t-s}v_x^iP_sF\right)\dif s=\int_0^tP_{t-s}[v_x^i,\cL]P_sF \dif s.
		\end{align}
		Using \eqref{com:1} we then have
		\begin{align*}
			&|\nabla_x P_{t}F|(\Phi_x)\leq \sum_i|v_x^i P_tF|(\Phi_x)
			\\&\leq 
			\sum_i |P_tv_x^iF|(\Phi_x) 
			+(N-1)\int_0^t\Big(\sum_{y\sim x}a_{xy}\|\nabla_y P_sF\|_{L^\infty}+\sum_{e\sim x}a_{xe}\|\nabla_eP_sF\|_{L^\infty}\Big)\dif s.
		\end{align*}
Since $|P_tv_x^iF|(\Phi_x) \le \|\nabla_xF\|_{L^\infty}$,
we  take sup w.r.t. $\Phi_x$ and obtain
\begin{equ}[e:DxPF]
	\|\nabla_x P_{t}F\|_{L^\infty}\leq  (N-1)\|\nabla_xF\|_{L^\infty} +(N-1)\int_0^t\Big(\sum_{y\sim x}a_{xy}\|\nabla_y P_sF\|_{L^\infty}+\sum_{e\sim x}a_{xe}\|\nabla_e P_sF\|_{L^\infty}\Big)\dif s.
\end{equ}
		On the other hand using \eqref{e:1} of Lemma~\ref{lem:cG} we obtain
\begin{equ}[e:DePF]
			\|\nabla_e P_tF\|_{L^\infty}\leq d(\mfg)\|\nabla_e F\|_{L^\infty}+\int_0^t\Big(\sum_{\bar e} D_{e\bar e} \|\nabla_{\bar e} P_sF\|_{L^\infty}+\sum_{x} D_{ex} \|\nabla_{x} P_sF\|_{L^\infty}\Big)\dif s,
\end{equ}
with a matrix $D$ such that 
\begin{itemize}
\item
$D_{e\bar e}= d(\mfg) a_{e\bar e}$ if $e\sim \bar e$ and $D_{e\bar e}=0$ otherwise
\item
  $D_{ex}= d(\mfg) a_{ex}$ if $x\in e$ and $D_{ex}=0$ otherwise,
 \end{itemize} 
 where $a_{e\bar e}$ and $a_{ex}$ are introduced in Lemma \ref{lem:cG}. 
We can also write \eqref{e:DxPF} in terms of 
\begin{itemize}
\item
 $D_{xy}=(N-1)a_{xy}$ if $x\sim y$ and $D_{xy}=0$ otherwise
\item
 $D_{xe}=(N-1)a_{xe}$ if $x\in e$ and $D_{xe}=0$ otherwise.
  \end{itemize}
		Since $x\notin \Lambda_F$ we get $\nabla_x F=0$,   by iteration we arrive at
		\begin{align*}
			\|\nabla_x P_tF\|_{L^\infty}\leq d(\mfg)\sum_{n=N_x}^\infty \frac{t^n}{n!}\Big(\sum_{e}D_{xe}^{(n)} \|\nabla_e F\|_{L^\infty}+\sum_{y}D_{xy}^{(n)} \|\nabla_y F\|_{L^\infty}\Big),
		\end{align*}
		with $N_x=d(x,\Lambda_F)$ and $D_{xe}^{(n)}\vee D_{xy}^{(n)}\leq C^n_0$, where
		\begin{equ}[e:C0]
		C_0=d(\mfg) \Big(a_{xx}+2da_{xy}+2da_{xe}+a_{ee}+6(d-1)a_{e\bar e}+2a_{ex}\Big).
		\end{equ}
		As a result, using $n!\geq e^{n\log n-2n}$, for any $c>0$ with $B>0$ satisfying $2-\log B+\log C_0+\frac{C_0}{B}\leq -2c$ and $d(x,\Lambda_F)\geq Bt$, we have
		\begin{align}\label{na:x}
			\|\nabla_x P_tF\|_{L^\infty}
			\leq\sum_{n=N_x}^\infty \frac{t^n}{n!}C_0^n d(\mfg) \$ F \$_{\infty}
			\leq \frac{(C_0t)^{N_x}}{N_x!}d(\mfg)e^{tC_0} \$ F \$_{\infty}
			\leq  d(\mfg)e^{-2c\, d(x,\Lambda_F)}\$ F \$_{\infty}.
		\end{align}
		Similarly, we get for $d(e,\Lambda_F)\geq Bt$,
		\begin{align}\label{na:e}
			\|\nabla_e P_tF\|_{L^\infty}
			\leq  d(\mfg)e^{-2c \,d(e,\Lambda_F)}\$ F \$_{\infty}.
		\end{align}
		In the following we choose $t= d(\Lambda_F,\Lambda_H)/(2B)$. For $x\in \Lambda_H$ we use \eqref{na:x} and \cite[Theorem 5.6.1]{Wang} (see also (4.9) in \cite{SZZ22}) to have
		\begin{equs}
			\|\<\nabla_xP_{t-s}F,\nabla_xP_{t-s}H\>\|_{L^\infty}
			&\leq \|\nabla_xP_{t-s}F\|_{L^\infty}\|\nabla_xP_{t-s}H\|_{L^\infty}
			\\
			&\leq c_1e^{-c\, d(\Lambda_F,\Lambda_H)-c\,d(\Lambda_F,x)}\$ F \$_{\infty}\$H\$_{\infty}
		\end{equs}
	for some $c_1>0$, which may differ from line to line.
		The analogous bound also holds for $x\in \Lambda_F$,
		and  for $x\notin \Lambda_F\cup \Lambda_H$, 
	as in \eqref{e:three-exp} we get a factor 
	$e^{-\frac{c}2\, d(\Lambda_F,\Lambda_H)-\frac{c}2\,(d(x,\Lambda_F)\wedge d(x,\Lambda_H))}$.
%
Similar bound holds for $	\|\<\nabla_eP_{t-s}F,\nabla_eP_{t-s}H\>\|_{L^\infty}$.
		Combining the above calculations we obtain 
		\begin{align*}
			\|P_t(FH)-P_tFP_tH\|_{L^\infty}\leq c_1e^{-\frac{c}4\, d(\Lambda_F,\Lambda_H)}\$ F \$_{\infty}\$H\$_{\infty}.
		\end{align*}
		Hence the result follows.
	\end{proof}
	
	\br
	From the above proof one can see that $c_N \gtrsim K_M / C_0$
	where $K_M\in \{K_{\mS^{N-1}}, G\}$ and $C_0$ is as in \eqref{e:C0},
	but this is not necessarily optimal.
	More explicitly for $M= G$,
	$$
	c_N \gtrsim \frac{K_G}{d(\mfg)(d|\kappa| N\sqrt{d(\mfg)}+dN\sqrt{d(\mfg)}|\beta|)},
	$$
	and from Lemma \ref{lem:com} below for $M= \mS^{N-1}$
	$$
	c_N \gtrsim \frac{K_{\mS^{N-1}}}{d(\mfg)(|\kappa| 2^{N} N^{3}+d|\kappa| N\sqrt{d(\mfg)}+dN\sqrt{d(\mfg)}|\beta|)}.
	$$
	\er
	
	In the following we prove \eqref{com:1}. 
	This will span Lemma~\ref{lem:w} -- Lemma~\ref{lem:com}.
	On the Lie group $G$, the right-invariant vector fields can be chosen to commute with the Beltrami--Laplacian and at the same time be orthonormal at every point. The second property is convenient in obtaining uniform bounds on various gradient terms.
	On the sphere, however, these two properties can not be achieved simultaneously,
	which makes \eqref{com:1} harder to prove.
	
	To this end we  choose suitable orthonormal vector fields $\{ w_x^j\}$ for a suitable domain in $\mS^{N-1}$, different from the Killing fields $\{v_x^i\}$ in  Corollary \ref{c:ex}.
	More precisely we use spherical coordinates $(\theta_1,\dots,\theta_{N-1})$ and restrict us for the construction in the domain
	$D=\{\theta_i\in [\pi/4,3\pi/4], i=1,\dots,N-1\}\subset \mS^{N-1}$. 
	For other regions of $ \mS^{N-1}$ one could simply alter the range of $\theta_i$ or select alternative local coordinates.

Since  \eqref{com:1} is an $L^\infty$ estimate, we first prove the following.
	
	\bl\label{lem:w} 	
	There exist vector fields $\{w^j\}_{j=1}^{N-1}$ on $D$, which are orthonormal  at every point of $D$ 
	(i.e., for every $\phi\in D$, the vectors $\{w^j(\phi)\}_{j=1}^{N-1}$ form an orthonormal basis of $T_{\phi}\mS^{N-1}$),  and 
	\begin{equs}[bd:c]
		\|\nabla_{w^i}w^i&\|_{L^\infty(D)}\leq (i-1)2^{i-1},\qquad
		\|\nabla_{w^i}w^j\|_{L^\infty(D)}\leq 2^{\frac{j-1}2}\quad (i>j),
		\\
		&\|\nabla_{w^i}w^j\|_{L^\infty(D)}\leq 2^{\frac{i+1}2} \quad (i< j).
	\end{equs}
	\el

\begin{proof}
We represent  $\phi=(\phi_1,\dots,\phi_N)\in \mS^{N-1}$ in spherical coordinates $(\theta_1,\cdots,\theta_{N-1})$ as
\begin{equs}[e:spherical]	
\phi_i =  r \Big(\prod_{k=1}^{i-1} \sin\theta_k\Big)
		\cos\theta_i, \quad (1\le i \le N-1), \qquad
\phi_{N} =  r \prod_{k=1}^{N-1} \sin\theta_k
\end{equs}
with $r=1$, with $\prod_{i=1}^0=1$ as usual.
		We then have
		$$\dif s^2 =  \dif\theta_1^2 +  \sin^2\theta_1 \dif\theta_2^2 +  \sin^2\theta_1 \sin^2\theta_2 \dif\theta_3^2 + \dots +  \sin^2\theta_1 \sin^2\theta_2\dots \sin^2\theta_{N-2} \dif\theta_{N-1}^2,$$
	 and the vector fields $\{\p_{\theta_i}\}_{i=1}^{N-1}$ are orthogonal: 
$\<\p_{\theta_i},\p_{\theta_j}\>=\delta_{i,j} \prod_{j=1}^{i-1}\sin^2\theta_j$.
		By direct calculation, the Christoffel symbols for $\{\theta_i\}$ satisfy that  for $k<i$
		\begin{equ}[e:Christ]
		0< |\Gamma_{\theta_i\theta_i}^{\theta_k}|\leq 1,\quad 0< |\Gamma_{\theta_i\theta_k}^{\theta_i}|\leq 1,\quad 0< |\Gamma_{\theta_k\theta_i}^{\theta_i}|\leq1.
		\end{equ}
		Other Christoffel symbols are zero.
		The  desired orthonormal basis 
		is given by
		$w^j=\partial_{\theta_j} /\prod_{i=1}^{j-1}\sin \theta_i$ for $1\le j\le N-1$.
		We then have
		\begin{align*}
			\nabla_{w^i}w^j=\Big(\prod_{k=1}^{i-1}\sin\theta_k\Big)^{-1}\Big[\frac1{\prod_{k=1}^{j-1}\sin\theta_k}\sum_l\Gamma^{\theta_l}_{\theta_i\theta_j}\p_{\theta_l}+\p_{\theta_i}\Big(\frac1{\prod_{k=1}^{j-1}\sin\theta_k}\Big)\partial_{\theta_j}\Big].
		\end{align*}
		By straightforward calculation and using  \eqref{e:Christ} we find that \eqref{bd:c} holds.
	\end{proof}


	\bl\label{lem:vw}  
	For the Killing fields $\{v^i\}_{i=1}^{N-1}$ as in Corollary \ref{c:ex} and the orthonormal fields $\{w^j\}_{j=1}^{N-1}$ in Lemma \ref{lem:w}, there exist  smooth functions $a_{ij}$ on $D$ such that 
\begin{equ}\label{viw}
	v^i=\sum_{1\leq j\leq (i+1)\wedge (N-1)} \!\!\!\!\!\!\!\! a_{ij}w^j,
	\qquad 
	\|a_{ij}\|_{L^\infty(D)}\leq 2.
\end{equ}
\el

\begin{proof}
The proof is essentially elementary, but we include necessary details for completeness.
We write $(\phi_1,\dots,\phi_N)\in\mR^N$ in terms of spherical coordinates $(r,\theta_1,\cdots,\theta_{N-1})$ as in
\eqref{e:spherical}, now for $r>0$.
We first compute the Jacobian $\frac{\partial(r,\theta_1,\cdots, \theta_{N-1})}{\partial (\phi_1,\cdots, \phi_N)}$. Since $r=(\phi_1^2+\cdots +\phi_N^2)^{1/2}$, we have
$$
\frac{\partial r}{\partial \phi_i} = \frac{\phi_i}{r}
= \frac{\phi_i}{(\phi_1^2+\cdots +\phi_N^2)^{1/2}}.
$$
Since $\theta_1= \cos^{-1} (\phi_1/r)$, we obtain
\begin{equs}
\frac{\partial \theta_1}{\partial \phi_1}
=
-\frac{1}{\sin\theta_1}
\frac{\partial}{\partial \phi_1} (\phi_1/r)
=
-\frac{1}{\sin\theta_1} \frac{r^2-\phi_1^2}{r^3}
= -\frac{\sin\theta_1}{r},
\end{equs}
and similarly, for $2\le i\leq N-1$, we derive (recall $\prod_{k=2}^1=1$)
\begin{equs}
\frac{\partial \theta_1}{\partial \phi_i}
=
\frac{1}{\sin\theta_1}  \frac{\phi_1\phi_i}{r^3}
=
 \frac{\cos\theta_1 \cos \theta_i \prod_{k=2}^{i-1}\sin \theta_k}{r},\quad  \frac{\partial \theta_1}{\partial \phi_N}
 =\frac{1}{\sin\theta_1}  \frac{\phi_1\phi_N}{r^3}=
 \frac{\cos\theta_1  \prod_{k=2}^{N-1}\sin \theta_k}{r}.
\end{equs}
We then observe that $\theta_{N-1} =\tan^{-1}(\phi_N/\phi_{N-1})$ and
$$
\theta_2=\cos^{-1} \Big( \frac{\phi_2}{(r^2-\phi_1^2)^{1/2}}\Big),
\quad\cdots\quad,\theta_{N-2}=\cos^{-1} \Big( \frac{\phi_{N-2}}{(r^2-(\sum_{i=1}^{N-3}\phi_i^2))^{1/2}}\Big).
$$
From this we obviously have
$$
\frac{\partial \theta_k}{\partial \phi_m}=0,\qquad (k>m).
$$
By direct calculation, we also have
\begin{equs}
\frac{\partial \theta_2}{\partial \phi_2}
&=
\frac{\partial }{\partial \phi_2}
\cos^{-1} \Big( \frac{\phi_2}{(r^2-\phi_1^2)^{1/2}}\Big)
=
-\Big(
\frac{r^2-\phi_1^2}{r^2-\phi_1^2-\phi_2^2}
\Big)^{1/2}
\frac{\partial }{\partial \phi_2}
\Big( \frac{\phi_2^2}{r^2-\phi_1^2}\Big)^{1/2}
\\
&=
- \frac{(r^2-\phi_1^2-\phi_2^2)^{1/2}}{ r^2-\phi_1^2}
=
-\frac{r\sin\theta_1\sin\theta_2}{ (r \sin\theta_1)^2}
=
-\frac{\sin\theta_2}{ r \sin\theta_1}\;,
\end{equs}
and similarly for $1\leq i<N-1$
\begin{equs}
\frac{\partial \theta_i}{\partial \phi_i}
=
- \frac{(r^2-(\sum_{k=1}^i\phi_k^2))^{1/2}}{ r^2-(\sum_{k=1}^{i-1}\phi_k^2)}
=-\frac{r\prod_{k=1}^i\sin\theta_k}{r^2\prod_{k=1}^{i-1}\sin^2\theta_k}
=-\frac{\sin\theta_i}{r\prod_{k=1}^{i-1}\sin\theta_k},
\end{equs}
and for $i<j\leq N-1$
\begin{equs}
	\frac{\partial \theta_i}{\partial \phi_j}
	&=
	\frac1{\sin \theta_i} \frac{\phi_i\phi_j}{( r^2-(\sum_{k=1}^{i-1}\phi_k^2))^{3/2}}=\frac{\cos \theta_i\cos \theta_j\prod_{k=i+1}^{j-1}\sin\theta_k}{r\prod_{k=1}^{i-1}\sin\theta_k},\\	
	\frac{\partial \theta_i}{\partial \phi_N}
	&=	\frac1{\sin \theta_i} \frac{\phi_i\phi_N}{( r^2-(\sum_{k=1}^{i-1}\phi_k^2))^{3/2}}=\frac{\cos \theta_i \prod_{k=i+1}^{N-1}\sin\theta_k}{r\prod_{k=1}^{i-1}\sin\theta_k},
	\\
	\frac{\p\theta_{N-1}}{\partial \phi_{N-1}}
	&=-\frac{\phi_N}{\phi_{N-1}^2}\cos^2\theta_{N-1}=-\frac{\sin \theta_{N-1}}{r\prod_{k=1}^{N-2}\sin\theta_{k}},
	\\
	\frac{\p\theta_{N-1}}{\partial \phi_{N}}
	&=\frac{1}{\phi_{N-1}}\cos^2\theta_{N-1}=\frac{\cos \theta_{N-1}}{r\prod_{k=1}^{N-2}\sin\theta_{k}}.
\end{equs}
We have for $1\leq i\leq N$
\begin{align}\label{tra}
\partial_{\phi_i}
=\sum_{l=1}^{N-1}
 \frac{\partial \theta_l}{\partial \phi_i} \partial_{\theta_l}
+ \frac{\partial r}{\partial \phi_i} \partial_r
=
\sum_{l=1}^{N-1}
\frac{\partial \theta_l}{\partial \phi_i} \partial_{\theta_l}
+ \frac{\phi_i}{r} \partial_r.
\end{align}
So, recalling $r=1$ on the unit sphere,
\begin{equs}
v^1 &= -\phi_2\partial_{\phi_1} + \phi_1 \partial_{\phi_2}
\\
&=
0\cdot  \partial_r
+( -\phi_2 \frac{\partial \theta_1}{\partial \phi_1}
 + \phi_1 \frac{\partial \theta_1}{\partial \phi_2})\partial_{\theta_1}
 +
( -\phi_2 \frac{\partial \theta_2}{\partial \phi_1}
 + \phi_1 \frac{\partial \theta_2}{\partial \phi_2})\partial_{\theta_2}
\\
&=
(\phi_2 \sin\theta_1 + \phi_1 \cos\theta_1\cos\theta_2 )\partial_{\theta_1}
+
( - \phi_1\frac{\sin\theta_2}{\sin\theta_1} )\partial_{\theta_2}
\end{equs}
which is a linear combination of $w^1$, $w^2$
with coefficients being polynomials
in $\sin$, $\cos$.
In a similar fashion, we substitute the expression for $\frac{\partial\theta_i}{\p \phi_j}$ into \eqref{tra} and replace $\p_{\theta_i}$ by $\prod_{k=1}^{i-1}\sin \theta_k w_i$ in \eqref{tra}. We note that the denominator in $\frac{\partial\theta_i}{\p \phi_j}, j\geq i$ is given by  $\prod_{k=1}^{i-1}\sin\theta_k$,  which conveniently cancels out. Hence, we obtain
\begin{align*}
	\p_{\phi_i}=\sum_{j=1}^{i}b_{i,j}w_j+\frac{\phi_i}{r}\p_r,
\end{align*}
with $b_{i,j}$ being product of $\sin$ and $\cos$ and $\|b_{i,j}\|_{L^\infty(D)}\leq 1$.
Hence, we have for $1\leq i\leq N-1$
\begin{align*}
	v^i=-\phi_{i+1}\p_{\phi_1}+\phi_1\p_{\phi_{i+1}}=-\phi_{i+1}b_{1,1}w_1+\phi_1\sum_{j=1}^{i+1}b_{i+1,j}w_j=\sum_{j=1}^{i+1}a_{ij}w_j,
\end{align*}
with $\|a_{ij}\|_{L^\infty(D)}\leq 2$. 
\end{proof}

\bl 
For each $x\in \Lambda$, let
 $\{w_x^j\}=\{w^j\}$ as in Lemma \ref{lem:w}. Then
\begin{align}\label{bd:w}
	w_x^iw_x^j\cS=0 \quad (\text{for } i<j),
	\qquad
	 |w_x^iw_x^j\cS|\leq 2^{N/2}N^{3/2}|\kappa|
	 \quad (\text{for } i\geq j).
\end{align}
\el
\begin{proof}
One has for $\Phi_x$ written in spherical coordinates
$$-\cS_2/(2\kappa N)=u^t\Phi_x=u_1\cos\theta_1+\sum_{i=2}^{N-1}u_i\Big(\prod_{k=1}^{i-1}\sin\theta_k\Big)\cos\theta_i+u_N\Big(\prod_{k=1}^{N-1}\sin\theta_k\Big),$$
for $u=Q_e\Phi_y$.
Hence,
$$w_x^1\cS/(2\kappa N)=-u_1\sin\theta_1+\sum_{i=2}^{N-1}u_i\cos\theta_1\Big(\prod_{k=2}^{i-1}\sin\theta_k\Big)\cos\theta_i+u_N\cos\theta_1\Big(\prod_{k=2}^{N-1}\sin\theta_k\Big),$$
and for $2\leq j\leq N-1$
\begin{align*}w_x^j\cS/(2\kappa N)&=\frac1{\prod_{i=1}^{j-1}\sin\theta_i}\Big[-u_j\Big(\prod_{k=1 }^{j-1}\sin\theta_k\Big)\sin\theta_j+\sum_{i=j+1}^{N-1}u_i\Big(\prod_{k=1,k\neq j}^{i-1}\sin\theta_k\Big)\cos\theta_j\cos\theta_i\\&\qquad\qquad+u_N\Big(\prod_{k=1,k\neq j}^{N-1}\sin\theta_k\Big)\cos\theta_j\Big]
	\\&= -u_j\sin\theta_j+\sum_{i=j+1}^{N-1}u_i\Big(\prod_{k=j+1}^{i-1}\sin\theta_k\Big)\cos\theta_j\cos\theta_i+u_N\Big(\prod_{k=j+1}^{N-1}\sin\theta_k\Big)\cos\theta_j.\end{align*}
By straightforward calculation we obtain $w_x^iw_x^j\cS=0 $ for  $i<j$ and
\begin{align*}
	\qquad |w_x^iw_x^j\cS|\leq 2^{N/2}N\sum_{k=1}^N|u_k||\kappa|\leq 2^{N/2}N^{3/2}|\kappa| \qquad \text{ for } i\geq j,
\end{align*}
where we used $|1/\sin\theta_i|\leq \sqrt 2$ on $D$ and H\"older's inequality.
\end{proof}

	\bl\label{lem:com}
The bound \eqref{com:1} holds 
 with
 $$a_{xe}=2|\kappa| N\sqrt{d(\mfg)},
 \qquad
 a_{xy}=2|\kappa| N \sqrt{N-1},
 \qquad
 a_{xx}=CN^{3}2^{N}|\kappa|
 $$
where $x\neq y$,  for some universal constant $C>0$.
	\el
	\begin{proof}
		Since $v_x^i$ commutes with the Beltrami-Laplacian $\Delta_x$, we also have \eqref{e:J123} as in the case $M=G$. 
		
Our aim is to estimate $L^\infty$-norm of each $J_i$ from \eqref{e:J123}. We will use the orthonormal basis $\{w_x^j\}$ constructed in Lemma \ref{lem:w} and restrict to the domain $D$ (introduced above Lemma \ref{lem:w}). For other domains we may  alter the domain of $\theta_i$ or select alternative local coordinates, and the calculation for each term follows in the same way.
		
		Since $\nabla_x f=\sum_jw_x^jf w_x^j$, we use $\nabla_{v_x^i}w_y^j=0$ for $x\neq y$ to  have
\begin{equs}
J_1 & = \sum_{x\sim y, y\neq x}\sum_jv_x^iw_y^j\cS \cdot w_y^jF+\sum_jv_x^iw_x^j\cS\cdot w_x^jF+\sum_jw_x^j\cS\cdot \nabla_{v_x^i}w_x^jF,
\\
J_2&=\sum_j\Big(v_x^iw_x^j F\cdot w_x^j\cS+w_x^jF \nabla_{v_x^i}w_x^j\cS-w_x^jv_x^iF \cdot w_x^j\cS\Big)
			\\&=\sum_{j}\Big([v_x^i,w_x^j]F\cdot w_x^j\cS+w_x^j F\cdot \nabla_{v_x^i}w_x^j\cS\Big),
\end{equs}
		and we use $\nabla_{v_x^i}v_e^j=0$ to have
		\begin{align*}
			J_3=\sum_{x\sim e}\sum_jv_x^iv^j_e\cS\cdot v_e^jF.
		\end{align*}

		In the following we estimate $J_1, J_2$ and $J_3$. For the terms that $x\neq y$ in $J_1$ and $J_3$ are easy since $\nabla_{v_x^i}v_e^j=0$ and $\nabla_{v_x^i}w_y^j=0$.
		
		We first use H\"older's inequality w.r.t. $j$ to have
		\begin{align}\label{est:j3}
			|J_3|\leq \sum_{x\sim e}\Big(\sum_j|v_x^iv^j_e\cS|^2\Big)^{1/2}\cdot |\nabla_eF|\leq 2|\kappa| N\sqrt{d(\mfg)}\sum_{x\sim e}|\nabla_eF|,
		\end{align}
		where we used $|v_x^iv^j_e\cS|\leq 2|\kappa| N$.
		
		For $J_1$ we write $v_x^i=\sum_ja_{ij}w_x^j$ for  smooth functions $a_{ij}$ by \eqref{viw} and we have
		\begin{align*}
			|J_1|\leq& \sum_{x\sim y,x\neq y}\Big(\sum_j|v_x^iw_y^j\cS|^2\Big)^{1/2} |\nabla_yF|+\Big(\sum_j|v_x^iw_x^j\cS|^2\Big)^{1/2} |\nabla_xF|+\sum_{j,l}|w_x^j\cS||a_{il}||\nabla_{w^l_x}w^j_x||\nabla_x F|
			\\\leq&\sum_{x\sim y,x\neq y}2|\kappa| N\sqrt{N-1} |\nabla_yF|+\Big(\sum_j\Big|\sum_l|a_{il}|\cdot|w_x^lw_x^j\cS|\Big|^2\Big)^{1/2} |\nabla_xF|+C_a8d|\kappa| N|\nabla_xF|,
		\end{align*}
		where  $$C_a^2=\sum_{j}\Big(\sum_l\|a_{il}\|_{L^\infty(D)}\|\nabla_{w^l}w^j\|_{L^\infty(D)}\Big)^2,$$
		and we used $|v_x^iw^j_y\cS|\leq 2|\kappa| N$ for the first term on the RHS and $|\nabla_x\cS|\leq 8d|\kappa| N$ by Lemma \ref{lem:nab_x} for the last term. Using \eqref{bd:c} we find $C_a\sim 2^{N}N^{3/2}$.
		
		Using \eqref{bd:w} and $\|a_{ij}\|_{L^\infty(D)}\leq 2$ by Lemma \ref{lem:vw},  we obtain
		\begin{align}\label{est:j1}
			|J_1|\leq \sum_{x\sim y,x\neq y}2|\kappa| N\sqrt{N-1} |\nabla_yF|+CN^{3}2^{N}|\kappa||\nabla_xF|,
		\end{align}
		for some universal constant $C>0$.
		
		
		For $J_2$, we  write
		\begin{align*}
			[v_x^i,w_x^j]=\sum_k(a_{ik}[w_x^k,w_x^j]-w_x^ja_{ik}\cdot w_x^k)
			=\sum_k\Big(a_{ik}(\nabla_{w_x^k}w_x^j-\nabla_{w_x^j}w_x^k)-w_x^ja_{ik}\cdot w_x^k\Big),
		\end{align*}
		which implies that
		\begin{align*}
			|J_2| &\leq \sum_{k,j} \|a_{ik}\|_{L^\infty(D)}(\|\nabla_{w_x^k}w_x^j\|_{L^\infty(D)}+\|\nabla_{w_x^j}w_x^k\|_{L^\infty(D)})|\nabla_xF|\cdot |w_x^j\cS|
			\\&\qquad +\sum_{j,k}|w_x^ja_{ik}|\cdot |w_x^kF|\cdot |w_x^j\cS|+C_a8|\kappa|d N|\nabla_xF|
			\\
			&\leq (2C_a+C_{a,1}+C_w)8|\kappa| d N|\nabla_xF|,
		\end{align*}
		where
		$$C_{a,1}^2=\sum_{j}\Big(\sum_k\|a_{ik}\|_{L^\infty(D)}\|\nabla_{w^j}w^k\|_{L^\infty(D)}\Big)^2,\qquad
		C_w^2=\sum_{j,k}\|w_x^ja_{ik}\|_{L^\infty(D)}^2,
		$$
		 and we bound $|\nabla_x\cS|$ by $8d|\kappa|  N$. Using \eqref{bd:c} we find $C_{a,1}\sim  2^{N}N^{3/2}$. Using the expression for $w^j$ and $a_{il}$ in Lemma \ref{lem:w} and Lemma \ref{lem:vw} we have $\|w_x^ja_{il}\|_{L^\infty(D)}\lesssim 2^{j/2}$, and  $C_w\leq 2^{N/2}N$.
		Hence, we have
		\begin{align}\label{est:j2}
			|J_2|\leq CN^{5/2}2^{N}|\kappa||\nabla_xF|,
		\end{align}
		for some universal constant $C>0$.
		Hence, combining \eqref{est:j3}, \eqref{est:j1}, \eqref{est:j2}, the result follows.
	\end{proof}

 \br\label{rem:Other-M}
 (i) Our proofs for the Poincar\'e inequality and mass gap also apply to other choices of the Higgs target spaces. In Section \ref{sec:Lattice YMH},
 we listed some additional options. Specifically, in Case \eqref{case4} in Section \ref{sec:Lattice YMH}, where $\Phi_x$ takes values in $G$ with matrix multiplication replaced by the adjoint representation, one can verify the Bakry--\'Emery conditions for small values of $\beta$ and $\kappa$. Consequently,  the log-Sobolev inequality,  mass gap, and  the uniqueness of the infinite volume limit also hold.

(ii) For the unbounded Cases  \eqref{case5} and  \eqref{case6}, wherein $\Phi_x$ assumes values in $N\times N$ skew-symmetric or symmetric matrices with the adjoint representation of $G$ on $M$, the disintegration method, previously applicable when $\Phi_x\in\mR^N$, remains effective. The Poincar\'e inequality and mass gap persist  for small values of $\beta$ and $\kappa$.

(iii) For $\Phi_x\in\mR^N$, we can also consider the case where $V(\Phi_x)=m |\Phi_x|^2+\lambda |\Phi_x|^4$, with $\lambda\geq0$ and $m\in\mathbb{R}$. When $m>0$, obtaining the Poincar\'e inequality and mass gap for small values of $\beta$ and $\kappa$ is straightforward, requiring only modifications to the proof of the mass gap for $\mu_Q$ (see Remark \ref{rem:V}). However, when $m<0$, additional conditions and calculations, as in \cite{Yosida} for $\Phi^4$ model, may be necessary to establish the mass gap of $\mu_Q$, which poses an interesting problem for future study.
 \er

\subsection{Improved mass gap condition for $M=G$ by gauge fixing}\label{sec:gauge fixing}

In this section we exploit the gauge symmetry to improve the conditions for mass gap of gauge invariant observables, in the case  $\Phi_x\in G$.
Recall \eqref{e:gauge} for the definition of gauge transformations.

Given any configuration $(Q,\Phi)$, since $\Phi_x\in G$, we can find a unique gauge transformation $g$
(in fact, $g_x=\Phi_x^{-1}$) such that
$$
g_x\Phi_x=I_N,
$$
with $I_N\in G$ being the identity matrix.
Such a gauge fixing that makes the Higgs field $\Phi$ into a fixed constant field
is usually called the U-gauge, c.f. \cite[Sec.~3]{Seiler}.
In this book the U-gauge fixing requires the so called ``complete breakdown of symmetry'' condition
i.e. the stability group of the constant $\Phi$ is trivial:
examples satisfying this condition includes 
$U(1)$, $SU(2)$ with defining representation, 
or $SO(3)$ with $\Phi$ consisting of a pair of fields in the vector representation,
see \cite[Theorem~3.18]{Seiler}.
In our case where the Higgs field is defined via the
representation on $G=SO(N)$ by left multiplication, the  stability group of $I_N$ is always trivial for any $N$.

To make our gauge fixing argument more precise, we follow the set-up in \cite[Appendix B]{MR3982691} by Driver.
 For fixed lattice $\Lambda$, we write
$$
\cQ=G^{E_{\Lambda}^+}\times G^{\Lambda}, \qquad
\sG=G^\Lambda
$$
for the manifold of field configurations and the gauge group, where we omit the subscript $\Lambda$ in $\cQ$ for simplicity.
For $(Q,\Phi)\in \cQ$ with $Q=(Q_e)_{e\in E^+_{\Lambda} }$ , $\Phi=(\Phi_x)_{x\in \Lambda}$, and $g=(g_x)_{x\in \Lambda}\in \sG$ we write
$$
g\cdot (Q,\Phi) = (g\cdot Q,g\cdot \Phi)\in \cQ
$$
for the gauge transformation, namely
\begin{align*}
	 (g\cdot Q)_e=g_xQ_eg_y^{-1} \;\;  \text{for } e=(x,y),\qquad
	(g\cdot \Phi)_x=g_x\Phi_x.
\end{align*}

Let
$m$ be the Haar measure  on $\cQ$ and $\lambda$ be the Haar measure on $\sG$.
By definitions, $m$ is invariant under the transformation $A=(Q,\Phi)\to g\cdot A=(g\cdot Q,g\cdot \Phi)$.
Define a gauge $v:\cQ\to \sG$ such that for $A=(Q,\Phi)$
$$
v(A)^{-1}_x\Phi_x=I_N.
$$
We then have
$$v(A)_x=\Phi_x.$$
According to \cite[Definition~B.2]{MR3982691},
the map $v$ should be smooth which is obvious here,
and also
$$
v(g\cdot A)=g\cdot v(A)
$$
which is clear since both sides are equal to $g\cdot \Phi$. Thus $v$ is a well-defined gauge.

We then have
\begin{equ}[e:def-Av]
\cQ_v\eqdef\{v(A)^{-1} \cdot A\,:\, A\in \cQ\}
\cong
G^{E_{\Lambda}^+}\times \{\text{Id}\}
\cong
G^{E_{\Lambda}^+}.
\end{equ}
 Here $\text{Id}\in G^{\Lambda}$ is the Higgs field which
equals $I_N$ at every point $x\in \Lambda$.
\cite[Lemma~B.3]{MR3982691} tells us that $\cQ_v$ is an embedded submanifold of $\cQ$,
and that $A,B\in \cQ$ are in the same $\sG$-orbit if and only if
$v(A)^{-1} \cdot A = v(B)^{-1} \cdot B$ i.e. they project to the same element in $\cQ_v$.
The first identification $\cong$ above
is obviously a bijection between two sets.
Since all the maps are smooth, it is also a diffeomorphism between two manifolds.

By \cite[Prop. B.9]{MR3982691} (disintegration formula), there exists a unique measure $m_v$ on $\cQ_v$ such that
for all $f:\cQ\to[0,\infty]$ measurable,
\begin{align}\label{e:mv}
	\int_{\cQ} f(A)\dif m(A)=\int_{\cQ_v}\dif m_v(B)\int_{\sG}\dif \lambda(g)f(g\cdot B).
\end{align}

\bl
In the above setting, $m_v$ is the pushforward measure of the Haar measure $m_0$ on $G^{E_\Lambda^+}$ via \eqref{e:def-Av}.
\el

\begin{proof}
The measure $m_v$ can be expressed as $\dif m_v(A)=j(Q)\dif m_0(Q) $ for some smooth density $j:G^{E_{\Lambda}^+}\to(0,\infty)$. 
Using \eqref{e:mv} we have
	\begin{align*}
		\int_{\cQ} f(A)\dif m(A)=&\int_{\cQ_v}\dif m_v(B)\int_{\sG}\dif \lambda(g)f(g \cdot B)
		\\=&\int_{\sG}\dif \lambda(g)\int_{G^{E_{\Lambda}^+}}j(Q)\dif m_0(Q)f(g \cdot  Q,g\cdot \mathrm{Id}).
	\end{align*}
	For any $C\in \cQ_v$ we apply above with $f$ replaced by $f(C\cdot )$ with $(C\cdot Q)_e:=C_e Q_e$ and $(C\cdot \Phi)_x:=\Phi_x$ and use
	for fixed $e=(x,y)$,
	$$(C\cdot (g \cdot  Q))_e=C_eg_xQ_eg_y^{-1}=g_xg_x^{-1}C_eg_xQ_eg_y^{-1}$$ to obtain
	\begin{align*}
		\int_{\cQ} f(C\cdot A)\dif m(A)=&\int_{\sG}\dif \lambda(g)\int_{G^{E_{\Lambda}^+}}\dif m_0(Q)j(Q)f(C\cdot (g \cdot Q),g\cdot\mathrm{Id})
		\\=&\int_{\sG}\dif \lambda(g)\int_{G^{E_{\Lambda}^+}}\dif m_0(Q)j(\tilde{Q}_g^C)f (g \cdot Q ,g \cdot \mathrm{Id} ),
	\end{align*}
	with $(\tilde Q_g^{C})_e=g_x^{-1}C_e^{-1}g_x Q_e$.

	On the other hand, since $m$ is Haar measure when restricted to $G^{E_{\Lambda}^+}$  we have
	\begin{align*}
		\int_{\cQ} f(C\cdot A)\dif m(A)=&\int_{\cQ} f(A)\dif m(A)
		\\=&\int_{\sG}\dif \lambda(g)\int_{G^{E_{\Lambda}^+}}j(Q)\dif m_0(Q)f(g \cdot Q,g \cdot \mathrm{Id}).
	\end{align*}
	Since $f$ is arbitary measurable function on $\cQ$, we have $j(\tilde Q_g^C)=j(Q)$ for any $C\in \cQ_v$ and $g\in \sG$. Hence, choosing $C_e=Q_e, g=\text{Id}$, we find that $j$ is a constant.
	Hence, the result follows.
\end{proof}

 Define a new action functional on $G^{E_\Lambda^+}$
  \begin{equ}
 	\mathcal S_{\YMH}'(Q) \eqdef
 	N\beta \, \Re \sum_{p\in \CP^+_\Lambda} \Tr(Q_p)
 	+2\kappa N \sum_{e\in E^+_\Lambda }\tr(Q_e).
 \end{equ}
and a new probability measure on $G^{E_\Lambda^+}$
\begin{equ}[e:gauge-fix-nu]
\widetilde\nu_\Lambda
=\frac1{Z_\Lambda} \exp\Big(\cS_{\YMH}'(Q)\Big)\prod_{e\in E_\Lambda^+}\dif \sigma_N(Q_e),
\end{equ}
where $\sigma_N$ is the Haar measure on $G$, $Z_\Lambda$ is a normalization constant. The following lemma allows us to analyze gauge-invariant observables under the measure $\widetilde\nu_{\Lambda}$, which depends only on the Yang--Mills component and consequently exhibits a simpler structure than $\mu_\Lambda$.

\bl\label{fo:gauge}
For any gauge invariant observable $h$ on $\cQ$, it holds that
\begin{align*}
	\E_{\mu_\Lambda}(h)=\E_{\widetilde\nu_\Lambda}(h(\cdot,\mathrm{Id})).
\end{align*}
Here $\widetilde\nu_\Lambda$ is given by \eqref{e:gauge-fix-nu}, 
and   $\mu_\Lambda$ is
our original Yang--Mills--Higgs measure.
\el
\begin{proof}
We choose $f$ in \eqref{e:mv} as
$\frac{1}{Z_{\Lambda}}\exp(\mathcal S_{\YMH} )\times h$.
Using \eqref{e:mv} we have
\begin{align*}
	\E_{\mu_\Lambda}(h)=	&\int_{\cQ} \frac1{Z_{\Lambda}}\exp\Big(\mathcal S_{\YMH}(Q,\Phi)\Big)h(Q,\Phi)\dif m(Q,\Phi)
	\\=&\int_{G^{E_\Lambda^+}}\dif m_0(Q)\int_{\sG}\dif \lambda(g)\frac1{Z_{\Lambda}}\exp\Big(\mathcal S_{\YMH}(g \cdot Q,g\cdot\mathrm{Id})\Big)h(g \cdot Q, g \cdot \mathrm{Id} )
	\\=&\int_{G^{E_\Lambda^+}}\dif m_0(Q)\int_{\sG}\dif \lambda(g)\frac1{Z_{\Lambda}}\exp\Big(\mathcal S_{\YMH}(Q,\mathrm{Id})\Big)h(Q, \mathrm{Id})
	\\=&\int_{G^{E_\Lambda^+}}\dif m_0(Q)\frac1{Z_{\Lambda}}\exp\Big(\mathcal S_{\YMH}'(Q)\Big)h(Q,\mathrm{Id})=\E_{\widetilde\nu_\Lambda}(h(\cdot,\mathrm{Id})),
\end{align*}
where we used gauge invariance of $\cS_{\YMH}$ and $h$ in the third step and $S_{\YMH}(Q,\mathrm{Id})=\mathcal S_{\YMH}'(Q)$ in the last step.
\end{proof}


We also know that $\widetilde\nu_\Lambda$ is a Gibbs  measure on $G^{E_\Lambda^+}$ and we
can consider the invariant Markov dynamics for $\nu_\Lambda$ and calculate $\Hess_{\cS'_{\YMH}}$.
We also define the following constant $K_{G}^{\widetilde\nu}$ which will play the role of the lower bound of the Bakry--\'Emery condition \eqref{con:B-E} for  the measure $\widetilde\nu_\Lambda$ with the action functional $\cS'_{\YMH}$:  
\begin{equ}[e:KS-G1]
 K_{G}^{\widetilde\nu}
 \eqdef \frac14(N-2)-N\Big(8(d-1)|\beta|+2|\kappa|\Big).	\tag{B\'E-G${}^\prime$}
\end{equ}

\bl\label{lem:lognu}
Suppose that $K_{G}^{\widetilde\nu}>0$,
then
the log-Sobolev and the Poincar\'e inequalities hold for $\widetilde\nu_\Lambda$,
i.e. for $F\in C^\infty(G^{E_\Lambda^+})$
\begin{align*}
	\text{ent}_{\widetilde\nu_{\Lambda}}(F^2)\leq \frac2{K_{G}^{\widetilde\nu}}\sum_e\E_{\nu_\Lambda}(|\nabla_eF|^2),\quad
		\text{var}_{\widetilde\nu_\Lambda}(F)\leq \frac1{K_{G}^{\widetilde\nu}}\sum_e\E_{\nu_\Lambda}(|\nabla_eF|^2).
\end{align*}
\el
\begin{proof}
	We have
	\begin{align*}
		\cS'_{\YMH}(Q)=	N\beta \, \Re \sum_{p\in \CP^+_\Lambda} \Tr(Q_p)
		- 2\kappa N\sum_{e\in E^+_\Lambda }(-\tr(Q_e))
		:=\cS_1'-\cS_2'.
	\end{align*}
	Note that the state space is $G^{E^+_\Lambda}$. Write $v=(v_e)_{e\in E^+}$,  and $|v|^2=\sum_e |v_e|^2$.
	As in \eqref{e:YMHess}, one has
	\begin{equ}
		|\Hess_{\cS_1'}(v,v)|\leq 8(d-1)N|\beta||v|^2.
	\end{equ}
Also,
	\begin{equation}\label{eq:Hs-f}
		\aligned
		|\Hess_{\cS_2'}(v,v)|
		&=|2\kappa N\sum_e\tr(X_e^2Q_e)|\leq 2|\kappa|N\sum_e|X_e|^2.
		\endaligned
	\end{equation}
	The Ricci curvature is the same as in \eqref{e:RicciG}.
Since $K_{G}^{\widetilde\nu}>0$,	we have verified the Bakry--\'Emery condition \eqref{con:B-E}.
	Hence, the result follows.
\end{proof}

	\bc
	\label{cor:fix-gauge}
	Let $M=G$.
	Suppose that $K_{G}^{\widetilde\nu}>0$. 
	For gauge invariant observables $F, H\in C^\infty_{cyl}(\cQ)$, suppose that $\Lambda_F\cap \Lambda_H=\emptyset$.
Then one has for any infinite volume limit $\mu$ of $\{\mu_\Lambda\}$
\begin{align*}
	|\cov_\mu(F,H)|
	\leq c_{1} d(\mfg)e^{-c_N d(\Lambda_F,\Lambda_H)}
	(\$ F\$_{\infty}\$ H \$_{\infty}+\|F\|_{L^2}\|H\|_{L^2}),
\end{align*}
where $c_{1}$ depends on $|\Lambda_F|$, $|\Lambda_H|$,
and $c_N$ depends on $K_{G}^{\widetilde\nu}$,  $N$ and $d$.
\ec
\begin{proof}
Since $\widetilde\nu_\Lambda$ is a Gibbs measure on the compact space $G^{E^+_\Lambda}$, 
using Lemma \ref{lem:lognu}  and by the same arguments as in the proof of Corollary~\ref{c:ex} or \cite[Corollary 4.11]{SZZ22} we get
\begin{align}\label{eq:mass}
	|\cov_{\widetilde\nu_\Lambda}(F,H)|
	\leq c_{1} d(\mfg)e^{-c_N d(\Lambda_F,\Lambda_H)}(\$ F\$_{\infty}\$ H \$_{\infty}+\|F\|_{L^2}\|H\|_{L^2}),
\end{align}
 with constants $c_1, c_N$ independent of $\Lambda$. 
 Using  Lemma \ref{fo:gauge}  we obtain that for   gauge invariant observables $F,H$, \eqref{eq:mass} holds with $\widetilde\nu_\Lambda$ replaced by $\mu_\Lambda$. Letting $\Lambda\to\mZ^d$ the result follows.
\end{proof}

\br \label{rem:CompareGG}
(i) Comparing \eqref{e:KS-G} and \eqref{e:KS-G1},
we note that $K_G< K_{G}^{\widetilde\nu}$.
This is 
because  $K_{G}^{\widetilde\nu}$ is equal to the first term in \eqref{e:KS-G} with $\delta=0$,
and the two terms in \eqref{e:KS-G} are strictly increasing and decreasing respectively in $\delta$.
This can be also seen by observing that 
$K_{G}^{\widetilde\nu}>0$ is equivalent with 
$K_{\YM}>0$ and $|\kappa|<K_{\YM}/(2N)$ as in the last statement of Theorem~\ref{main4},
and the latter condition is just \eqref{e:small-kappa} with $\alpha=0$. 

In other words,
the condition $K_{G}^{\widetilde\nu}>0$ for mass gap is weaker than $K_{G}>0$ in Corollary \ref{c:ex}. 

(ii) As in Corollary \ref{cor:SG-erg} we also have the uniqueness of the infinite volume limit of $\widetilde\nu_\Lambda$. This  however does not imply the uniqueness of the infinite volume limit of $\mu_\Lambda$.
\er

\section{Large $N$ limit for  $M\in \{ \mS^{N-1}, G\}$}\label{sec:large} 

In this section we discuss the large $N$ limit of gauge invariant observables by applying the Poincar\'e inequality. For pure Yang--Mills, this was discussed in \cite{SZZ22}.
Recall the gauge invariant observables Wilson loops and Wilson lines defined in \eqref{de:wlo} and \eqref{de:wl}.


\bt[Large $N$ limit of gauge invariant observables]
\label{co:1} For $M=\mS^{N-1}$, assume $K_{\mS^{N-1}}>0$.
For every Wilson loop variable $W_l$  and the unique infinite volume limit $\mu$ from Corollary \ref{cor:SG-erg},
  one has
\begin{equ}[e:varW]
	\var_\mu\Big(\frac1NW_l\Big)\leq \frac{n(n-3)}{K_{\mS^{N-1}} N} \;.
\end{equ}
In particular, under the above assumptions, one has the convergence
\begin{equ}[e:EW]
	\Big|\frac{W_l}{N}-\E\frac{W_l}{N}\Big|\to0
	\qquad
	\mbox{as } N\to \infty
\end{equ}
in probability, and for any loops $l_1,\dots, l_m$,
  the factorization property holds
\begin{equ}[e:factorize]
	\lim_{N\to\infty}\Big|\E\frac{W_{l_1}\dots W_{l_m}}{N^m}-\prod_{i=1}^m\E\frac{W_{l_i}}{N}\Big|=0.
\end{equ}

For $M=G$, assume $K_{G}^{\widetilde\nu}>0$. 
For every infinite volume limit $\mu$,   every Wilson loop or Wilson line variable $W_l$, one has
\begin{equ}[e:varWG]
	\var_\mu\Big(\frac1NW_l\Big)\leq \frac{n(n-3)}{K_{G}^{\widetilde\nu} N} \;,
\end{equ}
and \eqref{e:EW} \eqref{e:factorize} hold for both Wilson loops and Wilson lines.
\et

\begin{proof}
When $M= \mS^{N-1}$,	
\eqref{e:varW} follows by exactly the same argument as in  \cite[Corollary 1.5]{SZZ22}, namely, we apply the Poincar\'e inequality for $\mu$ obtained in Corollary \ref{cor:SG-erg}, with $F=\frac1NW_l$.
	
	Concerning the case of Wilson loops and Wilson lines $W_l$ for $M= G$, we  first perform the gauge fixing as in Section \ref{sec:gauge fixing}  to rewrite $W_l$ into a variable $W_l'$, which only depends on the field $Q$ (since the field $\Phi$ is turned into identity via this gauge fixing). Then the variable $W_l'$ can be viewed as a Wilson loop variable, except that the two ends of the lattice curve $l$ may not coincide.
Recalling the gauge fixed measure $\nu_\Lambda$ defined in \eqref{e:gauge-fix-nu},
we then apply the Poincar\'e inequality for $\nu_\Lambda$ from Lemma \ref{lem:lognu} with $F=\frac1NW_l'$ as in the proof of \cite[Corollary 1.5]{SZZ22} (which does not depend on whether the curve $l$ is a loop or not)
to have
\begin{equ}
	\var_{\mu_\Lambda}\Big(\frac1NW_l\Big)=\var_{\widetilde\nu_\Lambda}\Big(\frac1NW_l'\Big)\leq \frac{n(n-3)}{K_{G}^{\widetilde \nu} N} \;,
\end{equ}
where the first identity is by Lemma \ref{fo:gauge}.
Letting $\Lambda\to\mZ^d$, we obtain \eqref{e:varWG}.
The convergence to the mean follows from Chebychev's inequality.  For the factorization we use $W_l\leq N$ and H\"older's inequality to have
		\begin{align*}
		\Big|\E\frac{W_{l_1}\dots W_{l_m}}{N^m}-\prod_{i=1}^m\E\frac{W_{l_i}}{N}\Big|\leq 	\sum_{i=1}^m\var_{\mu}\Big(\frac1NW_{l_i}\Big)^{1/2},
		\end{align*}
		which implies the result.
\end{proof}

\br It is natural to ask whether the factorization property holds for Wilson lines when $\Phi_x\in \mS^{N-1}$. In this case, for a given Wilson line $W_l$, we note that $W_l$ is of order $1$. We then take $F=W_l$ in the Poincar\'e inequality and find  by direct calculation that $\cE_\nu(F,F)$ is of order $N$.  Note that $K_{\mS^{N-1}}$ is also of order $N$. We would not expect factorization to hold in this case.
\er

\section*{Declarations}

The authors have no competing interests to declare that are relevant to the content of this article.

Data sharing not applicable to this article as no datasets were generated or analyzed during the current study.

\bibliographystyle{alphaabbr}
\bibliography{refs}

\end{document}